\documentclass[a4paper]{amsart}

\title[Dirichlet Stereohedra for Full Cubic Groups]
{On the Number of Facets of Three-Dimensional Dirichlet Stereohedra III: Full Cubic Groups}
\thanks{Research partially supported by  the Spanish Ministry of Education and Science, grant number
 MTM2005-08618-C02-02.}

\author{Pilar Sabariego
        \and
        Francisco Santos}
        \address{Departamento de Matem\'aticas, Estad\'{\i}stica y
        Computaci\'on, Universidad de Cantabria, Santander, Spain}
        \email{pilar.sabariego@unican.es, francisco.santos@unican.es.}

\usepackage{latexsym}
\usepackage{amssymb}
\usepackage{amsmath}
\usepackage[dvips]{graphicx}

\usepackage{color}

\usepackage{multirow}

\newtheorem{theorem}{Theorem}[section]
\newtheorem{cor}[theorem]{Corollary}
\newtheorem{proposition}[theorem]{Proposition}
\newtheorem{corollary}[theorem]{Corollary}
\newtheorem{lemma}[theorem]{Lemma}
\newtheorem{remark}[theorem]{Remark}

\def\Z{\mathbb Z}
\def\R{\mathbb R}

\newcommand{\Isom}{\operatorname{Isom}}
\newcommand{\Vor}{\operatorname{Vor}}
\newcommand{\VorExt}{\operatorname{VorExt}}
\newcommand{\Infl}{\operatorname{Infl}}
\newcommand{\conv}{\operatorname{conv}}
\newcommand{\nor}{\operatorname{\mathcal N}}

\def\figura(#1,#2,#3,#4,#5){%
\leavevmode \hbox to #4{\vrule height#5 width0mm \special{#1 #2
scaled #3}\hfill}}

\begin{document}

\begin{abstract}
We are interested in the maximum possible number of facets that Dirichlet stereohedra
for  three-dimensional crystallographic groups can have. 
The problem for non-cubic groups was
studied in two previous paper by D.~Bochis and  the second author. This paper deals with 
``full'' cubic groups, while ``quarter'' cubic groups are left for a subsequent paper.
Here,``full''  and ``quarter''  refers to the recent classification of three-dimensional
crystallographic groups by Conway, Delgado-Friedrichs, Huson
and Thurston.

Our main result in this paper is  that Dirichlet stereohedra for any of the 27 full 
groups  cannot have more than $25$ facets. We also find
stereohedra with 17 facets for one of these groups.
\end{abstract}

\maketitle

\section{Introduction}

Tiling problems are among the most classical problems in discrete geometry,
yet they include many basic ones that are still open. See~\cite[Chapter 4]{Brass-etal-2005}
for a recent survey. A fundamental question, explicitly asked in~\cite[p. 177]{Brass-etal-2005}, \cite[p.~959]{GruShe-art}, and \cite{SchSen-1997}, and related to Hilbert's 18th problem \cite{Milnor-1974}
is:
\begin{quote}
What is the maximum number of facets of a convex polytope $P$ such that
$\R^3$ can be tiled by congruent copies of  $P$?
\end{quote}

The fact is that there is no known upper bound for the number of facets of such a polytope.
Only in the special
case where the symmetry group of the tiling acts transitively over the tiles, a global
upper bound is known. In this case, the tiling is called \emph{stereohedral} and the
tile is a \emph{stereohedron}:

\begin{theorem}[Delone~\cite{Delone-1961}]
A stereohedron in Euclidean $d$-space with respect to a symmetry group $G$
cannot have more
than $2^d (a + 1) -2$ facets, where $a$ is the number of \emph{aspects} of $G$.
\end{theorem}

The group $G\subset\Isom(\R^d)$ that acts transitively on a stereohedral tiling is necessarily  a \emph{crystallographic group}.
That is, it contains a full-dimensional discrete lattice $T\cong\Z^d$ as a subgroup of finite index.
The quotient group $G/T$ is usually called the \emph{point group} of $G$, and its cardinality
is the  number of \emph{aspects} of $G$. Equivalently, the number of aspects of $G$  
equals the number of translational
lattices in which a generic orbit of $G$ decomposes.
In dimension three, crystallographic groups can have up to 
48 aspects. Hence, the above ``fundamental theorem of the theory of stereohedra'' implies that three-dimensional stereohedra cannot have more than 390 facets.

On the other hand,  the stereohedron
  with the maximum number of facets known so far has ``only'' 38 of them.
It was found by P.~Engel in 1980 (see \cite{Engel81a} and \cite[p. 964]{GruShe-art}),
 using a computer search, and it tiles space under the action of the cubic group $I4_132$, with 24 aspects.


It is also interesting to observe that ``unbounded stereohedra'', that is, 
unbounded polyhedra that tile $\R^d$ by congruent copies under the action of a discrete 
transformation group, can have arbitrarily many facets, even in dimension 3. 
See~\cite{Erickson-Kim}, or  Lemma~4.2 in ~\cite{Bochis-Santos-2006}.

\medskip

A special case of stereohedral tiling is the Voronoi diagram of
an orbit of points in a crystallographic group. The stereohedra obtained
are called \emph{Dirichlet stereohedra} or \emph{plesiohedra}. 
The fact that Engel's 38 face stereohedron is
actually a Dirichlet stereohedron seems to indicate that Dirichlet stereohedra can
have as many facets as general stereohedra (although there is  no proof of this).
Since Dirichlet stereohedra are much easier to analyze, 
in previous papers the second author together with
D.~Bochi\c{s} has initiated an exhaustive study of the number of
facets of Dirichlet stereohedra for the different three-dimensional crystallographic groups.
They proved:
 \begin{itemize}
 \item Within the 100 crystallographic groups which contain
 reflection planes, the exact maximum number of facets of
 Dirichlet stereohedra
 is 18. Within \emph{cubic} groups containing reflection planes
 the maximum number of facets is eight~\cite{Bochis-Santos-2001}.
  \item Within the 97 non-cubic crystallographic groups without
 reflection planes,  no Dirichlet stereohedron can have more than 80 
 facets
 ~\cite{Bochis-Santos-2006}.
 The same paper shows an explicit construction of Dirichlet stereohedra with 32
 facets for a tetragonal group.
 
 \end{itemize}
 
 For cubic groups without reflection planes (there are 22 of them), Bochi\c{s}
 and Santos were only able to prove an upper bound of 162 facets~\cite{Bochis-1999}. 
 Here and in~\cite{SabariegoSantos-quarter}
 we improve this bound, and hence the general 
 upper bound for the number of facets of three-dimensional Dirichlet stereohedra. 
For technical reasons, we do this separately for ``full groups'' (this paper) and ``quarter groups'' 
(paper~\cite{SabariegoSantos-quarter}). Among other things, 
the proof we present here for full groups is ``computer-free'', while the
proof for quarter groups is not.
 
 \medskip
 
   Our work takes great advantage of  the new classification of crystallographic groups developed 
 by Conway et at.~in~\cite{Conway-etal-2001}. 
 To enumerate cubic groups, they use the fact that every cubic group contains
 rotations of order three in the four diagonal directions of the orthonormal lattice $\Z^3$.
 They call \emph{odd subgroup} of a cubic crystallographic group $G$ the subgroup 
 generated by all these rotations and  show that:
  
 \begin{theorem}[Conway et al.~ \cite{Conway-etal-2001}]
 \label{thm:conway-etal}
For every cubic crystallographic group $G$:
 \begin{enumerate}
 \item Its odd subgroup is one of two possible groups $R_1$ or $R_2$.
 \item $R_1$ and $R_2$ are normal in \ $\Isom(\R^3)$, hence $G$ lies between $R_i$ and its normalizer $\nor(R_i)$, for
 the appropriate $i$.
 \end{enumerate}
 \end{theorem}
 
 In particular, the classification of cubic groups is equal to the classification of (conjugacy classes of)
 subgroups of the finite groups $\nor(R_1)/R_1$ and $\nor(R_2)/R_2$, which are of orders
 16 and 8. (The latter is isomorphic to the dihedral group $D_8$; the former contains $D_8$ as an index two subgroup). The odd subgroups $R_1$ and $R_2$
 differ in that $R_1$ contains rotation axes of order three that intersect one another while
 all rotation axes of order three in $R_2$ are mutually disjoint. 
 Conway et al.~call
 \emph{full groups} those with odd subgroup equal to $R_1$ and \emph{quarter} groups those
 with odd subgroup equal to $R_2$. (The reason is that $R_2$ is a subgroup of $R_1$ of index four.
 That is, quarter groups contain only one quarter of the possible rotation axes).
 There are 27 full groups (14 of them without reflection planes)
 and eight quarter groups (none with reflection planes, since $\nor(R_2)$ does not contain reflections).
 
 The main result in this paper is:
 
 \begin{theorem}
 \label{thm:main}
 Dirichlet stereohedra for full cubic groups cannot have more than
 25 facets. 
 \end{theorem}
 
 More precisely, for each of the 14 full groups without reflections we obtain the upper bound
 shown in Table~\ref{table:main}. Here and in the rest of the paper we use the
International Crystallographic Notation for crystallographic
groups, see \cite{Lockwood-1978}.
 Let us remind the reader that for cubic groups with reflections, Theorem~2.4 in \cite{Bochis-Santos-2001} gives a bound of eight facets.
 
 \begin{table}
 \[
\begin{array}{|c|c|c|c|}
\hline
\text{$|G:R_1|$} & \text{Group $G$} &  \text{Our bound} & \text{proved by} 
\\
\hline 
1 & R_1=F23\ ^* & 10 & \text{Cor.~\ref{coro:firstbound}}
\\
\hline
\multirow{5}*{2}
& F432\ ^* & 14 &\text{Cor.~\ref{coro:firstbound}}
\\
\cline{2-4}
& F\overline{4}3c &  14 & \text{Cor.~\ref{coro:firstbound}}
\\
\cline{2-4}
& F\frac{2}{d}\overline{3} & 14 & \text{Cor.~\ref{coro:firstbound}}
\\
\cline{2-4}
& P23\ ^* &  15  & \text{Section~\ref{sec:order2}}
\\
\cline{2-4}
& F4_{1}32 & 17 & \text{Cor.~\ref{coro:firstbound}}
\\
\hline
\multirow{6}*{4} 
& P432\ ^* &  11  & \text{Cor.~\ref{coro:order4}} 
\\ 
\cline{2-4}
& I23\ ^* &  21  & \text{Section~\ref{sec:order2}}
\\
\cline{2-4}
& P\frac{2}{n}\overline{3}  & 23  & \text{Section~\ref{sec:order2}}
\\
\cline{2-4}
& F\frac{4_{1}}{d}\overline{3}\frac{2}{n} &  25& \text{Cor.~\ref{coro:firstbound}}
\\
\cline{2-4}
& P\overline{4}3n &  23 & \text{Section~\ref{sec:order2}}
\\
\cline{2-4}
& P4_{2}32 & 25  &\text{Section~\ref{sec:p4232}}
\\
\hline 
\multirow{2}*{8}
 &  P\frac{4}{n}\overline{3}\frac{2}{n} & 23 & \text{Cor.~\ref{coro:order4}} 
 \\
 \cline{2-4}
 & I432\ ^*  & 22 & \text{Cor.~\ref{coro:order4}} 
 \\
\hline
\end{array}
\]
\label{table:main}
\caption{Our bound for the number of facets of Dirichlet stereohedra of full cubic groups without reflections. The six groups marked with an asterisk are \emph{symmorphic}. Their Dirichlet stereohedra were completely described by Engel (1981), and have at most 17 facets.}
\end{table}

 Although, as far as we know, our work represents the first significant improvement to Delone's global upper bound on the number of facets of general Dirichlet stereohedra, there are at least the following previous  classifications of special Dirichlet stereohedra:
 \begin{itemize}
 \item Engel~\cite{Engel81b} completely classifies Dirichlet stereohedra for the \emph{symmorphic cubic groups}. $G$ is called symmorphic if it is the semidirect product of its translation subgroup and its point group. Equivalently, if there is an orbit of $G$ on which the translational subgroup acts transitively. 
 
 Symmorphic cubic groups are obviously all full. There are nine of them with reflections six without reflections. Engel found that their Dirichlet stereohedra have up to 17 facets.
 
 \item Koch~\cite{Koch73} completely classified Dirichlet stereohedra for cubic orbits with less than three degrees of freedom. 
 She found  117 types of stereohedra, with up to  23 facets.

 \item Dress et al.~\cite{Dressetal93} completely classified the stereohedral tilings whose symmmetry group acts transitively on the two-dimensional faces. There are 88 types of tilings, and 7 types of tiles.
 \end{itemize}
 
  \begin{remark}
 \label{rem:stabilizers}
 \rm
 In the previous papers in this series we always assumed that the base point $p$ for the orbit has trivial stabilizer for the crystallographic group $G$ under study. 
 We believed this was not a loss of generality in the sense that every crystallographic orbit with non-trivial stabilizer for a group $G$ is also an orbit with smaller stabilizer for another group $G'$.
That is actually not true, but we prove in the appendix that
there are only four groups $G$ in which it fails, all of them cubic and without reflections. These cases were previously found by  Fischer~\cite{Fischer73,Fischer80} by an exhaustive study of cubic orbits with less than three degrees of freedom. 

We also show in the appendix that in the full case these special orbits produce stereohedra with at most 12 facets. In the quarter case we refer to the classification of Koch mentioned above, which gives a bound of 23 facets.
\end{remark}

 \section{The structure of full cubic groups}
 \label{sec:full-groups}
 %
The odd group  $R_1$ of full cubic groups is a crystallographic group of type $F23$ (in
  International Crystallographic Notation).  
  It is generated
 by rotations of order three in all the lines with directions
 $(\pm1, \pm 1, \pm 1)$ and passing through integer points. In particular, the points where two or
 more axes meet form a body centered lattice 
 \[
 I=\{(0,0,0),(\frac{1}{2},\frac{1}{2},\frac{1}{2})\} + \Z.
 \] 
To understand full groups, the Delaunay tesselation of this lattice is specially helpful.
 \begin{figure} [h]
    \begin{center}
      \includegraphics[width=7cm,height=7cm]{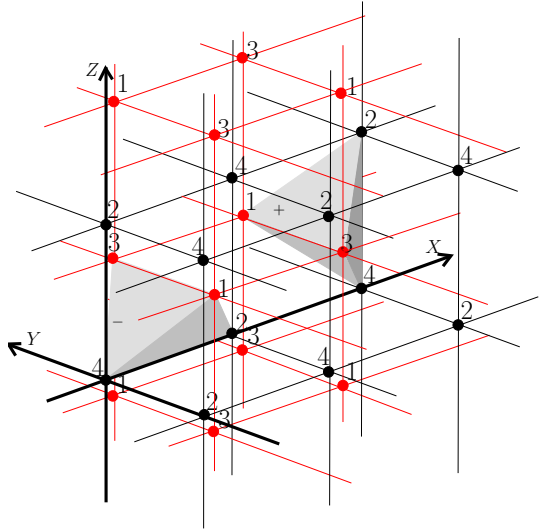}
    \end{center}
      \caption{The body-centered cubic lattice $I$. The Delaunay tetrahedra fall into two classes, according as they can be moved to coincide with the tetrahedron on the left or right, by a motion that preserves the labeling of the 
      points in $I$.
       }
     \label{fig:Red2}
\end{figure}
Each Delaunay tetrahedron has   two opposite dihedral
 angles equal to $\pi/2$,  and the other four dihedral angles equal to $\pi/3$. The latter
 lie along rotation axes of order three and the former on rotation axes of order two,
 parallel to coordinate axes. 
  Since every edge has an even number of incident Delaunay tetrahedra, 
 the Delaunay tesselation is 2-colorable.
 The group $R_1$ acts transitively and with trivial stabilizer on the Delaunay tetrahedra of a given color.

 Another way to define the coloring is as follows. The action of $R_1$ divides 
 the lattice $I$ into four orbits, namely:
\begin{align*}
F_i:=\{(a,b,c)\in I : a+b+c-\frac{i}{2} \in 2\Z\}, \qquad i=1,2,3,4.
\end{align*}
Each orbit is  a face centered cubic lattice. $F_2$ and $F_4$ together form a primitive lattice, 
and $F_1$ and $F_3$ form another.

Each Delaunay tetrahedron has a vertex in one of each of the $F_i$'s, and
two Delaunay tetrahedra have the same color if and only if the isometry that sends one to the other
and keeps the labels of vertices is orientation-preserving. 
In  Figure~\ref{fig:Red2}, lattice points are labelled according to their four orbits, and two
Delaunay tetrahedra, which belong to different ``colors'',  have been shaded.

We fix, for the rest of the paper, the following ``base Delaunay tetrahedron'' $T$:
 \[
 T:=\conv(\{(1/2,-1/2,1/2),(1,0,0), (1/2,1/2,1/2), (0,0,0)\}.
 \]
 $T$ is one of the two shaded tetrahedra in Figure~\ref{fig:Red2}. We call its four vertices $v_1$, $v_2$, $v_3$ and $v_4$, in the order we have given them, so that $v_i\in F_i$. We also introduce the following
 notation:
 \begin{itemize}
 \item Let $T_i$, $i=1,\dots,4$ be the tetrahedron adjacent to $T$ and sharing with it the triangle opposite to
 $v_i$. Let $v'_i$ the vertex of $T_i$ not in $T$.
 \item Similarly, let $T_{ij}$, with $i\ne j$, denote the neighbor of $T_i$ and sharing with it the triangle opposite
 to $v_j$. Its vertices of $T_{ij}$ are $v_k$, $v_l$, $v'_i$ and a certain vertex $v''_{ij}$, where 
 $\{i,j,k,l\}=\{1,2,3,4\}$.
   \end{itemize}
 Observe that:
 \begin{itemize}
 
 \item If $i-j$ is odd, then  $T$, $T_{ij}$ and $T_{ji}$ are related by the rotation of order three on the axis passing through $v_k$ and $v_l$.
 
 \item If $i-j$ is even, then $v''_{ij}=v'_i$. That is, 
$T_{13}=T_{31}$ and $T_{42}=T_{24}$. Moreover, 
 $T_{ij}$ is related to $T$ by an order two rotation on the axis
 passing through $v_k$ and $v_l$. The rotation of order four on this same axis relates the two of them with 
 $T_i$ and $T_j$.
 \end{itemize}

 The normalizer $\nor(R_1)$ of $R_1$ 
 is the symmetry group of the lattice $I$. It acts transitively over all the Delaunay tetrahedra and with a stabilizer of order eight, isomorphic to the dihedral group $D_8$ (this stabilizer
 is the symmetry group of the Delaunay tetrahedron itself). Hence, the index of $R_1$ in $\nor(R_1)$,
 that is, the order of the quotient group $\nor(R_1)/R_1$, is 16.
 
Knowing $\nor(R_1)$ and $R_1$ we can give the whole list of all full groups. They are pictured
in Table~\ref{table:fullgroups}. We classify them with two parameters: the order $s$ of the stabilizer
of $T$, and a number $m=0$ or $m=1$ depending on whether the group
acts on the Delaunay tesselation of the lattice $I$ with one or two orbits of tetrahedra.
 
\begin{figure}[htb]
\begin{center}
\[
\begin{array}{ccccc}
&&
\includegraphics[width=1.65cm,height=1.65cm]{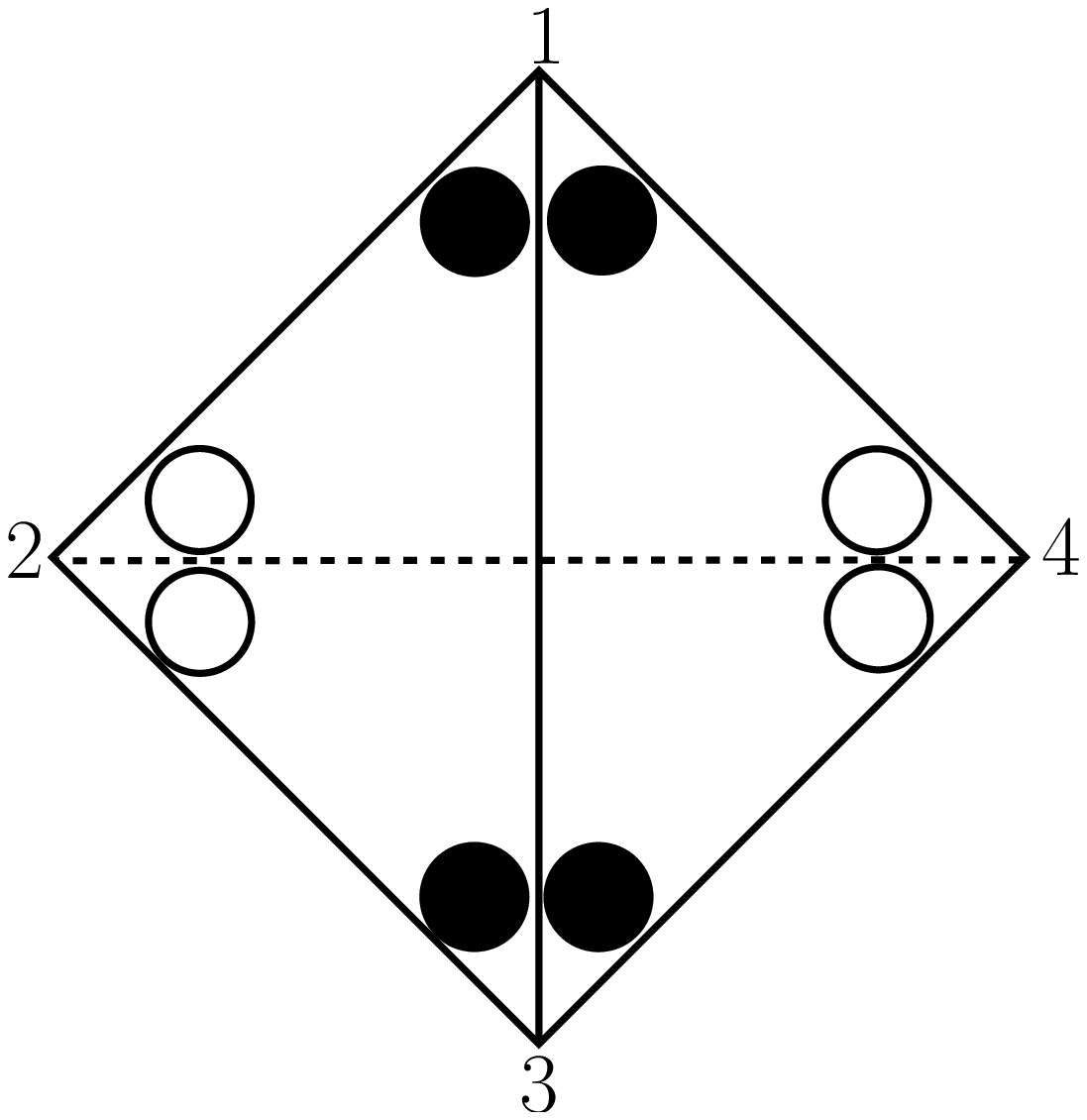}
&\\
&\swarrow&\downarrow&\searrow \medskip\\
\includegraphics[width=1.65cm,height=1.65cm]{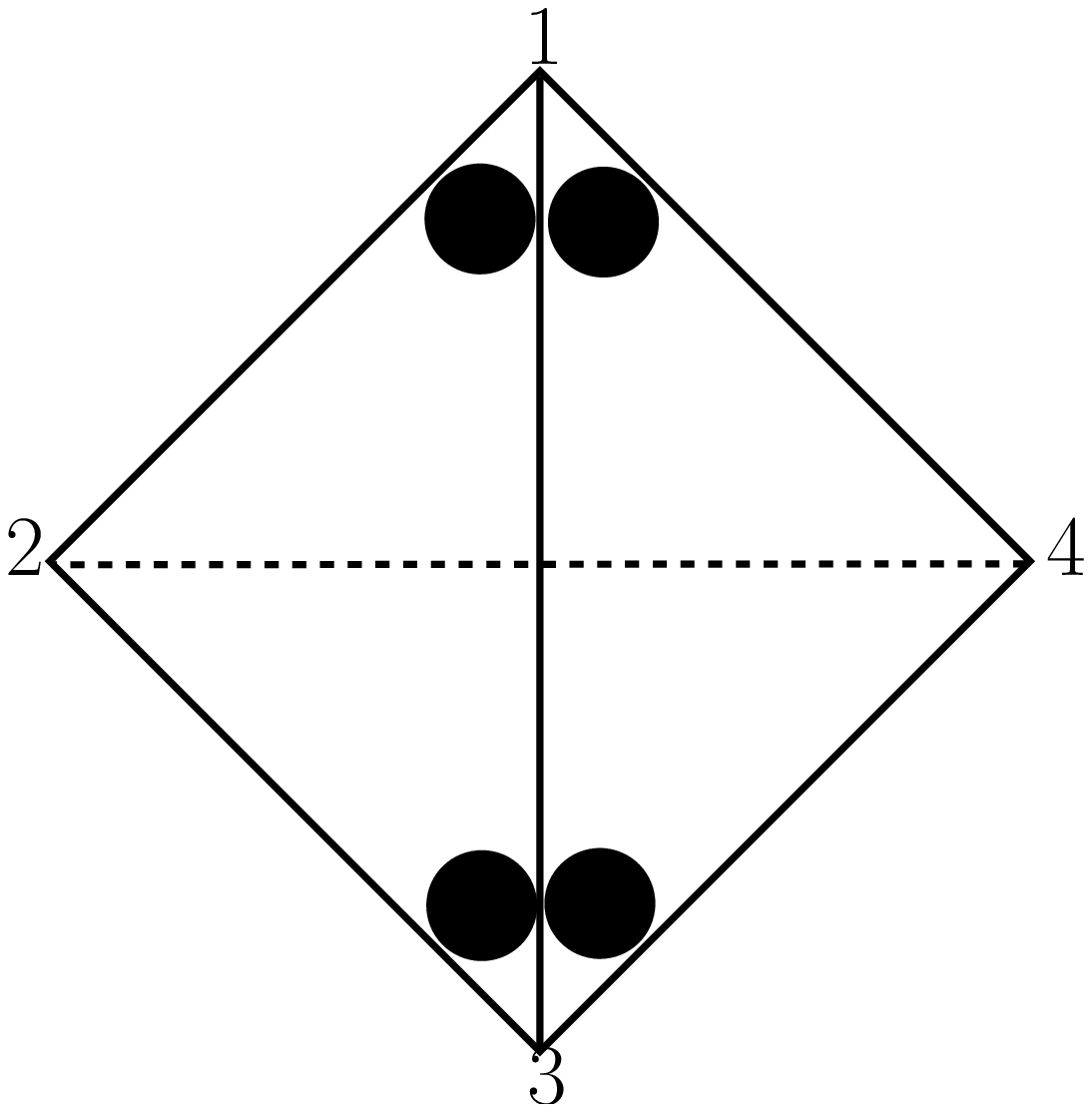}
&&
\includegraphics[width=1.65cm,height=1.65cm]{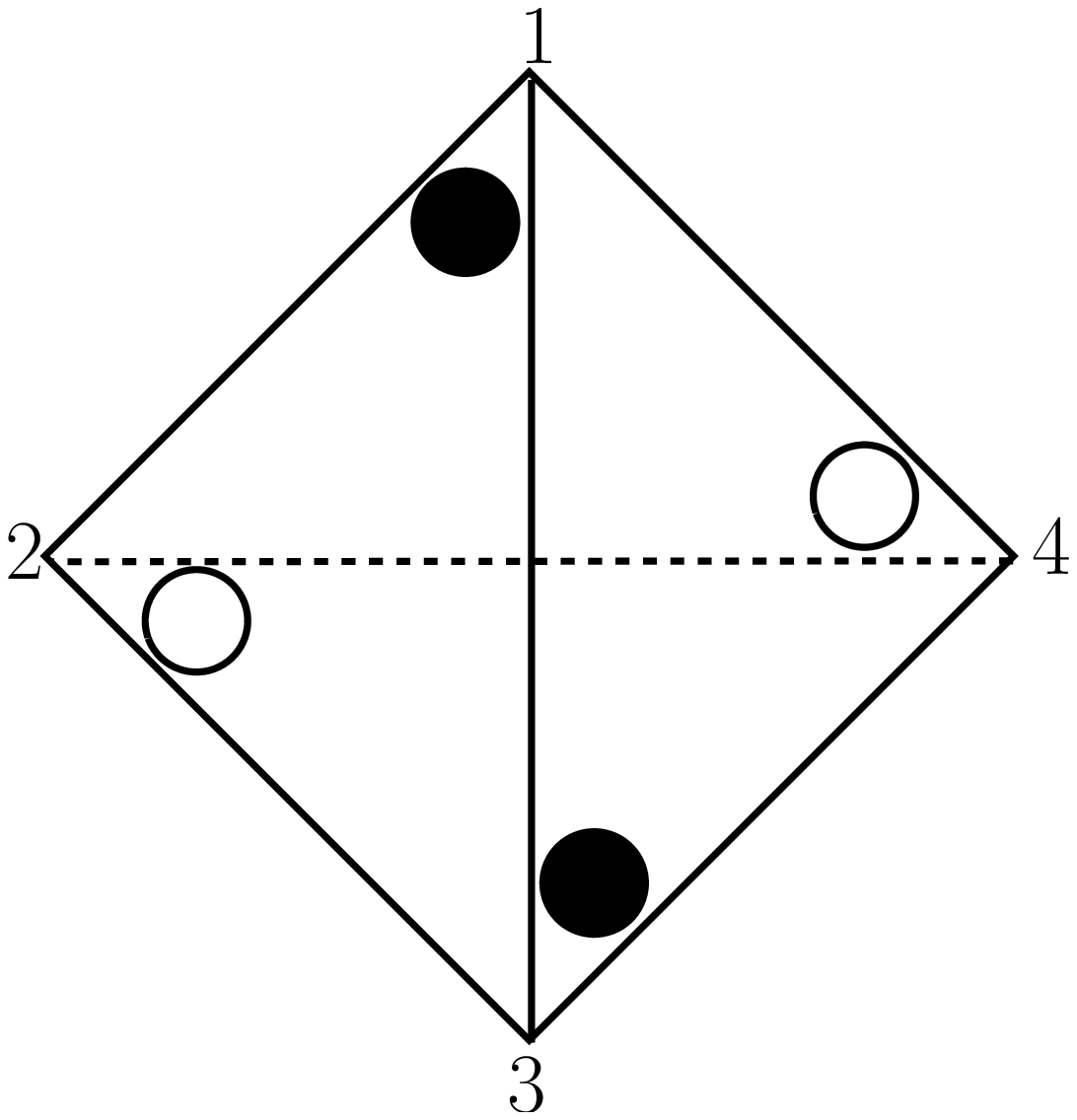}
&&
\includegraphics[width=1.65cm,height=1.65cm]{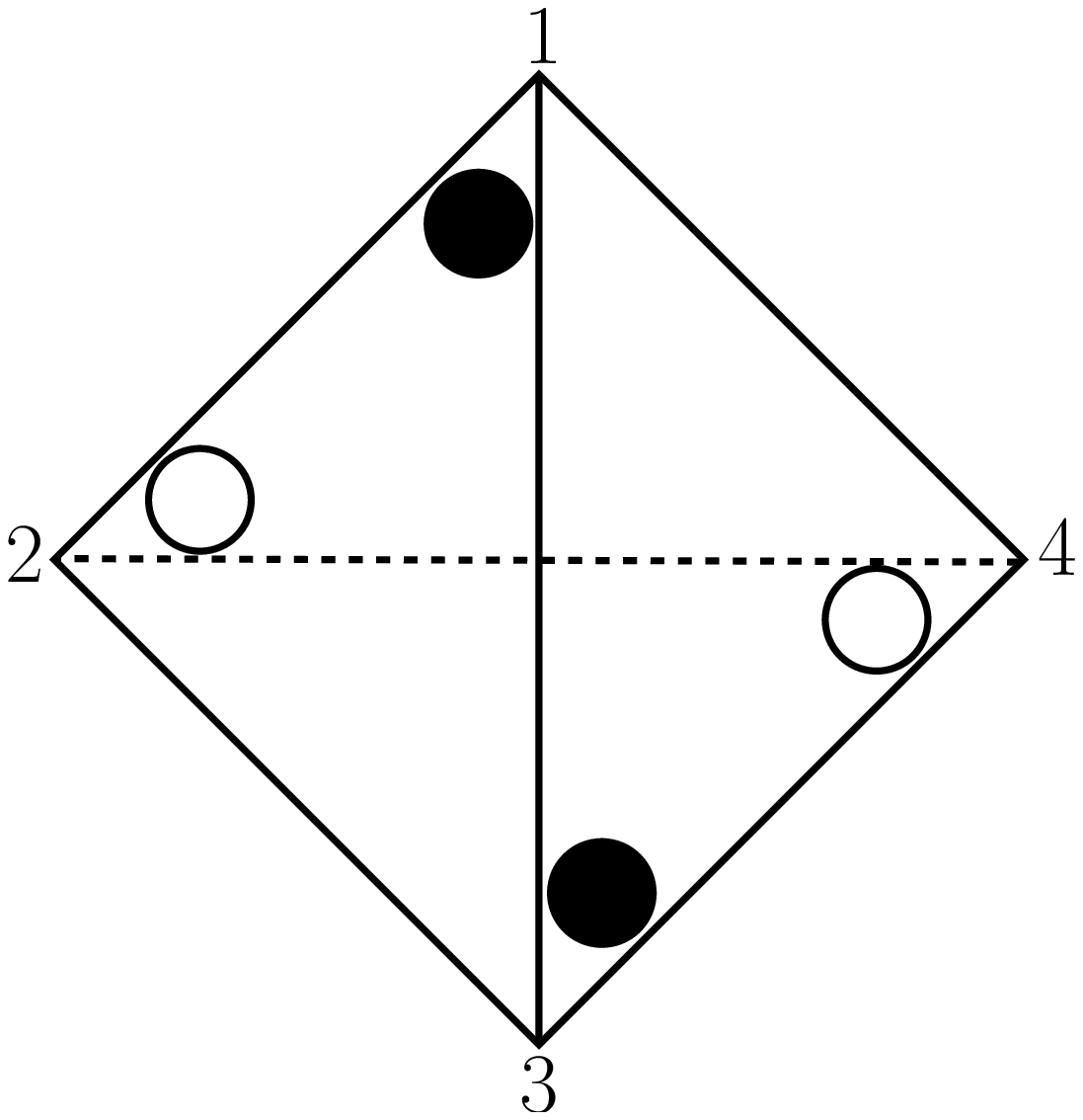}
\\
\downarrow&\searrow&\downarrow&\swarrow&\downarrow \medskip\\
\includegraphics[width=1.65cm,height=1.65cm]{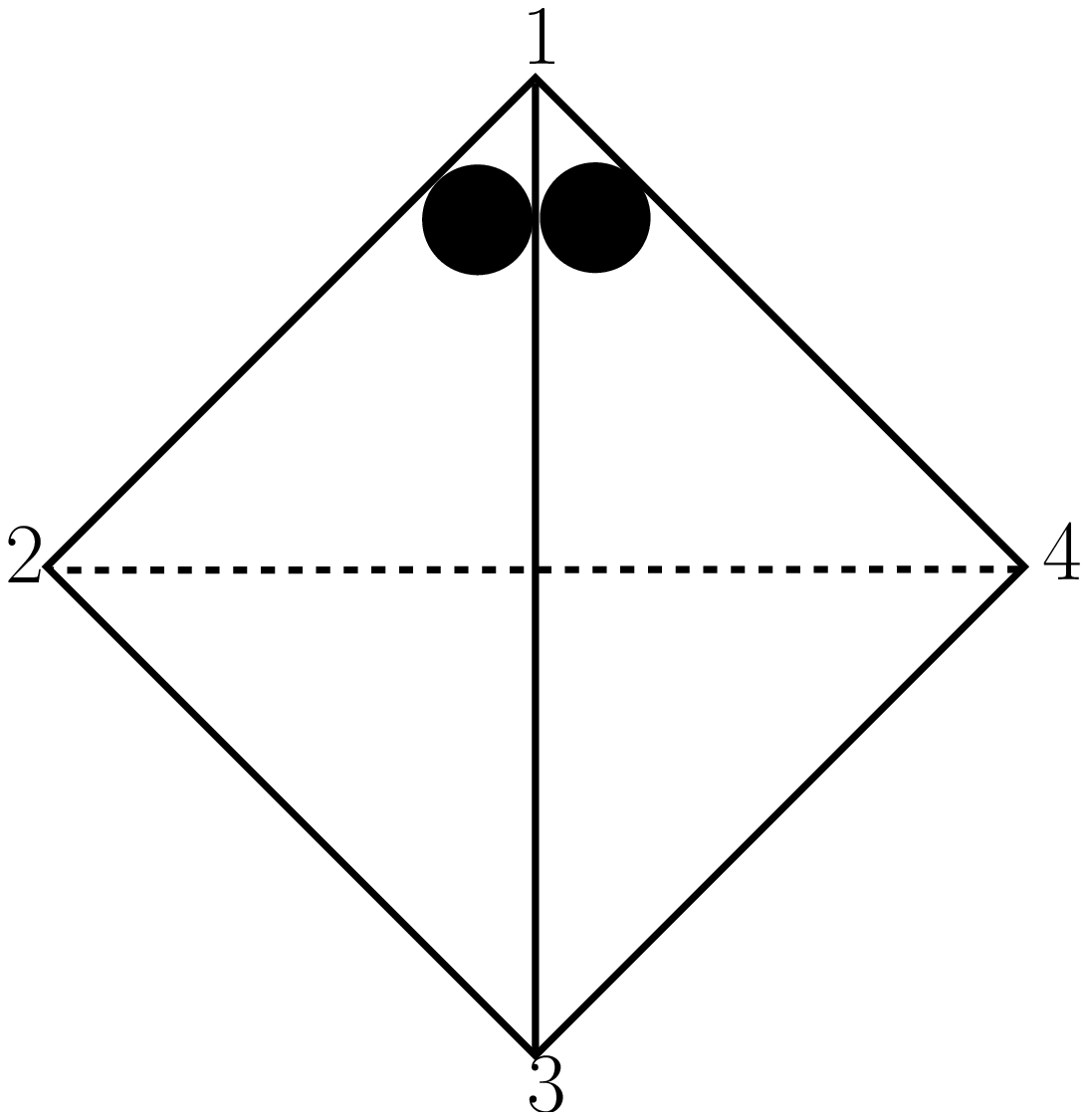}
&&
\includegraphics[width=1.65cm,height=1.65cm]{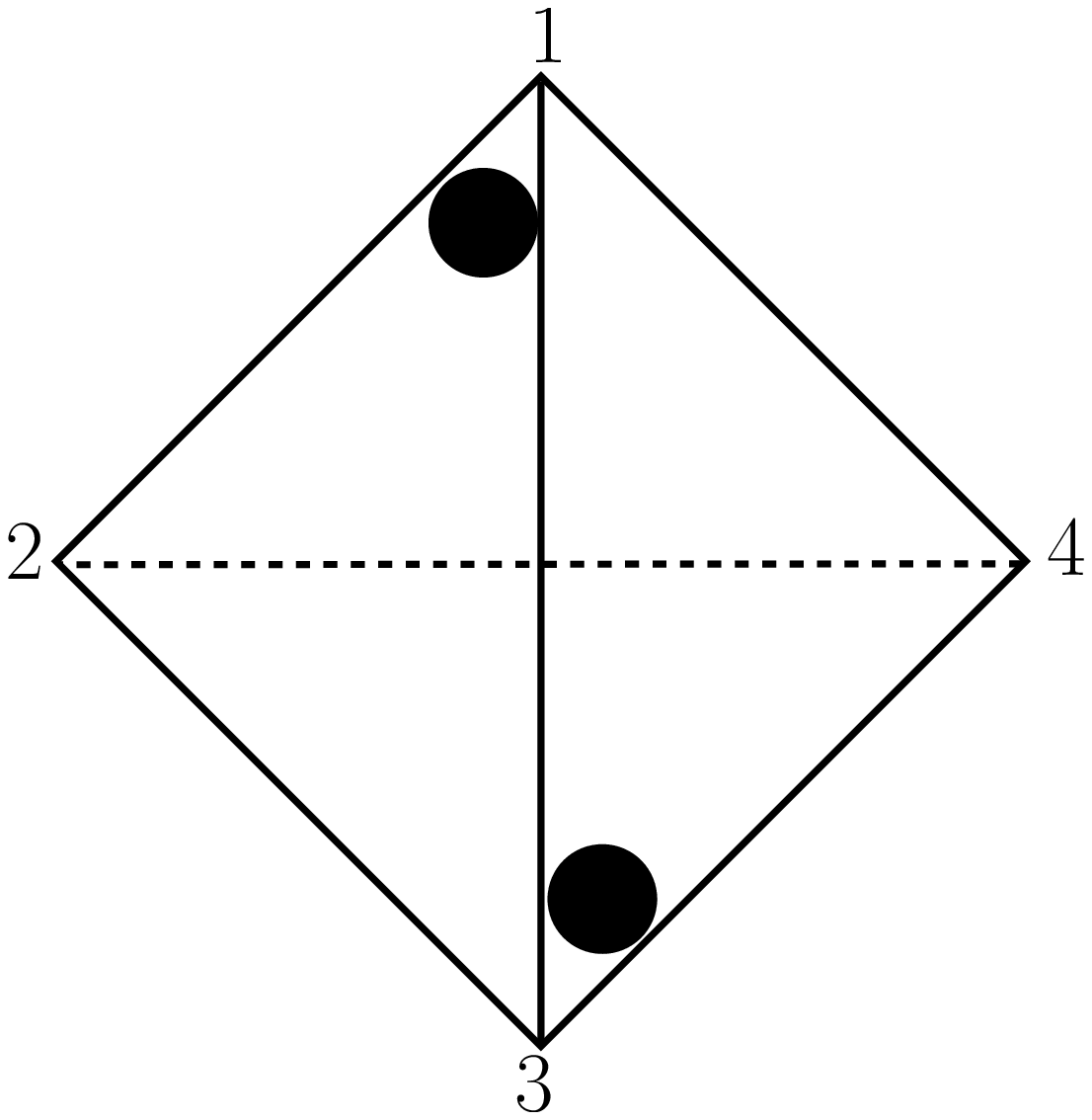}
&&
\includegraphics[width=1.65cm,height=1.65cm]{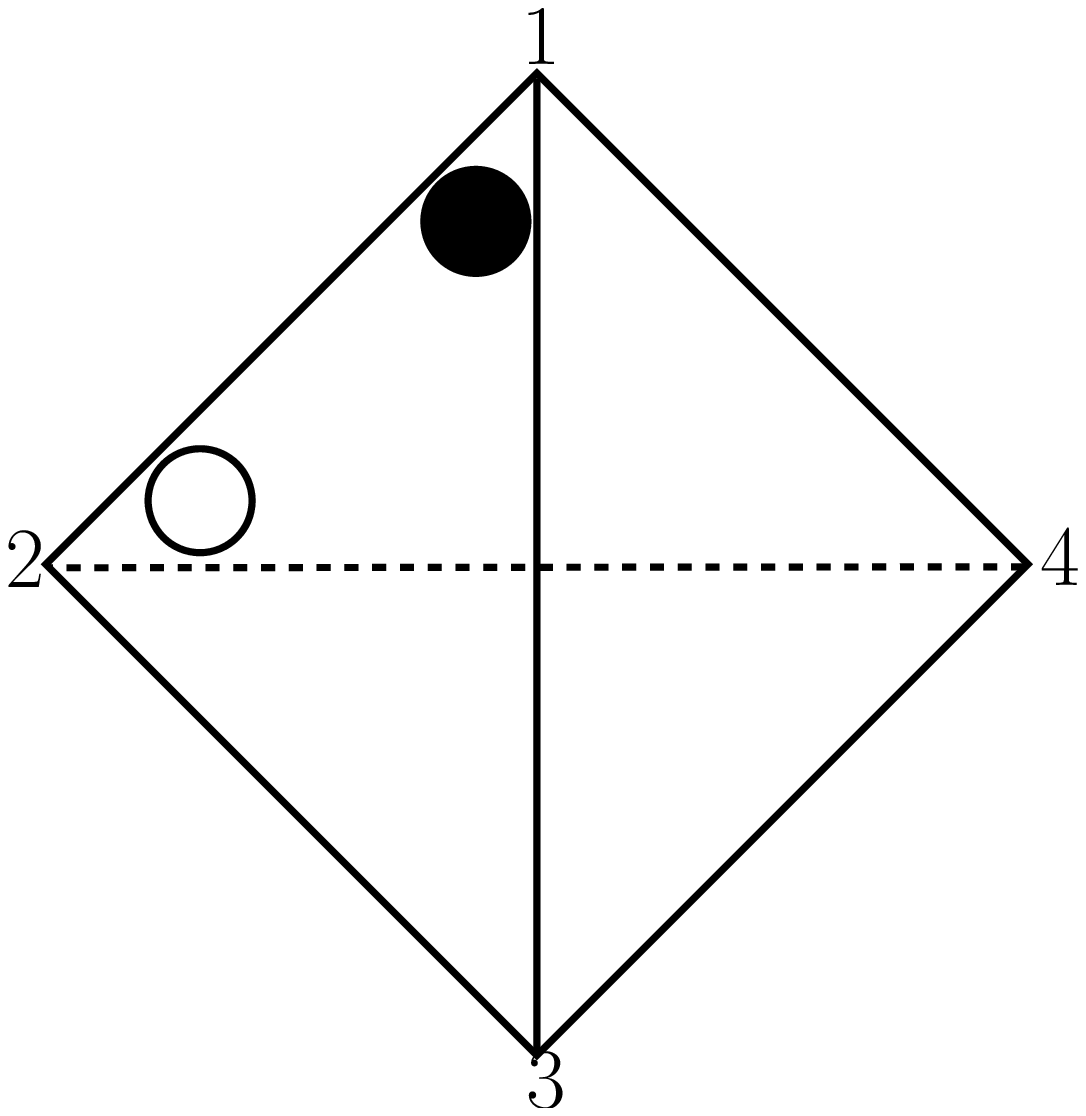}
\\
&\searrow&\downarrow&\swarrow& \medskip\\
&&
\includegraphics[width=1.65cm,height=1.65cm]{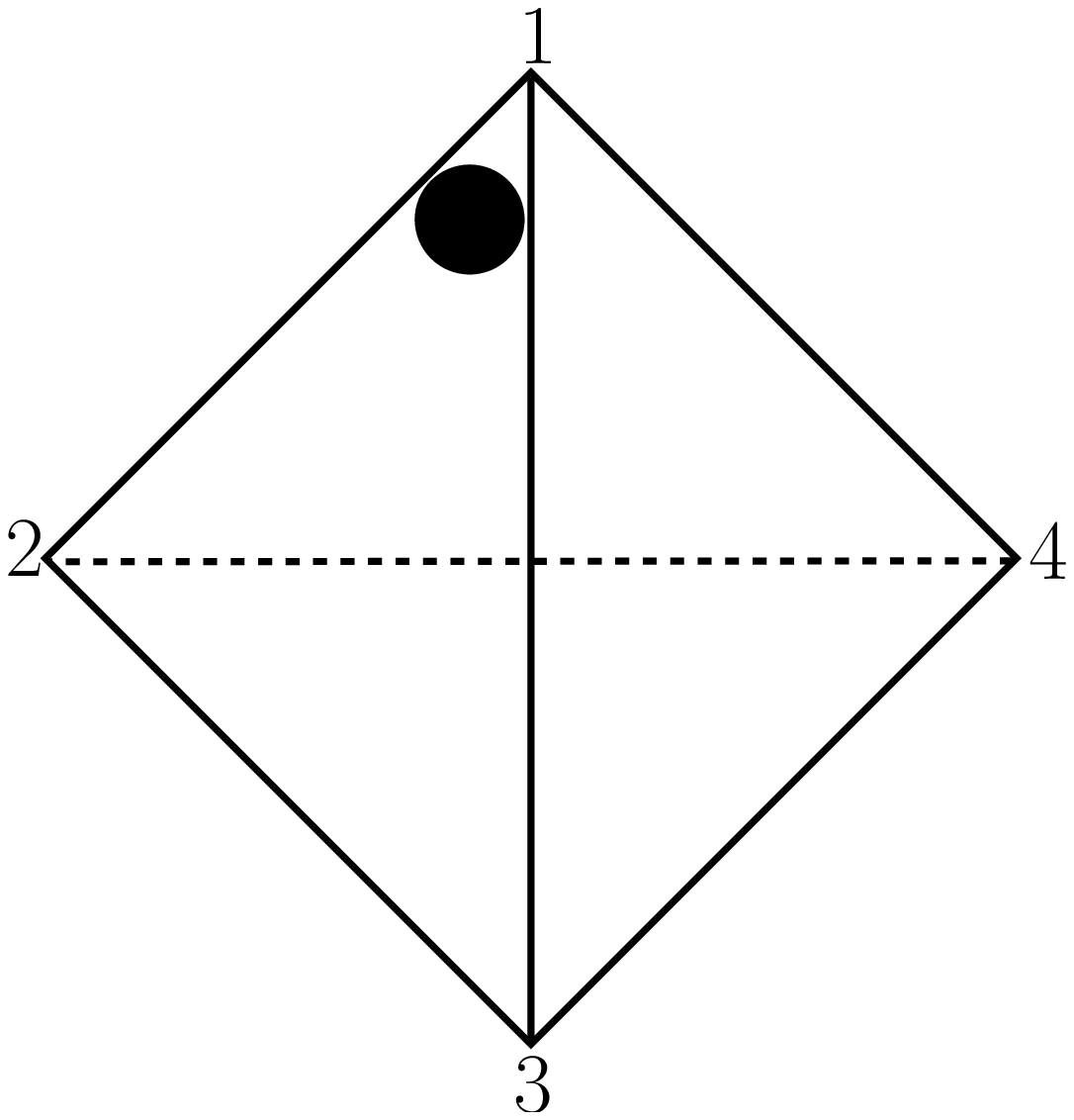}
&&
\\
\end{array}
\]
\end{center}
\caption{The subgroups of symmetries of the base tetrahedron $T$.
\label{fig:symmetries}}
\end{figure}

A group $G$ with $m=0$ is completely characterized by 
the stabilizer of the base tetrahedron $T$. Indeed,
$G$ is generated by $R_1$ together with this stabilizer. 
Hence, such groups are in bijection
to the eight (conjugacy classes of)
subgroups of  $D_8$, the symmetry group of $T$. The lattice of them
is drawn in Figure~\ref{fig:symmetries}. For each group we
have drawn a projection of $T$ in the vertical direction and, inside it,
the orbit of a generic point. 
The dots representing an orbit are drawn some black and some white, depending on whether
the corresponding point lies above or below the horizontal plane at height $1/4$, which bisects $T$.
The corresponding full cubic groups are listed in the first
column of Table~\ref{table:fullgroups}, with the same picture. 
 The three groups that contain reflection planes are shaded.
%
%
%
%
\begin{table}[htp]
\tiny
\begin{center}
\begin{tabular}{|c|c||cccc|}
\hline 
&{m=0} & \multicolumn{4}{||c|}{m=1}\\
\hline s=1 &
\begin{minipage}{1.65 cm}
 \begin{center}
 \includegraphics[width=1.65cm,height=1.65cm]{table27groups/1.eps}
 \\ $F23$
 \end{center}
\end{minipage} 
&
\definecolor{light}{gray}{.75}\colorbox{light}{
\begin{minipage}{1.55 cm}
 \begin{center}
 \includegraphics[width=1.65cm,height=1.65cm]{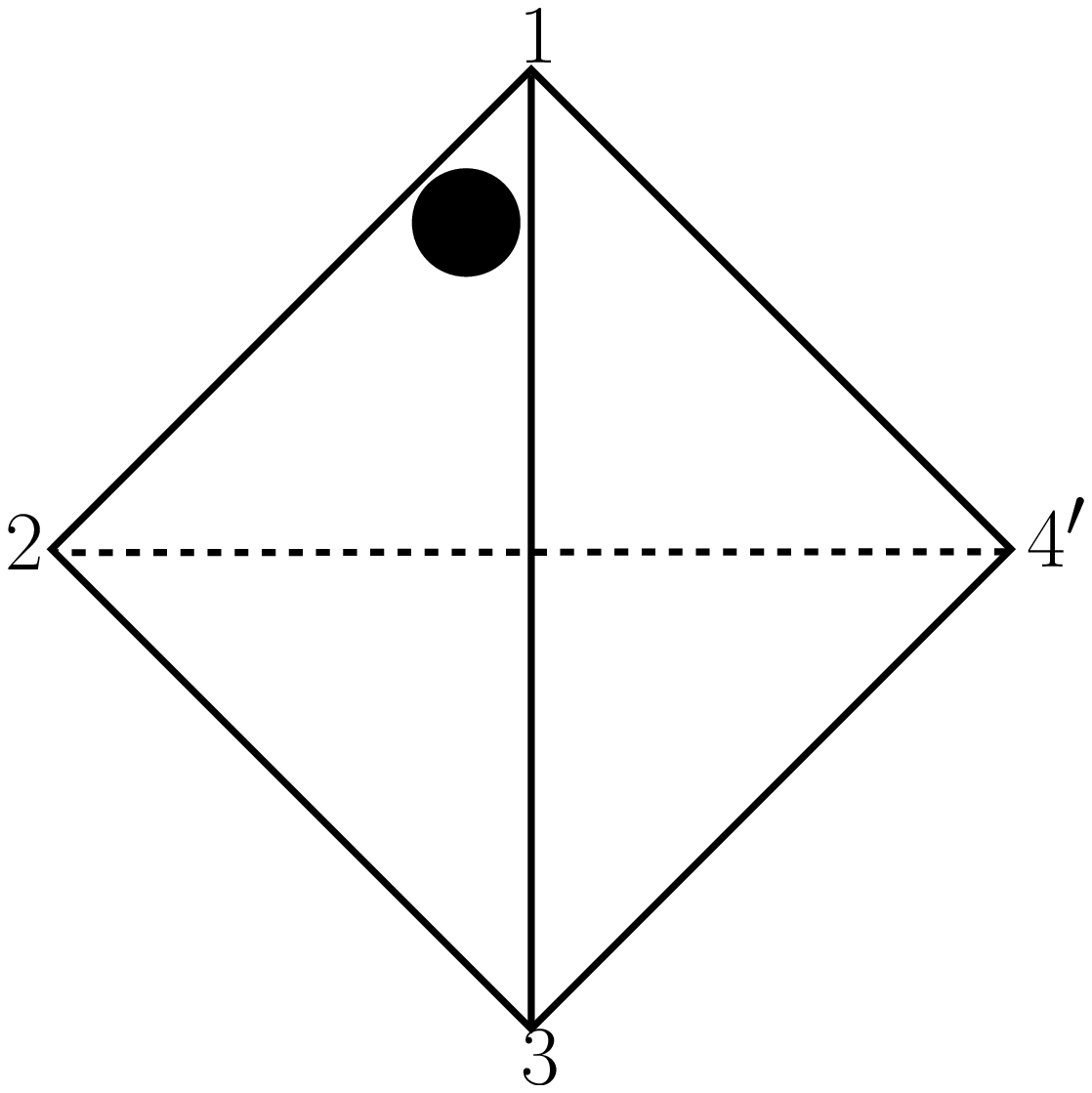}
 \\ $F\overline{4}3m$
 \end{center}
\end{minipage}
}
&
\begin{minipage}{1.65 cm}
 \begin{center}
 \includegraphics[width=1.65cm,height=1.65cm]{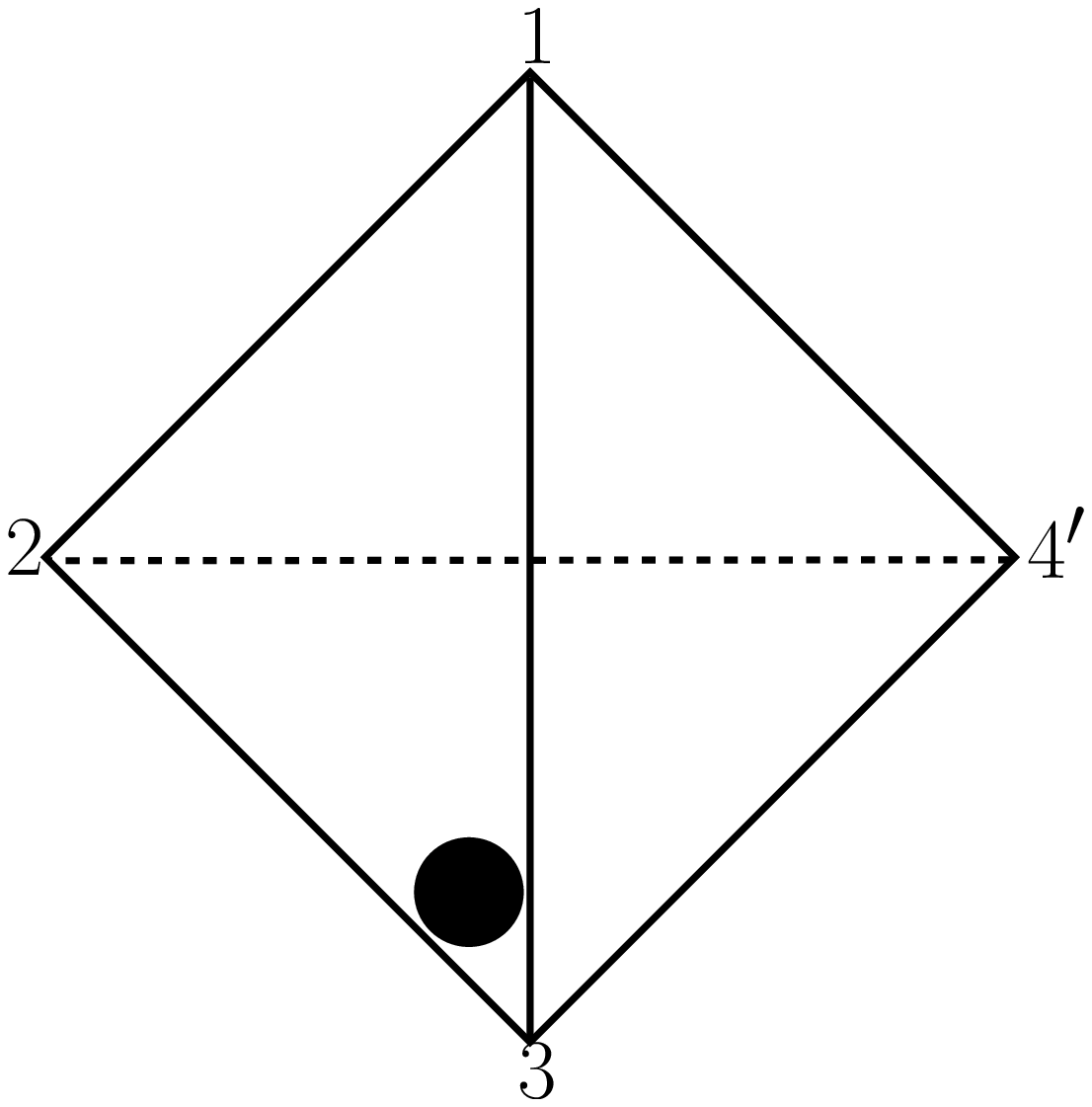}
 \\ $F432$
 \end{center}
\end{minipage} 
&
\begin{minipage}{1.65 cm}
\begin{center}
 \includegraphics[width=1.65cm,height=1.65cm]{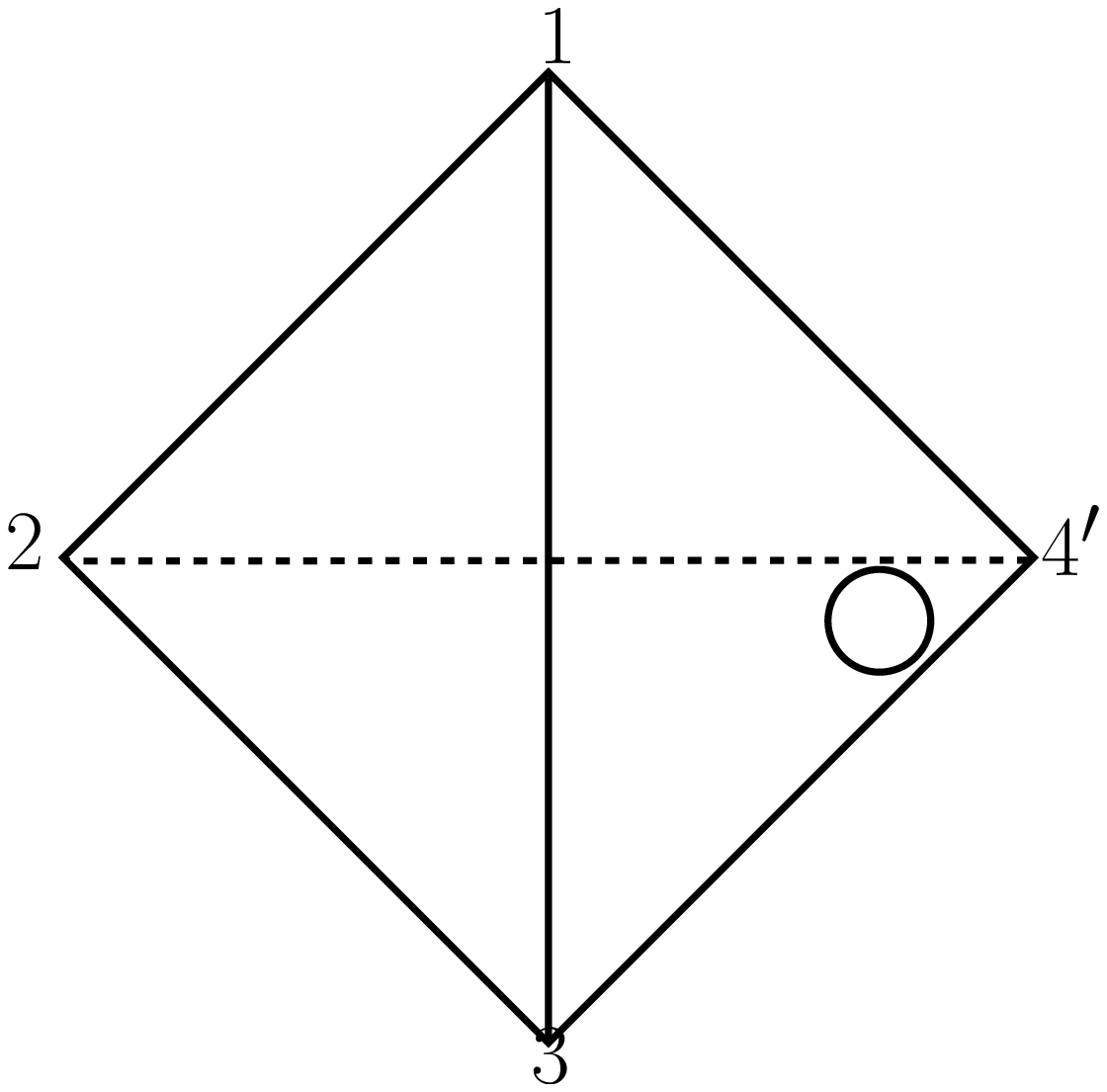}
 \\$F\frac{2}{d}\overline{3}$
 \end{center}
\end{minipage} 
&
\begin{minipage}{1.65 cm}
 \begin{center}
 \includegraphics[width=1.65cm,height=1.65cm]{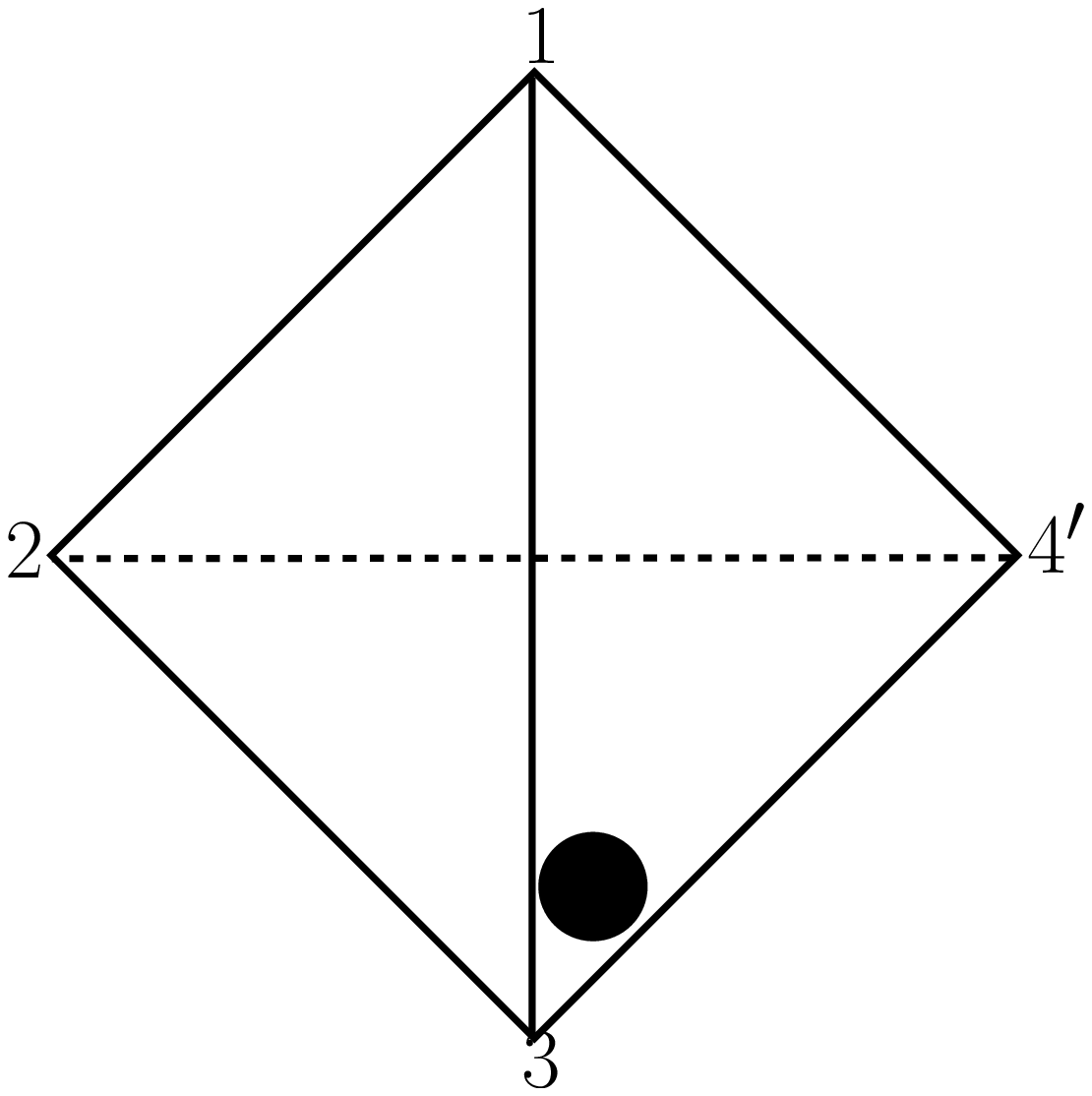}
 \\  $F\overline{4}3c$
 \end{center}
\end{minipage}
\\
\hline
\hline  &
\begin{minipage}{1.65 cm}
 \begin{center}
 \includegraphics[width=1.65cm,height=1.65cm]{table27groups/2v.eps}
 \\
 $P23$ 
 \end{center}
\end{minipage} &
\definecolor{light}{gray}{.75}\colorbox{light}{
\begin{minipage}{1.65 cm}
 \begin{center}
 \includegraphics[width=1.65cm,height=1.65cm]{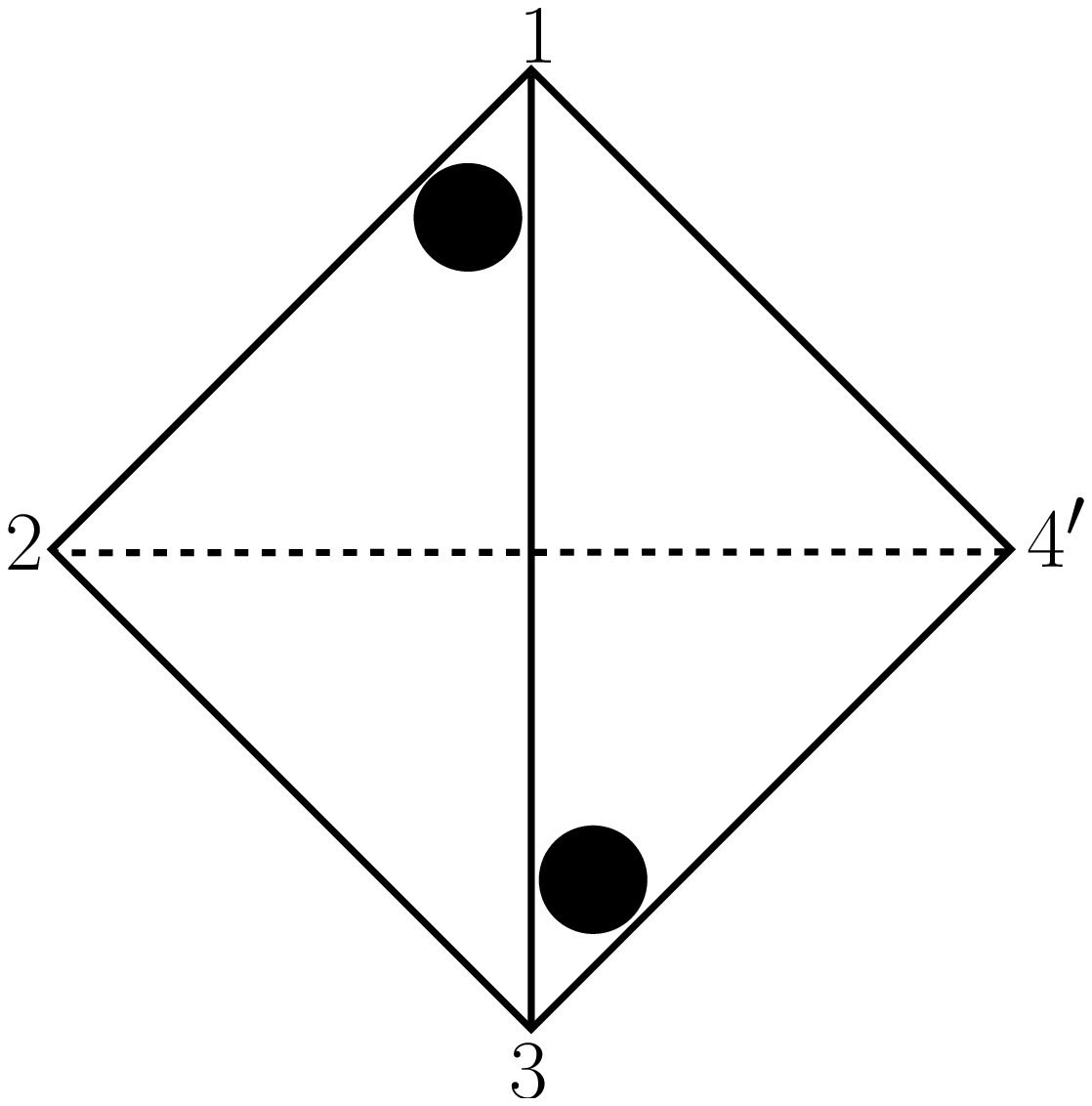}
 \\
 $P\overline{4}3m$
 \end{center}
\end{minipage}} &
\begin{minipage}{1.65 cm}
 \begin{center}
 \includegraphics[width=1.65cm,height=1.65cm]{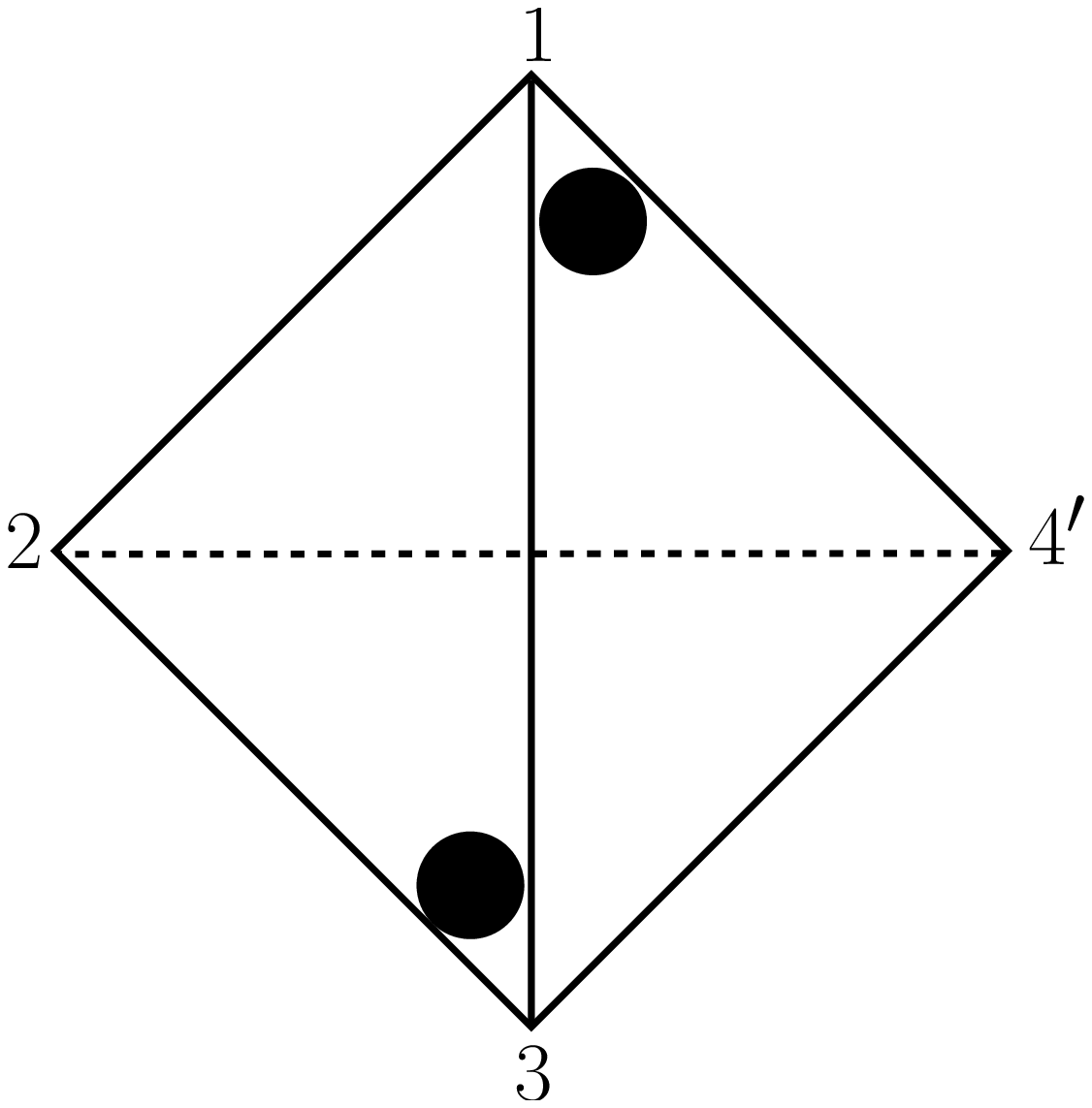}
 \\
 $P432$
 \end{center}
\end{minipage} &
\begin{minipage}{1.65 cm}
 \begin{center}
 \includegraphics[width=1.65cm,height=1.65cm]{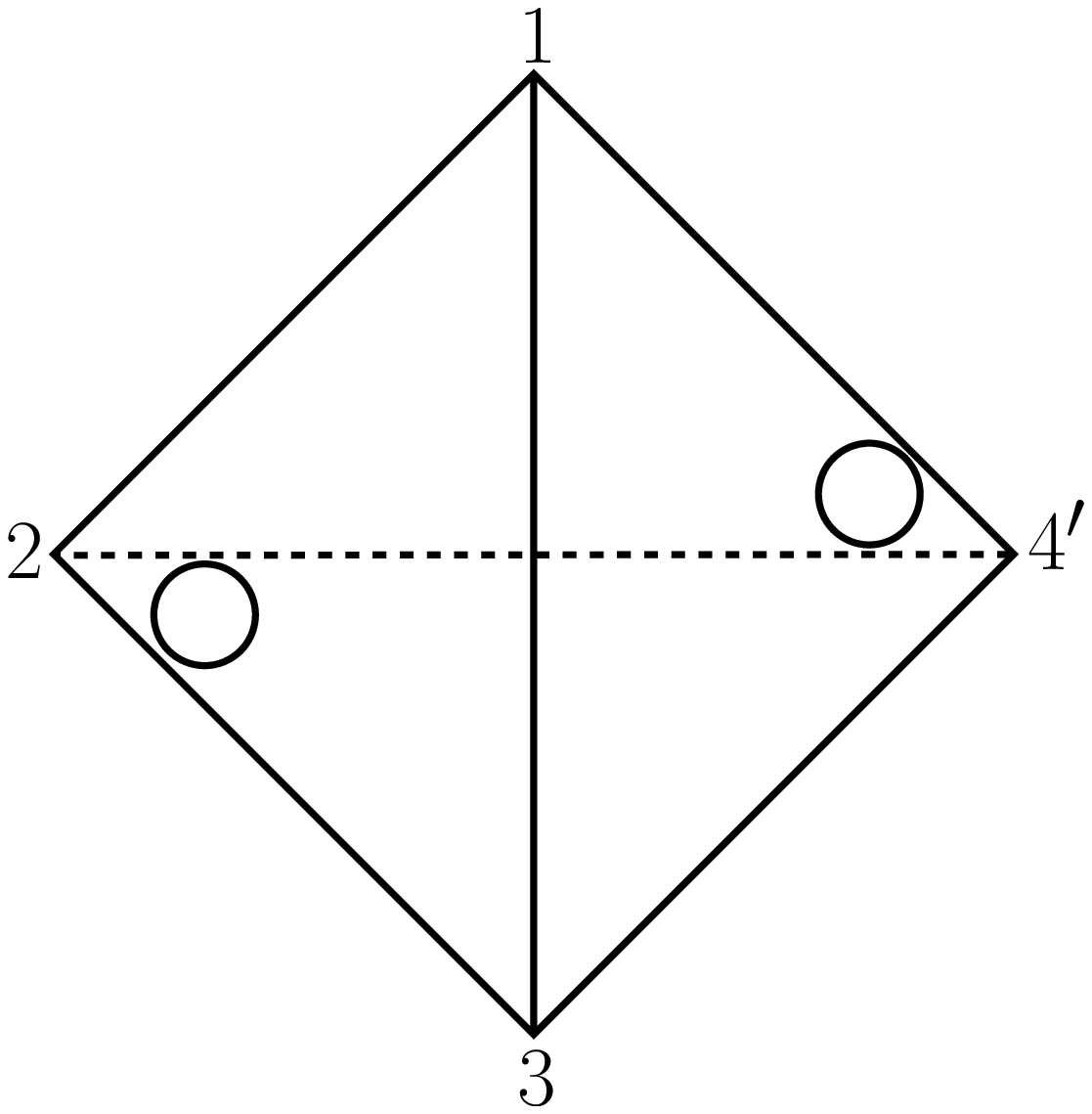}
 \\
 $I23$ 
 \end{center}
\end{minipage} &
\begin{minipage}{1.65 cm}
 \begin{center}
 \includegraphics[width=1.65cm,height=1.65cm]{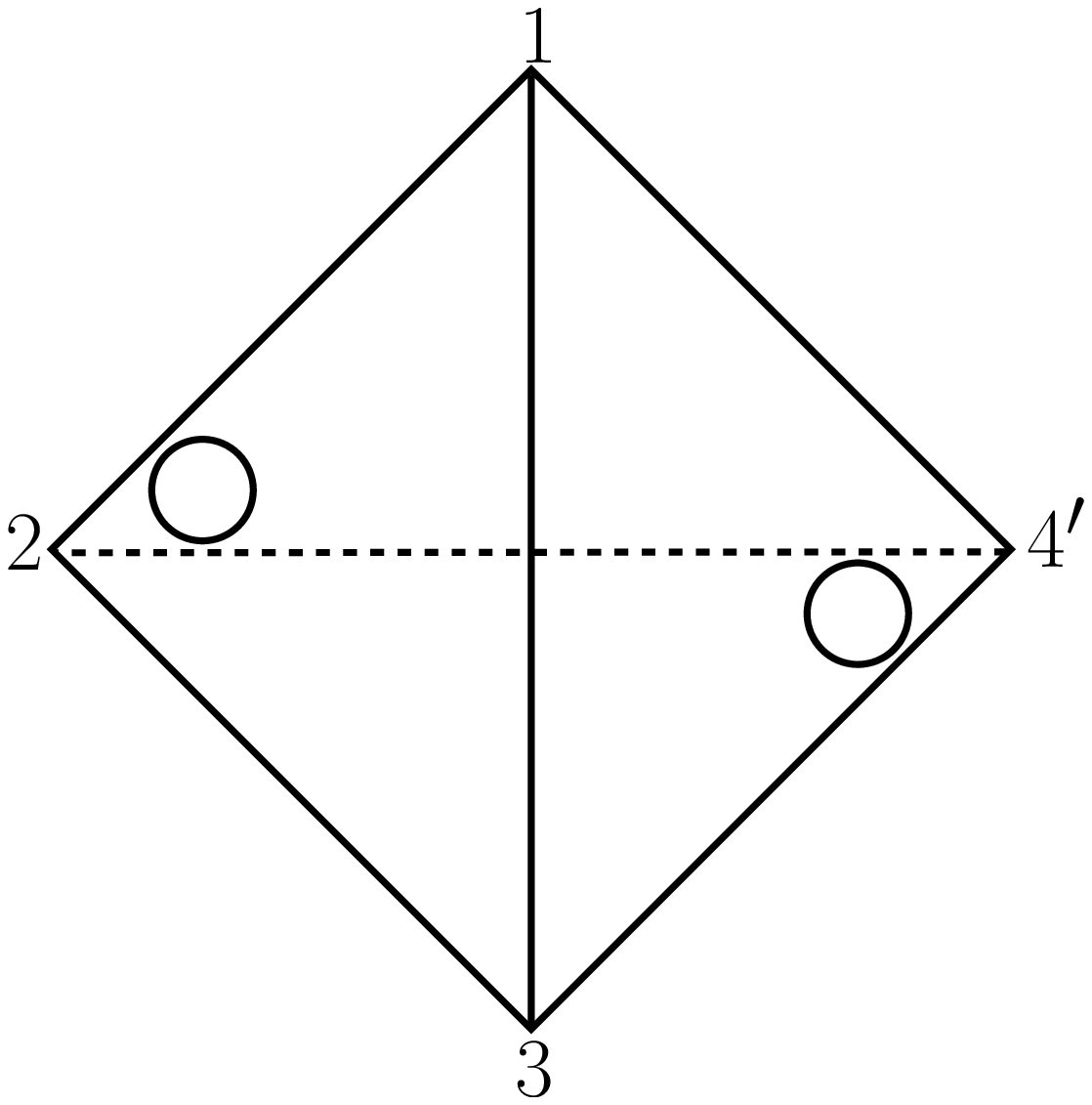}
 \\
 $P\frac{2}{n}\overline{3}$
 \end{center}
\end{minipage}
 \\
\cline{2-6}
 s=2 &
 \begin{minipage}{1.65 cm}
 \begin{center}
 \includegraphics[width=1.65cm,height=1.65cm]{table27groups/2d.eps}
 \\
 $F4_{1}32$
 \end{center}
\end{minipage} &
\definecolor{light}{gray}{.75}\colorbox{light}{
\begin{minipage}{1.65 cm}
 \begin{center}
 \includegraphics[width=1.65cm,height=1.65cm]{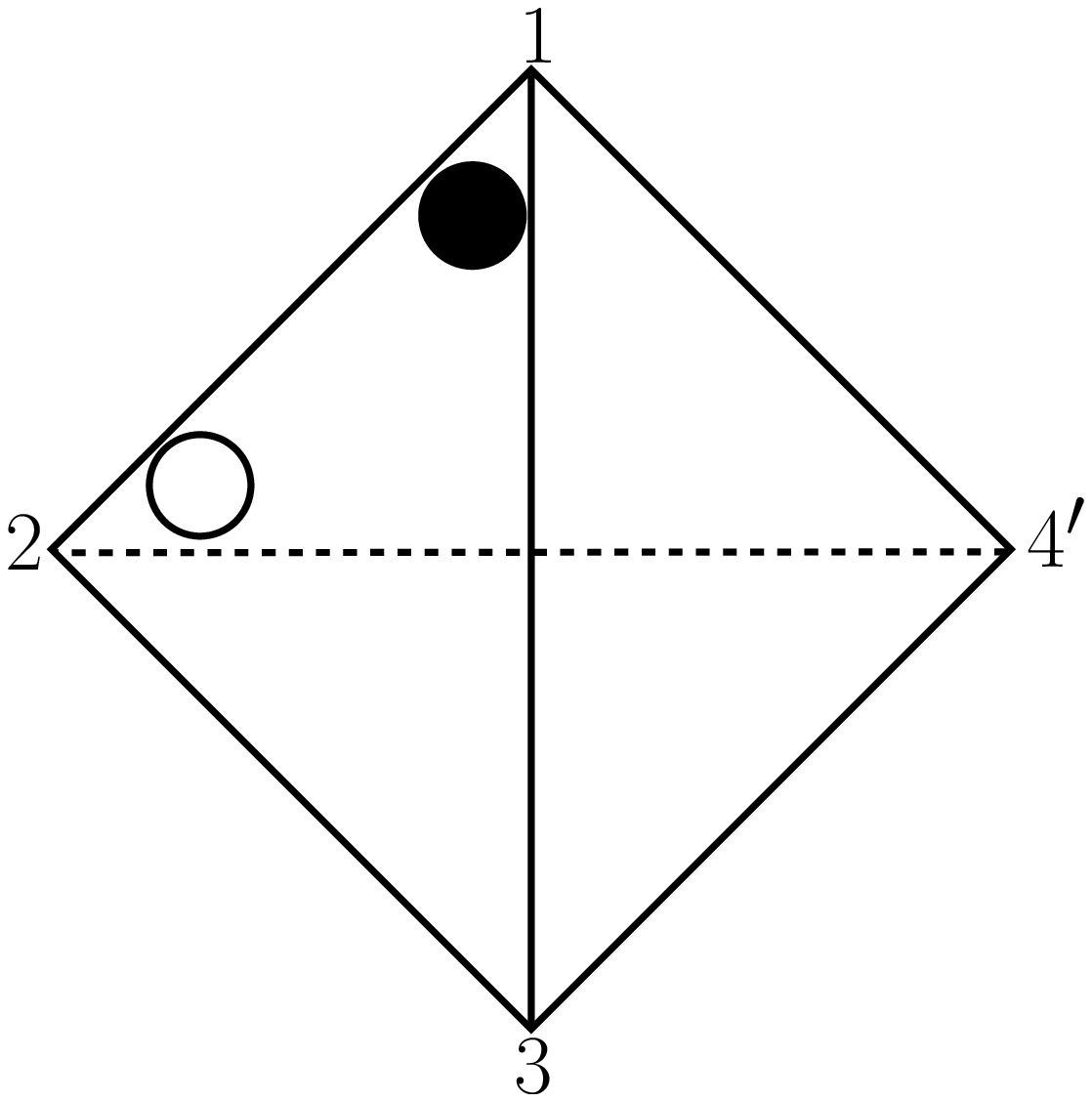}
 \\
 $F\frac{4}{m}\overline{3}\frac{2}{n}$
 \end{center}
\end{minipage}} 
&
\begin{minipage}{1.65 cm}
 \begin{center}
 \includegraphics[width=1.65cm,height=1.65cm]{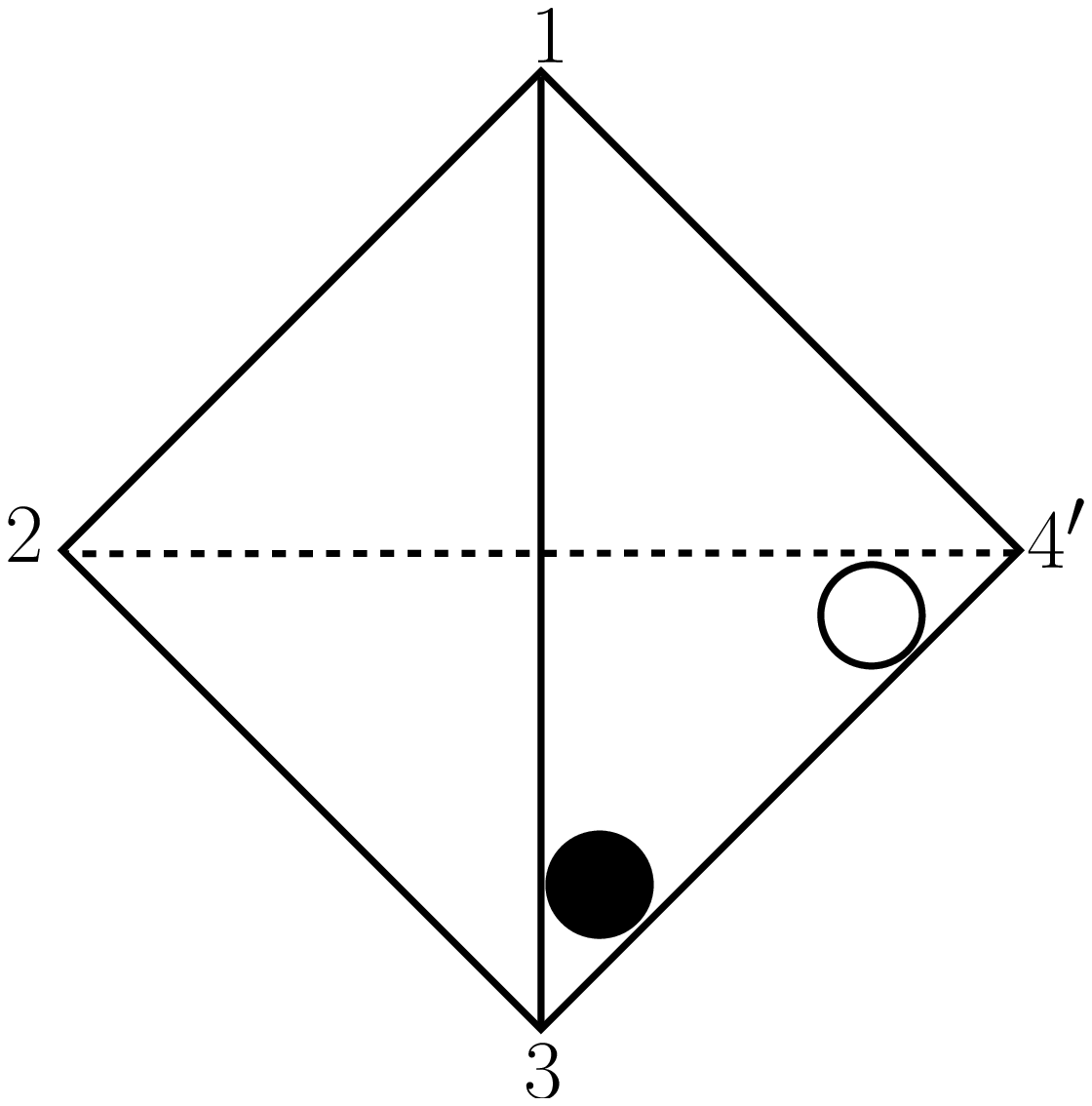}
 \\
 $F\frac{4_{1}}{d}\overline{3}\frac{2}{n}$
 \end{center}
\end{minipage} 
&
&
\\
\cline{2-6} &
\definecolor{light}{gray}{.75}\colorbox{light}{
\begin{minipage}{1.65 cm}
 \begin{center}
 \includegraphics[width=1.65cm,height=1.65cm]{table27groups/m.eps}
 \\
 $F\frac{2}{m}\overline{3}$
 \end{center}
\end{minipage}} &
\definecolor{light}{gray}{.75}\colorbox{light}{
\begin{minipage}{1.65 cm}
 \begin{center}
 \includegraphics[width=1.65cm,height=1.65cm]{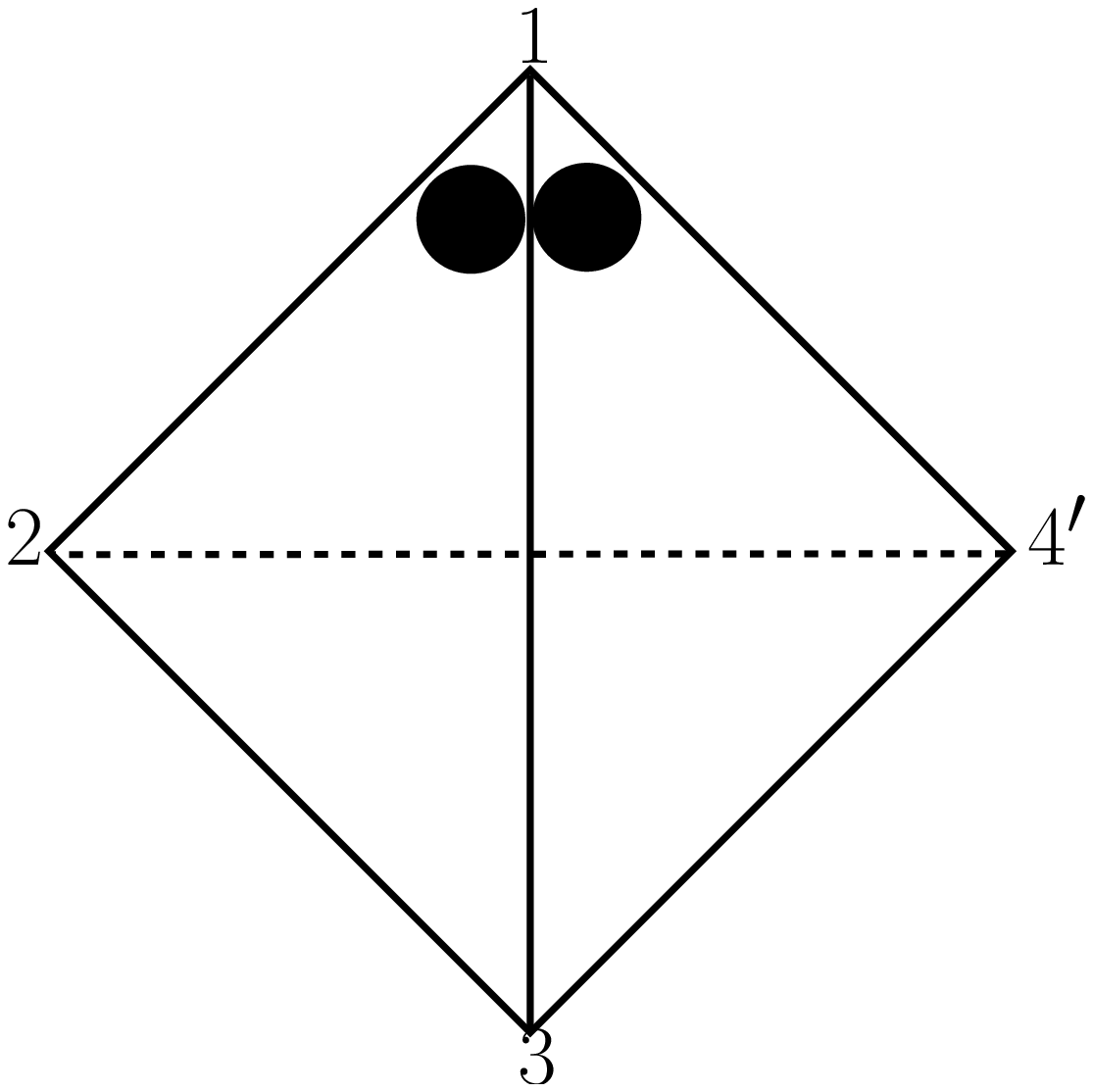}
 \\
 $F\frac{4}{m}\overline{3}\frac{2}{m}$
 \end{center}
\end{minipage}} &
\definecolor{light}{gray}{.75}\colorbox{light}{
\begin{minipage}{1.65 cm}
 \begin{center}
 \includegraphics[width=1.65cm,height=1.65cm]{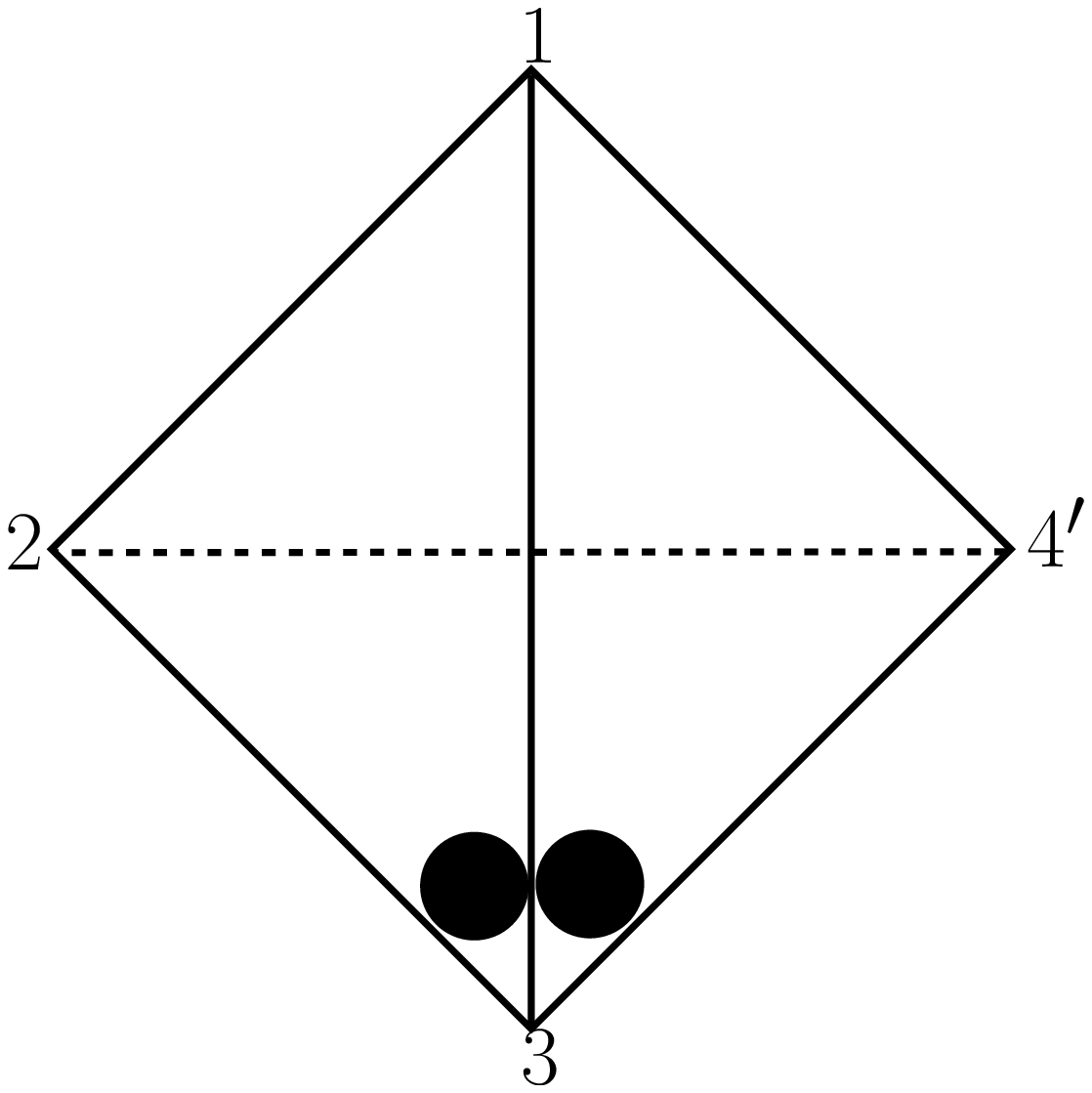}
 \\
 $F\frac{4_{1}}{d}\overline{3}\frac{2}{m}$
 \end{center}
\end{minipage}} & & \\
\hline
\hline &
\begin{minipage}{1.65 cm}
 \begin{center}
 \includegraphics[width=1.65cm,height=1.65cm]{table27groups/22.eps}
 \\
$P4_{2}32$
 \end{center}
\end{minipage} &
\definecolor{light}{gray}{.75}\colorbox{light}{
\begin{minipage}{1.65 cm}
 \begin{center}
 \includegraphics[width=1.65cm,height=1.65cm]{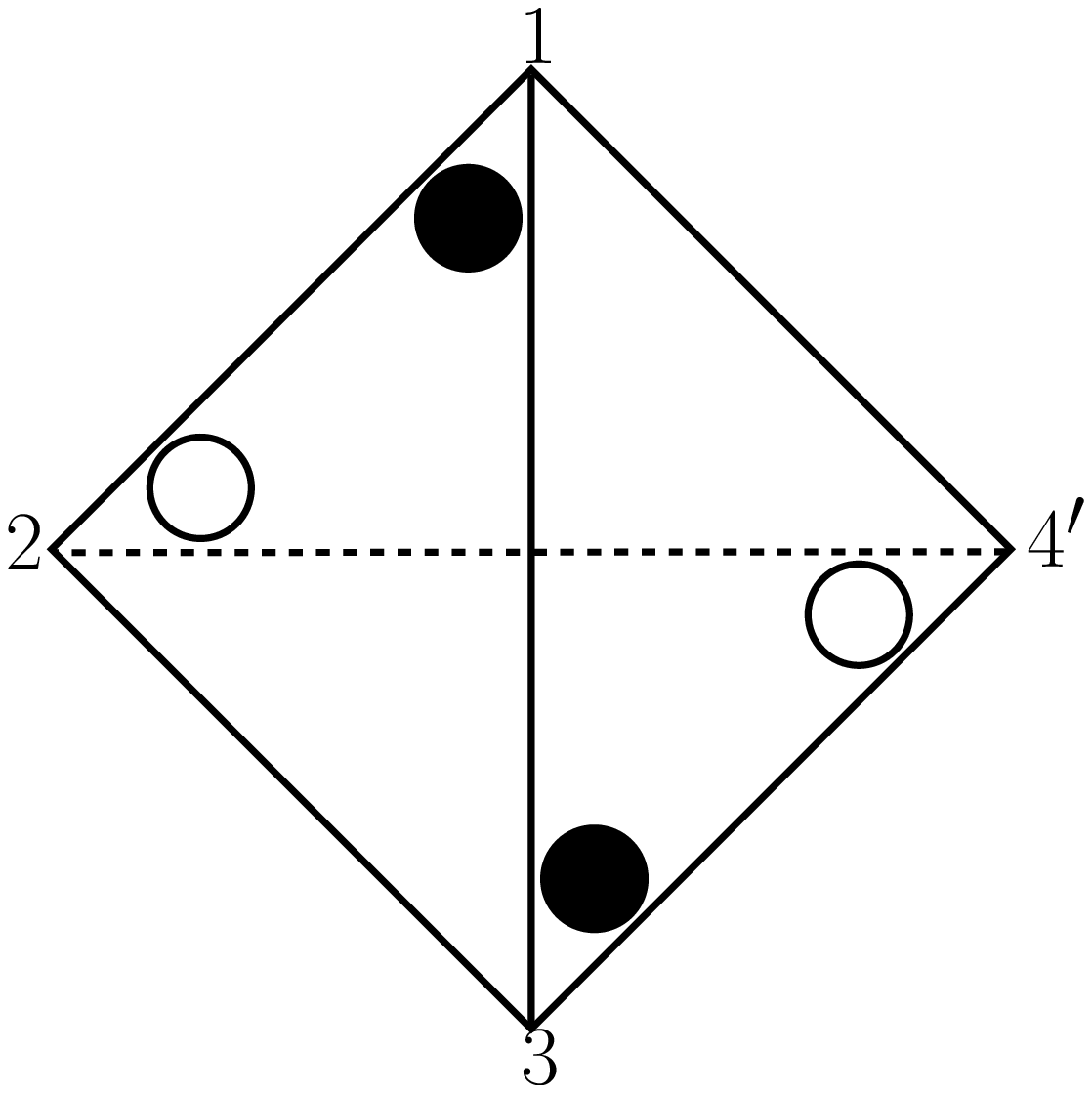}
 \\$P\frac{4_{1}}{n}\overline{3}\frac{2}{m}$
 \end{center}
\end{minipage}} &
\begin{minipage}{1.65 cm}
 \begin{center}
 \includegraphics[width=1.65cm,height=1.65cm]{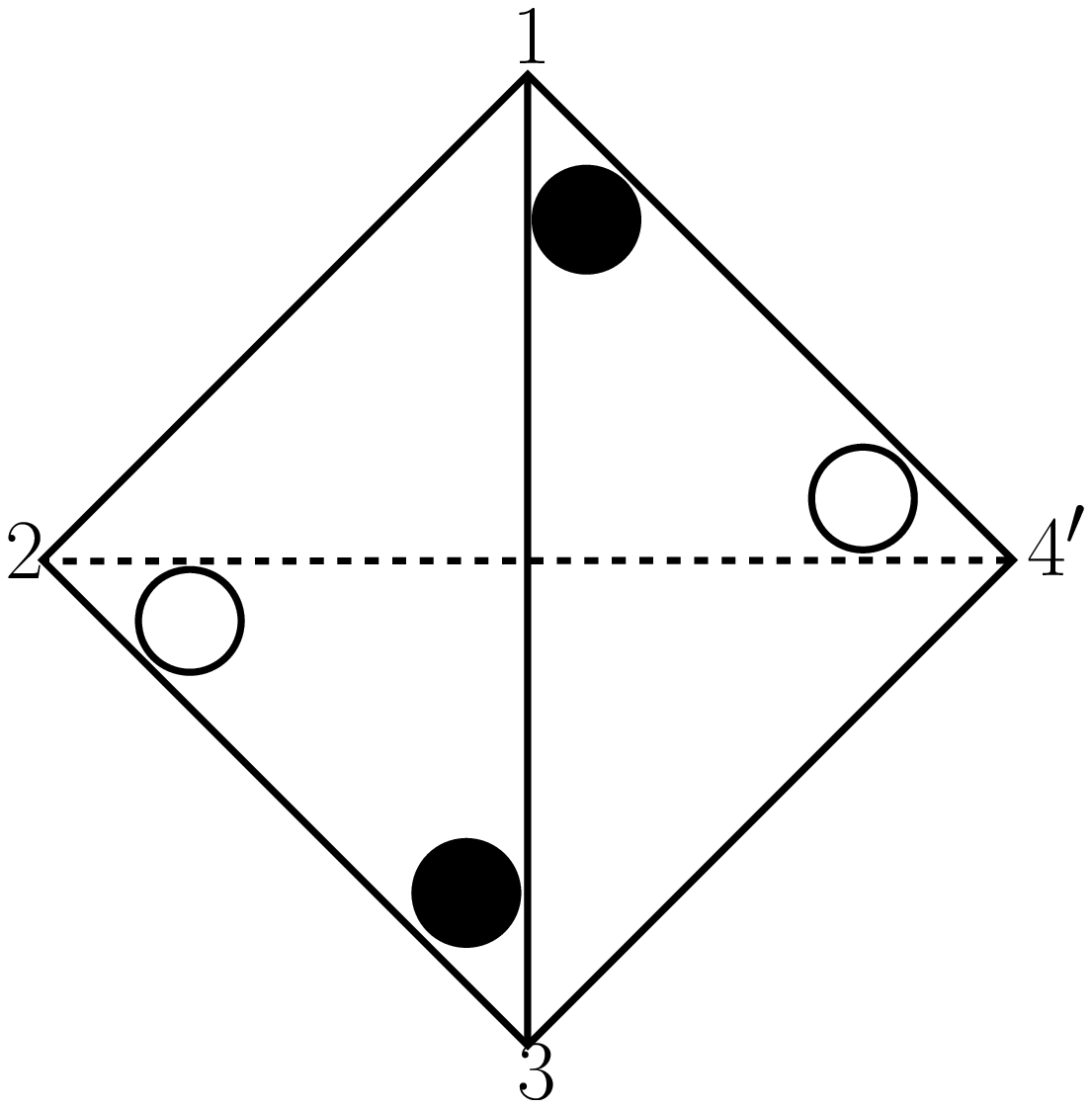}
 \\
 $I432$
 \end{center}
\end{minipage} & &
\\
\cline{2-6} s=4 &
\begin{minipage}{1.65 cm}
 \begin{center}
 \includegraphics[width=1.65cm,height=1.65cm]{table27groups/4.eps}
 \\
 $P\overline{4}3n$
 \end{center}
\end{minipage} &
\definecolor{light}{gray}{.75}\colorbox{light}{
\begin{minipage}{1.65 cm}
 \begin{center}
 \includegraphics[width=1.65cm,height=1.65cm]{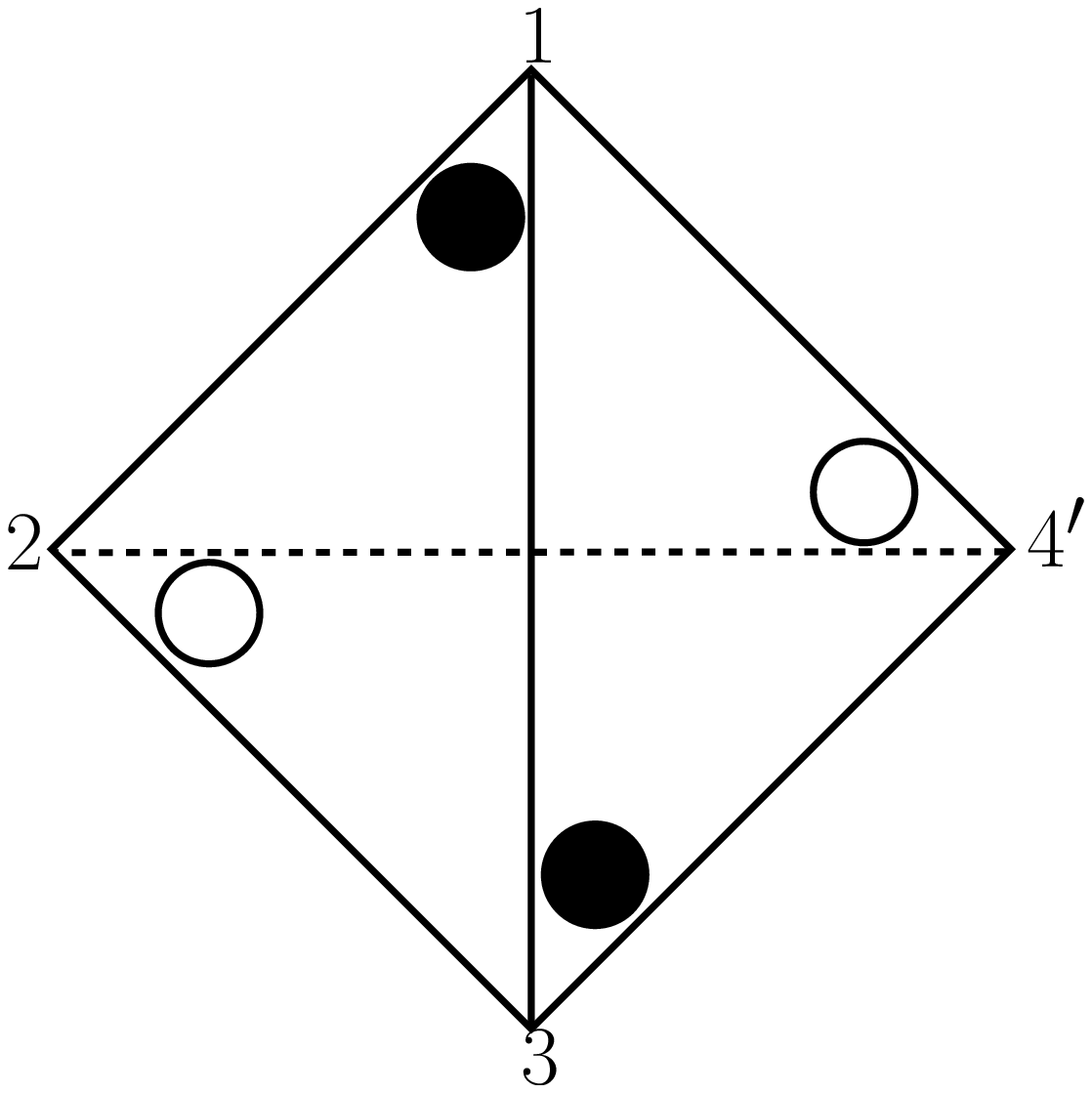}
 \\
 $I\overline{4}3m$
  \end{center}
\end{minipage}} &
\begin{minipage}{1.65 cm}
 \begin{center}
 \includegraphics[width=1.65cm,height=1.65cm]{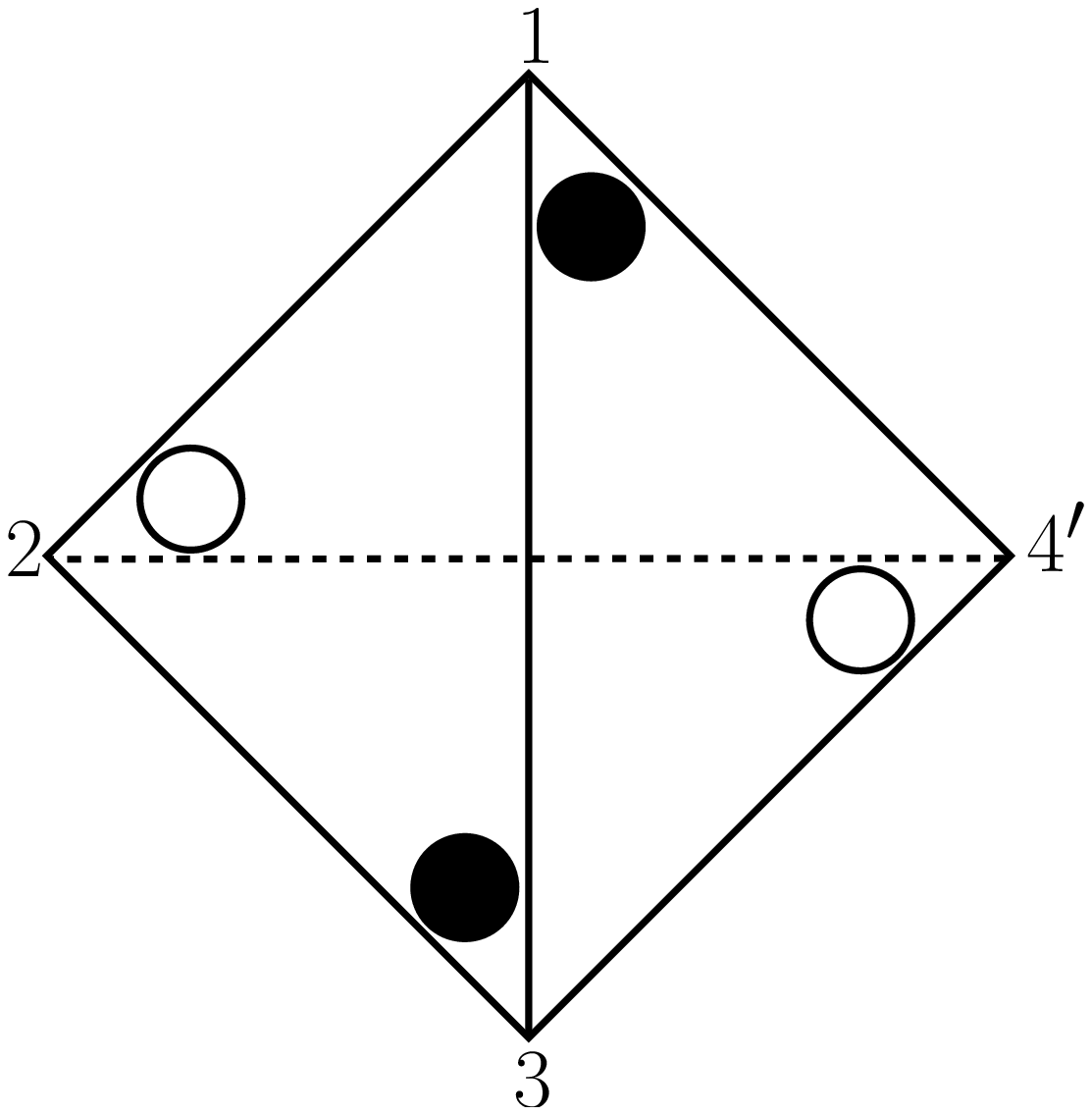}
 \\
 $P\frac{4}{n}\overline{3}\frac{2}{n}$ 
 \end{center}
\end{minipage} & &
 \\
\cline{2-6}  &
\definecolor{light}{gray}{.75}\colorbox{light}{
\begin{minipage}{1.65 cm}
 \begin{center}
 \includegraphics[width=1.65cm,height=1.65cm]{table27groups/2m.eps}
 \\
 $P\frac{2}{m}\overline{3}$
 \end{center}
\end{minipage}} &
\definecolor{light}{gray}{.75}\colorbox{light}{
\begin{minipage}{1.65 cm}
 \begin{center}
 \includegraphics[width=1.65cm,height=1.65cm]{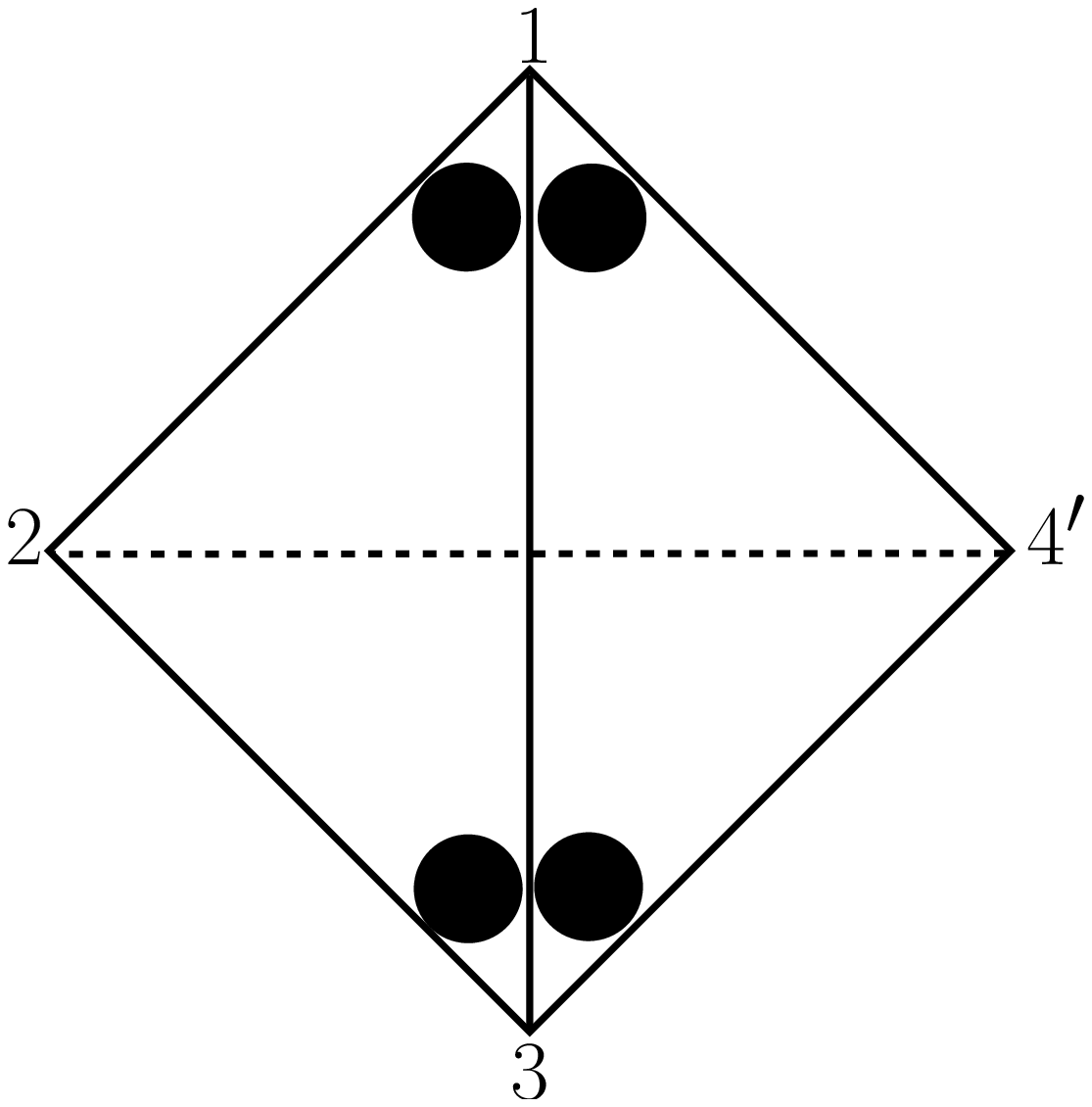}
 \\
 $P\frac{4}{m}\overline{3}\frac{2}{m}$
 \end{center}
\end{minipage}} &
\definecolor{light}{gray}{.75}\colorbox{light}{
\begin{minipage}{1.65 cm}
 \begin{center}
 \includegraphics[width=1.65cm,height=1.65cm]{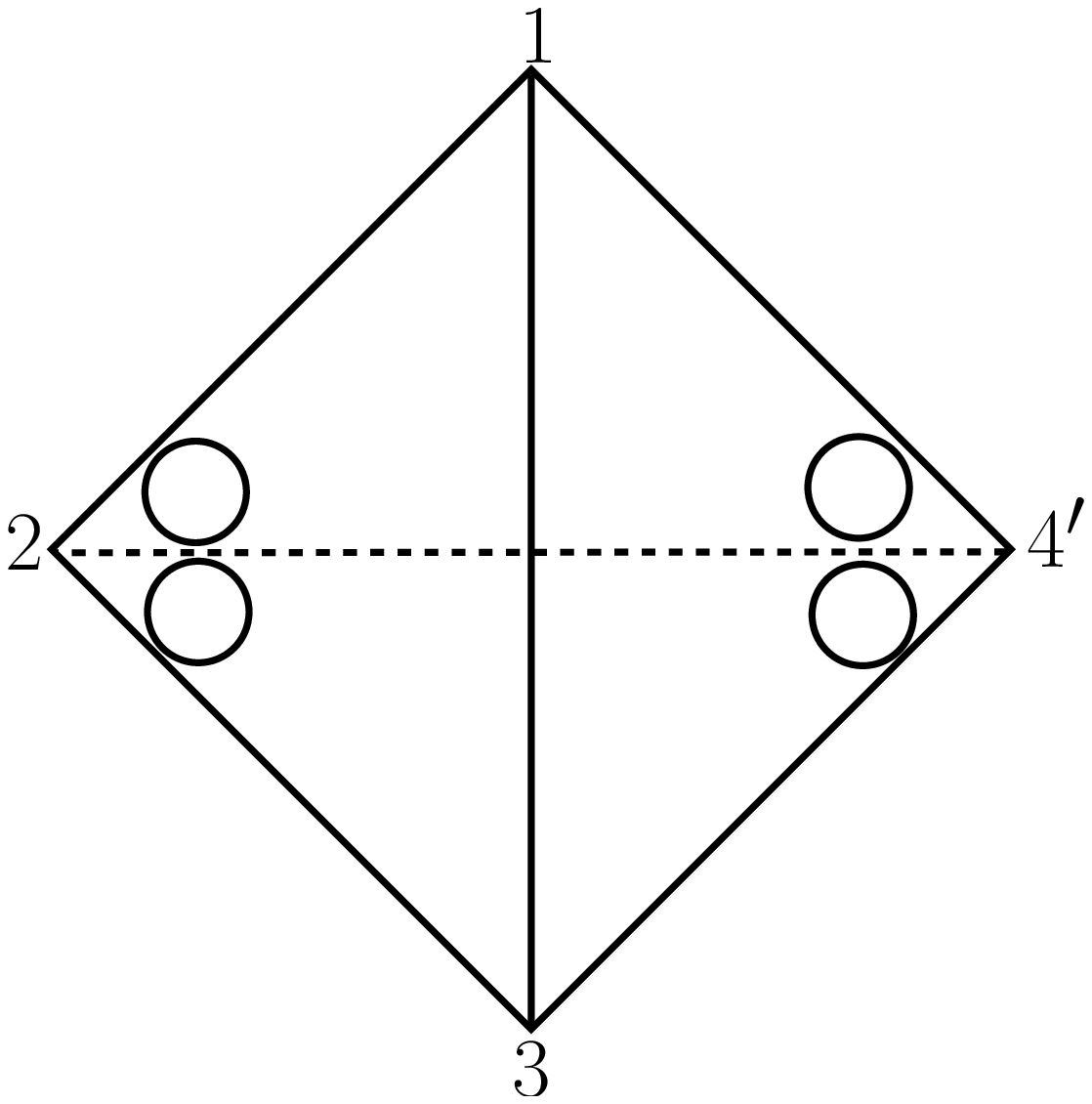}
 \\
 $I\frac{2}{m}\overline{3}$
 \end{center}
\end{minipage}} & &
\\
\hline
\hline s=8 &
\definecolor{light}{gray}{.75}\colorbox{light}{
\begin{minipage}{1.65 cm}
 \begin{center}
 \includegraphics[width=1.65cm,height=1.65cm]{table27groups/4m.eps}
 \\
 $P\frac{4}{m}\overline{3}\frac{2}{n}$
 \end{center}
\end{minipage}} &
\definecolor{light}{gray}{.75}\colorbox{light}{
\begin{minipage}{1.65 cm}
 \begin{center}
 \includegraphics[width=1.65cm,height=1.65cm]{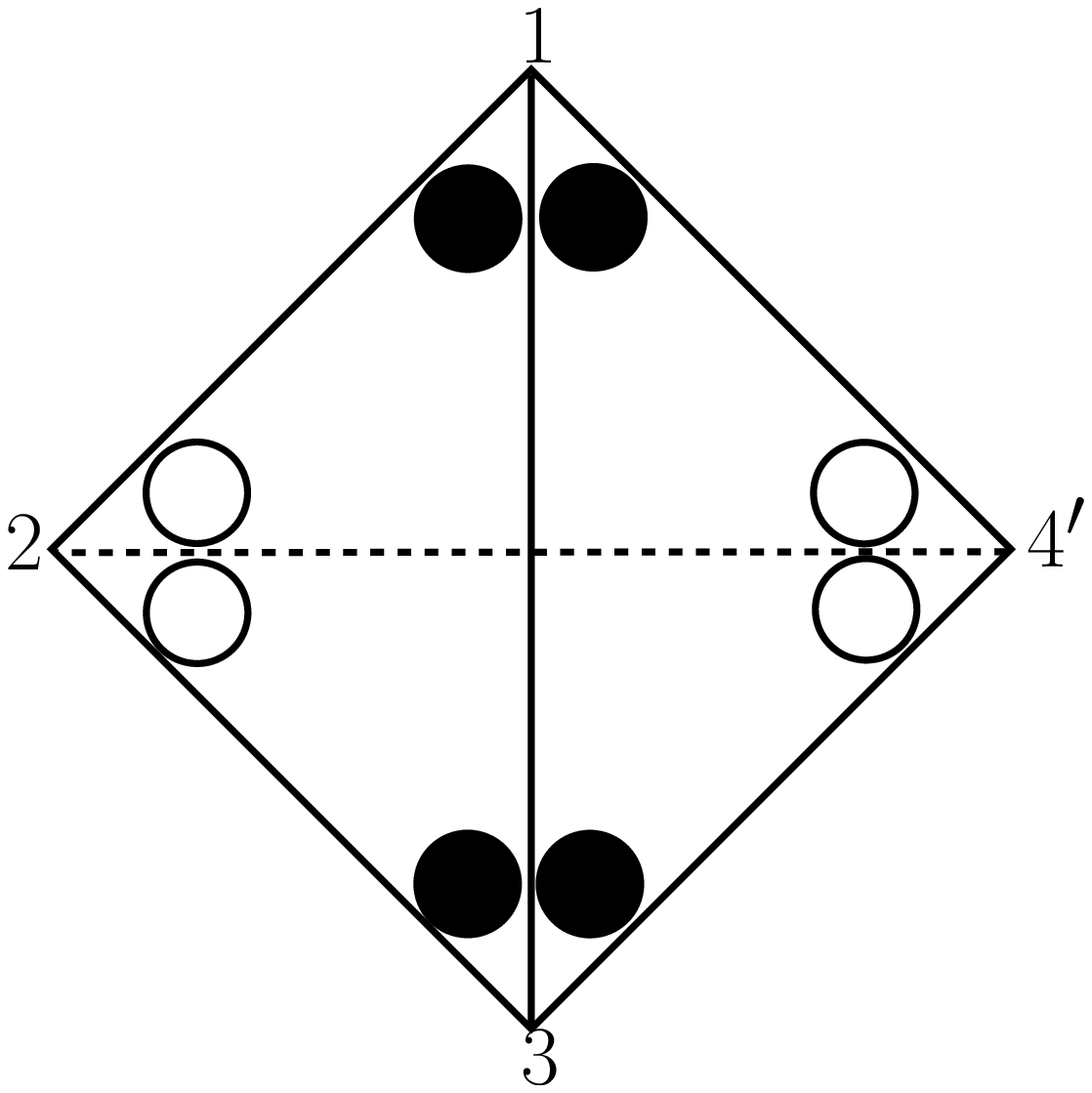}
 \\
 $I\frac{4}{m}\overline{3}\frac{2}{m}$
 \end{center}
\end{minipage}} & & &
\\
\hline
\end{tabular}
\end{center}
\medskip
\caption{The 27 full cubic groups.}
\label{table:fullgroups}
\end{table}
 
 For a group $G$ with $m=1$, we first look at the transformations in $G$ that
 send Delaunay tetrahedra to Delaunay tetrahedra of the same $R_1$-orbit. They form a 
 subgroup $H$ of index two in $G$, and $H$ is  one of the eight groups of the previous paragraph. Then, $G$ is generated
 by $H$ together with any single transformation $\rho\in G$ that sends $T$ to one of its neighbors, say  $T_4$. Let $H_T$ denote the stabilizer of $T$ in $H$, and let $\rho_T$
 be $\rho$ composed with the mirror reflection on the plane through $v_1$, $v_2$ and $v_3$ (that is, the unique  transformation that sends $T_4$ to $T$, keeping
 the labels of vertices). Let $G_T$ be the group generated by $H_T$ and $\rho_T$.
 Then, $G$ is completely characterized by $H_T$ and $G_T$ and the conditions on
 $H_T$ and $G_T$ are that  $H_T\le G_T \le D_8$ and $|G_T : H_T|$ equals 1 or 2.

 To picture these groups, in Table~\ref{table:fullgroups} we show the image in $T_4$ of the
 orbit of $H$ pictured in $T$. The rule of thumb is that the two orbits, when superimposed by 
 composing one with the mirror symmetry on the plane containing $v_1$, $v_2$ and $v_3$
 should together form another orbit of a subgroup of $D_8$. This happens automatically
 when $s=8$ or $s=4$. For $s=1$, it is equivalent to the transformation between
 the two orbits being an involution of $T$. For $s=2$, it is necessary but not sufficient
 that the transformation be an involution.
%
  
 Depending on the particular $H$, there may be one or more possibilities for $G$. There 
 is always at least one, in which $\rho$ is the mirror symmetry itself (and the involution
 is the identity). That is the first group shown after $H$ in each row of Table~\ref{table:fullgroups},
 and it is shaded since it obviously contains reflections. The other ones contain reflections if
 and only if $H$ contains reflections.

 \section{A first upper bound}
  \label{sec:first-bound}

 Let now $G$ be one of the 27 full groups, let $T$ be a Delaunay tetrahedron of the lattice $I$,
 and $p$ a point in $T$, whose Dirichlet stereohedron we want to study. We assume that $p$ has trivial stabilizer for the action of $G$ on $Gp$. That is to say, that the only motion in $G$ that fixes $p$ is the identity. In the appendix we treat the case of non-trivial stabilizers, and correct statements made about this in the previous papers~\cite{Bochis-Santos-2001,Bochis-Santos-2006}.
 %
 The main reason for this assumption is the following lemma:
 
 \begin{lemma}
 \label{lemma:stabilizer}
 If $p$ has trivial stabilizer, then a small perturbation of it cannot decrease the number of facets of the Dirichlet stereohedron $\Vor_{Gp}(p)$. That is, the maximum number of facets for a Dirichlet stereohedron  with trivial stabilizer is achieved when $p$  is generic.
 \end{lemma}
 
As a consequence, from now on we assume that $p$ is generic, in particular that it lies in the interior of~$T$.
 
 \begin{proof}
Let $p_1\in G p$ is a neighbor of $p$ in the Voronoi diagram. We claim that when $p$ is perturbed to a position $p'$ sufficiently close to $p$ the point $p'_1\in Gp'$ corresponding to $p_1$ will be a neighbor of $p'$.

For this, let $x\in \R^3$ be a point whose distance to $p$ and $p_1$ is equal, and smaller than the distance to any other  point in $Gp$. By continuity, a small perturbation of $p$ will make $p'$ and $p'_1$ still be the two orbit points closest to $x$, except perhaps not equidistant to it. Assume, without loss of generality, that after the perturbation we have $d(x,p') > d(x,p'_1)$. Then, moving $x$ in a straight-line towards $p'$ makes $d(x,p')$ decrease faster than the distance to any other point in $Gp'$. Hence, when $x$ reaches a position equidistant to $p'$ and $p'_1$, the distance to all other orbit points will be bigger, which proves that $p'$ and $p'_1$ are neighbors in $Vor_{Gp'}$.
 \end{proof}

 \begin{lemma}
  \label{lemma:5tetrahedra}
 The Dirichlet stereohedron $\Vor_{Gp}(p)$ is fully contained in the (non-convex) polyhedron
 obtained as the union of $T$ and its four neighbors $T_1$, $T_2$, $T_3$ and $T_4$. Moreover,
 no facet of $\Vor_{Gp}(p)$ lies in the boundary of this polyhedron. 
 \end{lemma}
 
 \begin{proof}
 We start with the following simple observation: if a group $G$ contains a rotation $\rho$ of angle $\alpha$
 and axis on the line $l$, then $\Vor_{Gp}(p)$ is fully contained on the dihedral sector of angle
 $\alpha$ with edge $l$ and centered around $p$. Indeed, this dihedral sector equals 
 $\Vor_{S}(p)$ where $S\subseteq Gp$ is the orbit of $p$ under the cyclic group generated by $\rho$.
 
 In our case, since the six edges of $T$ are rotation axes for every full group, $\Vor_{Gp}(p)$
 will be contained in the intersection of the six corresponding dihedral sectors.
 For each rotation of order
 three on edges of $T$, say the one on the edge $v_1v_2$, 
 the sector of angle $2\pi/3$ in question is in turn 
 contained in the half-space delimited by the plane containing the vertices $v_1$, $v_2$, $v'_3$ and 
 $v'_4$. The intersection of the four half-spaces is an infinite square prism circumscribed to $T$.
 Similarly, the rotations of order two on the axes $v_1v_3$ and $v_2v_4$, cut
 this infinite square prism  by two sectors of $3\pi/2$ with dihedral angle at the edges
 $v_1v_3$ and $v_2v_4$. The non-convex polyhedron obtained in this way coincides with the union of
 $T$ and its four neighbors $T_1$, $T_2$, $T_3$ and~$T_4$. 
 
 For the ``moreover'' observe that
 $p$ being interior to $T$ implies that the half-spaces and dihedral sectors can be taken open, except
 for the  corresponding axes. 
 \end{proof}
 
 \begin{lemma}
 \label{lemma:15tetrahedra}
 Let $q\in Gp$ be a point whose Dirichlet stereohedron $\Vor_{Gp}(q)$ shares a facet with 
 $\Vor_{Gp}(p)$. Then, $q$ lies in the interior of one of the following 15 Delaunay tetrahedra
 of $I$: $T$, its four neighbors $T_i$, or the ten additional neighbors of the latter $T_{ij}$.
 \end{lemma}

 \begin{proof}
 By definition, $p$ and $q$ are neighbors if and only if $\Vor_{Gp}(p)$ and $\Vor_{Gp}(q)$
share a facet. Lemma~\ref{lemma:5tetrahedra} implies that this facet (except perhaps some of its
edges) is contained in the interior of the union of $T$ and its neighbors. The same lemma applied to
$q$, implies that $q$ lies either in $T$, or in a neighbor of $T$, or in a neighbor of a neighbor of $T$.
 \end{proof}

Since each Delaunay tetrahedron contains at most eight points of $Gp$ (because $Gp\subseteq \nor(R_1)(p)$,
which has exactly eight points in each tetrahedron), this lemma gives already a global upper bound
of $8\times 15-1=119$ 
facets for Dirichlet stereohedra of full cubic groups. Even more, full groups that 
contain eight points per tetrahedron necessarily have reflection planes. So, every full group without
reflection planes has at most 4 orbit points in each tetrahedron, which gives  a bound of 
$4\times 15-1=59$.
%
 %
But a more refined bound can be given using the following ideas:
 \begin{itemize}
 \item First, we will give bounds for each group $G$ separately, using as parameters the numbers
 $s(G)$ and $m(G)$ introduced in the previous section. We do this to identify which groups are
 worth of a closer study. For example, in case $m(G)=0$ then the  list of 15 tetrahedra given in Lemma~\ref{lemma:15tetrahedra} goes down to only 11: $T$ and the ten ``neighbors of neighbors'' $T_{ij}$ 
 of $T$.

\item Second, the following lemma allows us to do almost as if the number of ``neighbors of neighbors'' of $T$ was
only six and not ten. 
 \end{itemize}

\begin{lemma}
\label{lemma:rotation}
Let $G\in \Isom(\R^3)$ be a crystallographic group and $\rho\in G$ a rotation. Let $p\in \R^3$  and let $q$ a point of the orbit of $p$, $Gp$. Let $S$ be the orbit of $q$ by the rotation $\rho$. Then, at most two points of $S$ are neighbors of $p$ in the Voronoi diagram of $Gp$, namely, those which make the smallest dihedral angle with $p$, in both directions, when seen from the axis of $\rho$.
\end{lemma} 
 
 \begin{proof}
Let $S=S_q$ be the orbit of $q$ by the rotation $\rho$, and let $S_p$ denote the orbit of $p$.
Clearly, $(S_p\cup S_q)\subseteq Gp$ implies $\Vor_{Gp}(p)\subseteq\Vor_{S_p\cup S_q}(p)$ and, in particular, points of $S_q$ that are neighbors of
$p$ in $\Vor_{Gp}$ must also be neighbors in $\Vor_{S_p\cup S_q}$. 

We now claim that if a point $q'\in S_q$ is a neighbor of $p$ in $\Vor_{S_p\cup S_q}$ then the relative interiors of the regions $\Vor_{S_p}(p)$ and $\Vor_{S_q}(q')$  intersect 
(this is actually a variation of Lemma 1.2 in~\cite{Bochis-Santos-2006} and can be stated more generally). The reason is that if a point $x$
lies in the common facet of the Voronoi regions of $p$ and $q'$ in $\Vor_{S_p\cup S_q}$ then $p$ and $q'$ are the unique closest points to $x$ in $S_p$ and $S_q$, respectively.

With this claim the proof is easy to finish. The Voronoi diagrams of both $S_p$ and $S_q$ consist of dihedra with an edge in the axis of $\rho$ and with angle $2\pi$ divided by the order of $\rho$. The dihedron $\Vor_{S_p}(p)$ can only overlap with (at most) two such dihedra of the diagram $\Vor_{S_q}$.
 
 Perhaps the special case where $p$ and $q$ lie in the same orbit modulo $\langle\rho\rangle$ needs an extra word.
 In this case, $S=S_p=S_q$ and the Voronoi diagram $\Vor_{S_p\cup S_q}(p)$ equals the cycle of dihedra around the axis of $\rho$ mentioned above. It is then clear that in it every cell (in particular, the one of $p$) has exactly two neighbors.
 \end{proof}


In the following statement, we say that a full cubic group $G$ ``mixes colors'' if it contains transformations that send black tetrahedra to white tetrahedra, and vice versa, in the 2-coloring of the Delaunay tesselation of the lattice $I$ described in Section~\ref{sec:full-groups}.

 \begin{corollary} 
 \label{coro:firstbound} 
 Let $G$ be a full group. Let $s(G)$ be the order of the stabilizer of $T$ in $G$, 
 and let $m(G)$ be equal to 1 if $G$ mixes colors of Delaunay tetrahedra and 0 if it does not.
 Then, the number of facets of any Dirichlet stereohedron for the group $G$ is bounded above by
 \[
 (7+4m(G))s(G)+3.
 \]
 \end{corollary}
   
 \begin{proof}
  The neighbors of $p$ in the Voronoi diagram of $Gp$ fall in one of the three following cases:
 \begin{itemize}
 \item Those in $T$, which are at most $s(G)-1$.
 \item Those in one of the four neighbor tetrahedra, which are at most $4s(G)$, but are present only if $m(G)=1$. That is, they are $4s(G)m(G)$.
 \item Those in one of the ten ``neighbors of neighbors'' of $T$. In principle the bound would be 
 $10 s(G)$, but it gets reduced to $6s(G)+4$ using Lemma~\ref{lemma:rotation}.
Indeed, let $S=\{p_1,\dots,p_{s(G)}\}$ be the set of points of $Gp$ inside the base tetrahedron $T$,
with $p=p_1$. Similarly, let $S_{ij}$ be the set of points of $Gp$ inside the tetrahedron $T_{ij}$.
Then, $S_{ij}$ is obtained from $S$ by the rotation of order two or
three that sends $T$ to $T_{ij}$. In particular, Lemma~\ref{lemma:rotation} says that for 
$(i,j)$ equal to one of $(1,2)$, $(1,4)$, $(3,2)$ and $(3,4)$, $S_{ij}$ and $S_{ji}$ contain, in total,
at most $s(G)+1$ neighbors of $p$ in the Dirichlet tesselation of $Gp$: the two obtained from $p_1$
by the rotation of order three, and only one of the two obtained from each of the other 
points in $S$. That is, the total number of neighbors in the ten $T_{ij}$'s is bounded above by
$2s(G)+4(s(G)+1) = 6 s(G) + 4$.
 \end{itemize}
 Adding up, the number of neighbors is at most
  \[
 s(G)-1 + 4s(G)m(G) + 6 s(G) + 4 =   4s(G)m(G) + 7 s(G) + 3.
 \]
 \end{proof}

 Table~\ref{table:firstbound} shows the bound given in this statement for each of the 27 full groups.
 Those that contain reflections appear in parentheses. As a conclusion, the global upper bound
 derived from this statement is 47.

  \begin{table}[htb]
\begin{center}
\begin{tabular}{|c|c|c|l|}
\hline s(G) & m(G) & \text{Bound} & \text{Groups} \\
\hline 
\multirow{2}*{1} & 0 & 10 & F23 \\
\cline{2-4}
& 1 & 14 & 
\begin{minipage}{6.5 cm}
\smallskip
$F432, F\frac{2}{d}\overline{3}, F\overline{4}3c, (F\overline{4}3m) $
\smallskip
\end{minipage}
\\
\hline 
\multirow{2}*{2} & 0 & 17 & 
\begin{minipage}{6.5 cm}
\smallskip
$F4_{1}32, P23, (F\frac{2}{m}\overline{3}) $
\smallskip
\end{minipage}
\\
\cline{2-4}
& 1 & 25 & 
\begin{minipage}{7 cm}
\smallskip
$F\frac{4_{1}}{d}\overline{3}\frac{2}{n}$, $I23$, $P432$, $P\frac{2}{n}\overline{3}$,
\smallskip
\\ 
$(F\frac{4_{1}}{d}\overline{3}\frac{2}{m})$, $(F\frac{4}{m}\overline{3}\frac{2}{n})$, $(P\overline{4}3m)$, $(F\frac{4}{m}\overline{3}\frac{2}{m}) $
\smallskip
\end{minipage}
\\
\hline 
\multirow{2}*{4} & 0 & 31 & 
\begin{minipage}{6.5 cm}
\smallskip
$P4_{2}32, P\overline{4}3n, (P\frac{2}{m}\overline{3})$
\smallskip
\end{minipage}
\\
\cline{2-4}
& 1 & 47 & 
\begin{minipage}{6 cm}
\smallskip
$I432, P\frac{4}{n}\overline{3}\frac{2}{n},$ 
\smallskip
\\
$(P\frac{4}{m}\overline{3}\frac{2}{m}),(I\frac{2}{m}\overline{3}), (P\frac{4_{1}}{n}\overline{3}\frac{2}{m}), (I\overline{4}3m) $
\smallskip
\end{minipage}
\\
\hline 
\multirow{2}*{8} & 0 & 59 & 
\begin{minipage}{6.5 cm}
\smallskip
$(P\frac{4}{m}\overline{3}\frac{2}{n})$
\smallskip
\end{minipage}
\\
\cline{2-4}
& 1 & 91 & 
\begin{minipage}{6.5 cm}
\smallskip
$(I\frac{4}{m}\overline{3}\frac{2}{m})$
\smallskip
\end{minipage}
\\
\hline
\end{tabular}
\end{center}
\caption{The upper bound derived from Corollary~\ref{coro:firstbound}.
\label{table:firstbound}
}
\end{table}

 \section{Groups that contain rotations of order four}
  \label{sec:order4}

 Corollary~\ref{coro:firstbound} is ultimately based on the fact that every full group
 contains the rotations of order three or two on all the edges of the Delaunay tetrahedra (see
 the proof of Lemma~\ref{lemma:5tetrahedra}). But, of course, the groups for which the  bound
 obtained 
 is worse are precisely those that have more transformations.
 Here and in the next section we take advantage of this fact by using additional rotations present in 
  some of the full groups.

In this section we look 
at the groups which contain the rotations of order four in the two edges of $T$
with dihedral angle of $\pi/2$. We call $\rho_{13}$ and $\rho_{24}$ these two rotations
($\rho_{ij}$ is meant to have axis on the edge $v_iv_j$). Certainly, every group containing them has $m=1$, since both rotations send $T$ to one of its neighbors. Moreover, to
 locate them in Table~\ref{table:fullgroups} we just need to observe that $\rho_{13}$
 and $\rho_{24}$ act on the four labels $1$, $2$, $3$ and $4$ of points in the lattice 
 $I$ as the transpositions $2\leftrightarrow 4$ and  $1\leftrightarrow 3$, respectively.
 So, a group from Table~\ref{table:fullgroups} contains these two rotations if and only if
 the corresponding pictures contains two specific dots, namely those in the
 picure for the group $P432$ (said differently, $P432$ is the group generated by $R_1$ and the rotations  $\rho_{13}$ and $\rho_{24}$).
 
 That is, the groups we are interested in are 
\[
P432, 
I432, 
P\frac{4}{n}\overline{3}\frac{2}{n},
(P\frac{4}{m}\overline{3}\frac{2}{m})
\text{ and } 
(I\frac{4}{m}\overline{3}\frac{2}{m} ).
\]
As usual, those with reflection planes are in parentheses. 

Our goal is to prove the analogs of Lemmas~\ref{lemma:5tetrahedra} and~\ref{lemma:15tetrahedra} in
this case. That is to say, we want to specify a region guaranteed to contain $\Vor_{Gp}(p)$ and derive
from it a region guaranteed to contain all neighbors of $p$. To do this we introduce the following
setting, which will be used again in the next section:

\begin{enumerate}

\item We divide the base tetrahedron $T$ in eight smaller tetrahedra $T^A$, $T^B$, $T^C$, $T^D$, $T^E$, $T^F$, $T^G$, and $T^H$,  via the two mirror reflection planes
of $T$ and their two  bisectors. These eight tetrahedra are congruent to one another.
The labels $A$ to $H$ are given cyclically, 
in the way shown in Figure~\ref{fig:minitetrahedra}. That is, $T^A$ is incident to
vertex $v_1$ and edge $v_1v_2$, $T^B$ is incident to vertex $v_2$ and edge $v_1v_2$,
and so on. 
\begin{figure}[htb]
\begin{center}
\includegraphics[width=9cm]{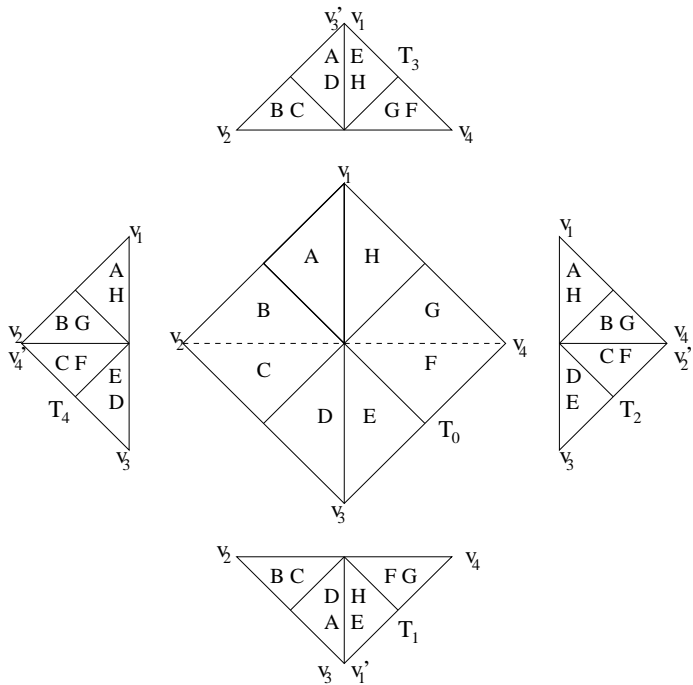}
\end{center}
\caption{The decomposition of a Delaunay tetrahedron into eight ``fundamental subdomains''.
\label{fig:minitetrahedra}}
\end{figure}
The same decomposition will be assumed in all the other Delaunay tetrahedra,
with similar notation. For example, 
$T_{12}^C$ denotes the image of $T^C$ by the unique transformation that sends $T$ to $T_{12}$
preserving the labels of vertices or, equivalently, the sub-tetrahedron in $T_{12}$
incident to its vertex $v''_{12}$ and to  the edge $v''_{12}v_3$.
Figure~\ref{fig:minitetrahedra} not only shows the decomposition of $T$, but also those of
its neighbors $T_1$, $T_2$, $T_3$ and $T_4$. The $T_i$ are drawn smaller
and using the same projection we used for $T$. Observe that this projection makes the
vertices $v'_i$ and $v_{i+2}$ coincide (indices assumed modulo four).

We call the eight tetrahedra in which each Delaunay tetrahedron is divided \emph{fundamental 
subdomain}. This is justified by the fact that they are 
fundamental domains of the normalizer $\nor(R_1)$. In particular, they are permuted
transitively by the normalizer $\nor(G)$ of any full group. 

\item Using the same arguments as in Lemma~\ref{lemma:5tetrahedra}, we will construct a region
$\VorExt(T^A)$ with the following properties:
it is a union of several fundamental subdomains and
for every $p$ in the interior of $T^A$, 
$\Vor_{Gp}(p)$ is contained in $\VorExt(T^A)$ and with no facet of $\Vor_{Gp}(p)$
in the boundary of $\VorExt(T^A)$.

 We call this region an \emph{extended Voronoi region} of $T^A$. Lemma~\ref{lemma:5tetrahedra} says that the union of $T$ and its four neighbors (which is the union of
 40 fundamental subdomains) would be a valid $\VorExt(T^A)$. But, of course, we want
 to find one as small as possible, since this will produce better bounds. 
  
 \begin{lemma}
Let $G$ be a full cubic group.
Let $S$ be a fundamental subdomain, obtained from $T^A$ by a certain transformation $g\in\nor(R_1)$.
Then, for every $q$ in the interior of $S$ we have
\[
\Vor_{Gq}(q) \subseteq g(\VorExt(T^A)).
\]
\end{lemma}

\begin{proof}
Since $g\in\nor(R_1)\subset \nor(G)$, we have that $\Vor_{Gq}(q)= g(\Vor_{Gp}(p))$, where $p=g^{-1}(q)$.
This, together with  $\Vor_{Gp}(p)\subseteq \VorExt(T^A)$,
finishes the proof.
\end{proof}

In the conditions of this lemma, we call $g(\VorExt(T^A))$ the extended Voronoi region of the subdomain $S$,
and denote it $\VorExt(S)$.

\item We call \emph{influence region} of $T^A$ the union of all the fundamental subdomains $S$
such that $\VorExt(S)$ and $\VorExt(T^A)$ overlap. We denote it $\Infl(T^A)$. We have that:

\begin{corollary}
\label{coro:influence}
 Let $p$ be a point in the interior of $T^A$ and
 let $q\in Gp$ be a point whose Dirichlet stereohedron $\Vor_{Gp}(q)$ shares a facet with 
 $\Vor_{Gp}(p)$. Then, $q$ lies in $\Infl(T^A)$.
\end{corollary}

\begin{proof}
Same as the proof of Lemma~\ref{lemma:15tetrahedra}.
\end{proof}

In particular, Lemma~\ref{lemma:15tetrahedra} implies that $\Infl(T^A)$ will be contained
in the union of $T$, its four neighbors, and the ten neighbors of the latter.

\end{enumerate}

Let us mention that this method of using fundamental subdomains, extended Voronoi regions, and influence regions was already used in~\cite{Bochis-Santos-2001} and~\cite{Bochis-Santos-2006}, with the same names,
although applied to planar crystallographic groups.

\medskip

\begin{figure}[htb]
\begin{center}
\includegraphics[width=10cm]{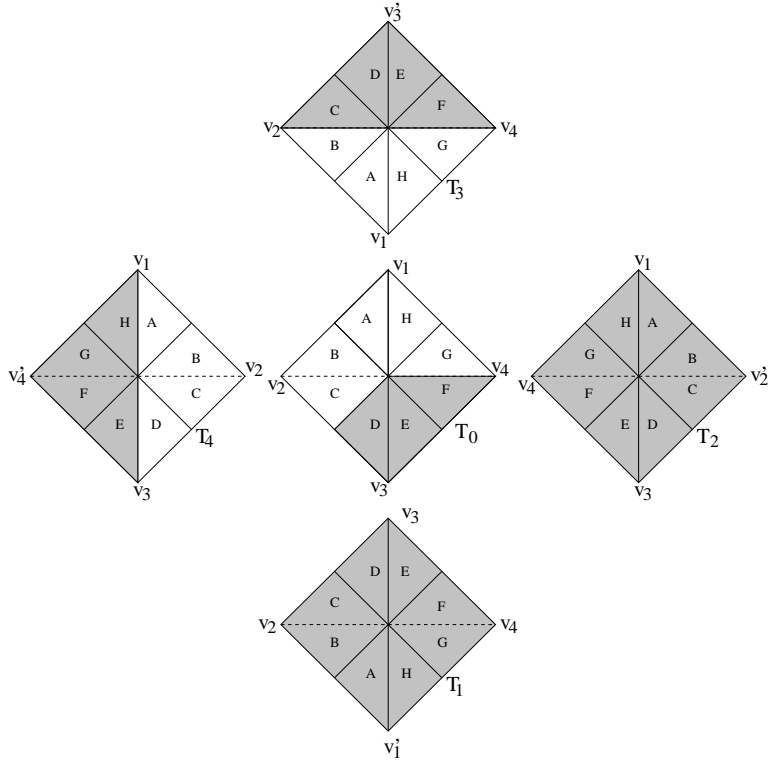}
\end{center}
\caption{The extended Voronoi region for groups that contain rotations of order four.
\label{fig:VorExt-order4}}
\end{figure}

Figure~\ref{fig:VorExt-order4} shows the Extended Voronoi region obtained for the groups
that contain the two rotations of order four. %
Our method to compute it is to start with the union of $T$ and its four neighbors and
cut some fundamental subdomains  from it via  the additional rotations present in these groups.
For example, the rotation $\rho_{13}$ allows us to exclude the whole of $T_2$ and
half of $T_4$ and the rotation $\rho_{24}$ allows us to exclude the whole of $T_1$ and
half of $T_3$. Finally, the order two rotation $\rho_0$ in the axis perpendicular to both $v_1v_3$ and
$v_2v_4$, which is an element in all these groups, allows us to exclude the fundamental subdomains $T^D$, $T^E$ and  $T^F$.
 So, the Extended Voronoi region consists of the following eleven
fundamental subdomains:
\[
\VorExt(T^A)=\{T^A, T^B, T^C, T^G, T^H, T_3^A, T_3^B, T_3^G, T_3^H, T_4^A, T_4^B, T_4^C, T_4^D  \}.
\]

\begin{figure}[htb]
\begin{center}
\includegraphics[width=12cm]{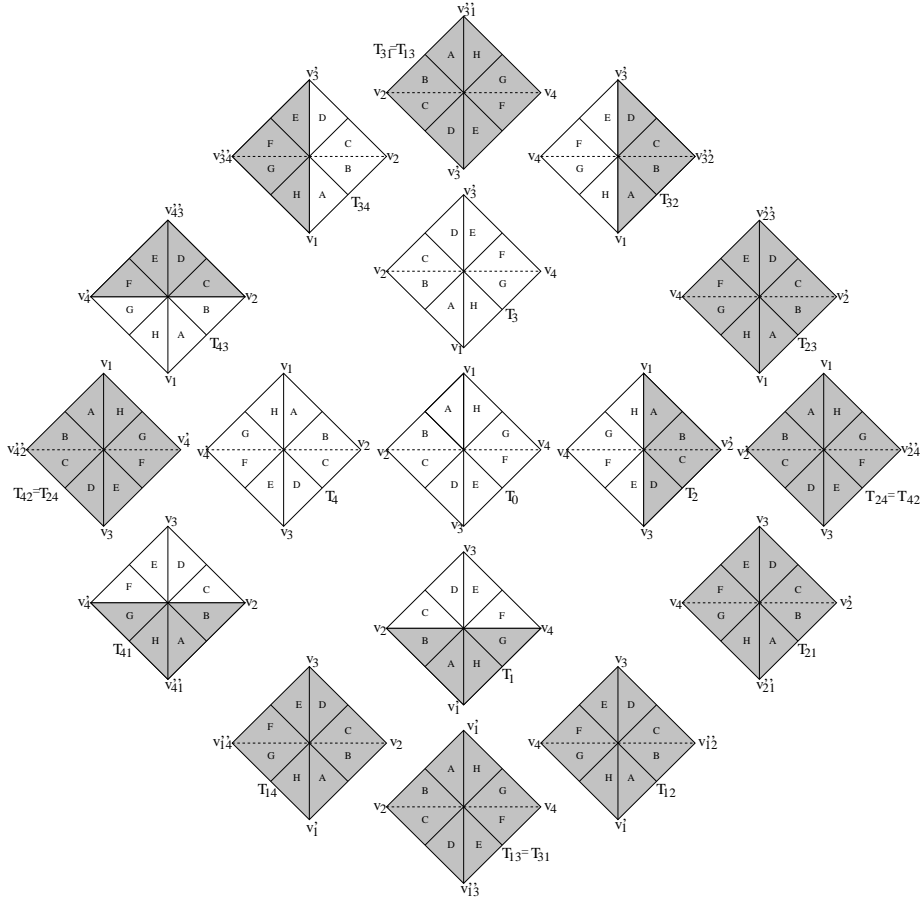}
\end{center}
\caption{The influence region for groups that contain rotations of order four.
\label{fig:influence-order4}}
\end{figure}

The influence region is shown in Figure~\ref{fig:influence-order4}. Its computation is a bit more
complicated. In principle, one should compute the 120 extended Voronoi regions for the fundamental subdomains contained in the 15 tetrahedra of Lemma~\ref{lemma:15tetrahedra} and check which ones 
overlap $\VorExt(T^A)$. But things can be slightly simplified. On the one hand, we have observed (a posteriori)
that  the influence region does not increase when adding $T^D$, $T^E$, and $T^F$ to $\VorExt(T^A)$.  
That is, we can act as if the extended Voronoi region of any given fundamental subdomain $S$ consisted of  the whole
Delaunay tetrahedron containing $S$ plus half of two of its neighbors, namely those closest to $S$. With this simplification, 
the reader should be able to check quickly that the influence region consists of the 48 fundamental subdomains
of   Figure~\ref{fig:influence-order4}, which are three full Delaunay tetrahedra plus six halves of them.

\begin{corollary}
\label{coro:order4}
Dirichlet sterohedra for groups of types 
$P432$,  $I432$ and   $P\frac{4}{n}\overline{3}\frac{2}{n}$ cannot have more than
11, 22 and 23 facets, respectively.
\end{corollary}

\begin{proof}
We start with the case of $G=P432$. If $p$ is a point in $T^A$, then $Gp$ contains 
exactly one pair of antipodal points in $T$, and in every Delaunay tetrahedron $T_i$ or
$T_{ij}$.  Since $\Infl(T^A)$
is made of 12 ``half-tetrahedra'', it contains exactly 12 points of $Gp$,
that is, eleven plus $p$ itself. 

Similarly, since $I432$ and $P\frac{4}{n}\overline{3}\frac{2}{n}$ contain \emph{two}
pairs of opposite points in each of $T$, $T_i$ or $T_{ij}$, $\Infl(T^A)$ contains
exactly 24 points of $Gp$, that is $23$ plus $p$. 
In $I432$ we subtract one, because in this case we have counted
points in both $T_{34}^B$ and $T_{43}^B$, of which only one can be 
a neighbor of $p$ by Lemma~\ref{lemma:rotation}.
\end{proof}

 \section{Groups that contain the transversal rotation of order two}
  \label{sec:order2}

We call ``transversal rotation of order two''  the symmetry $\rho_0$ of $T$ that exchanges
$v_1\leftrightarrow v_3$ and $v_2\leftrightarrow v_4$. That is, the rotation of order two
on the axis perpendicular to the edges $v_1v_3$ and $v_2v_4$, and crossing both of them.
(It is one of the transformations used in the previous section to cut out the extended Voronoi region).

Equivalently, in this section we are interested in the full groups that contain $P23$ as a subgroup. 
The list of them, excluding the ones already considered in the previous section and the ones that contain reflections, is:
\[
P23,
I23,
P\frac{2}{n}\overline{3},
P4_{2}32,
 P\overline{4}3n.
\]

\begin{figure}[htb]
\begin{center}
\includegraphics[width=10cm]{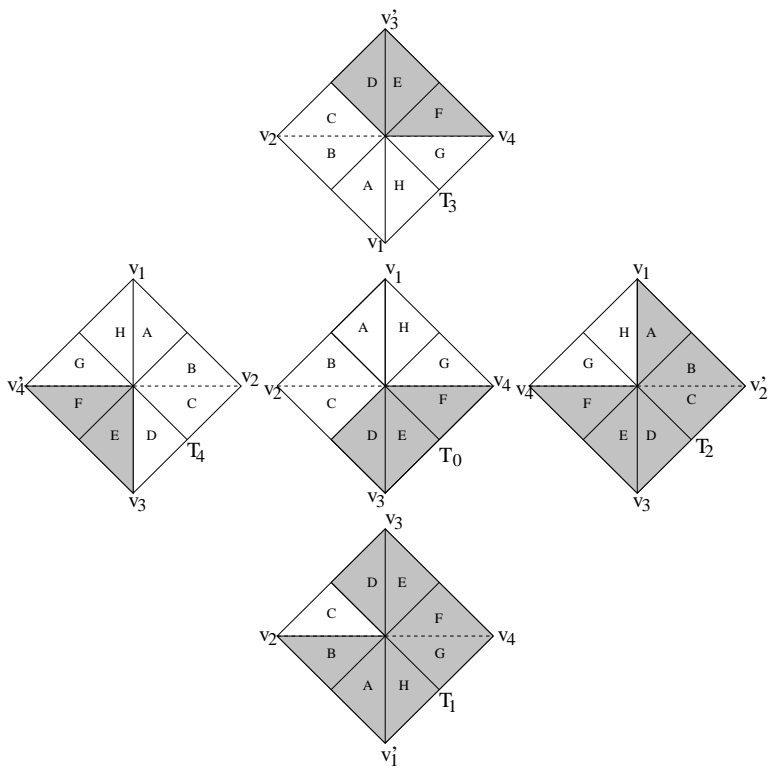}
\end{center}
\caption{The extended Voronoi region for groups that contain the transversal rotation of order two.
\label{fig:VorExt-order2}}
\end{figure}

\begin{figure}[htb]
\begin{center}
\includegraphics[width=12cm]{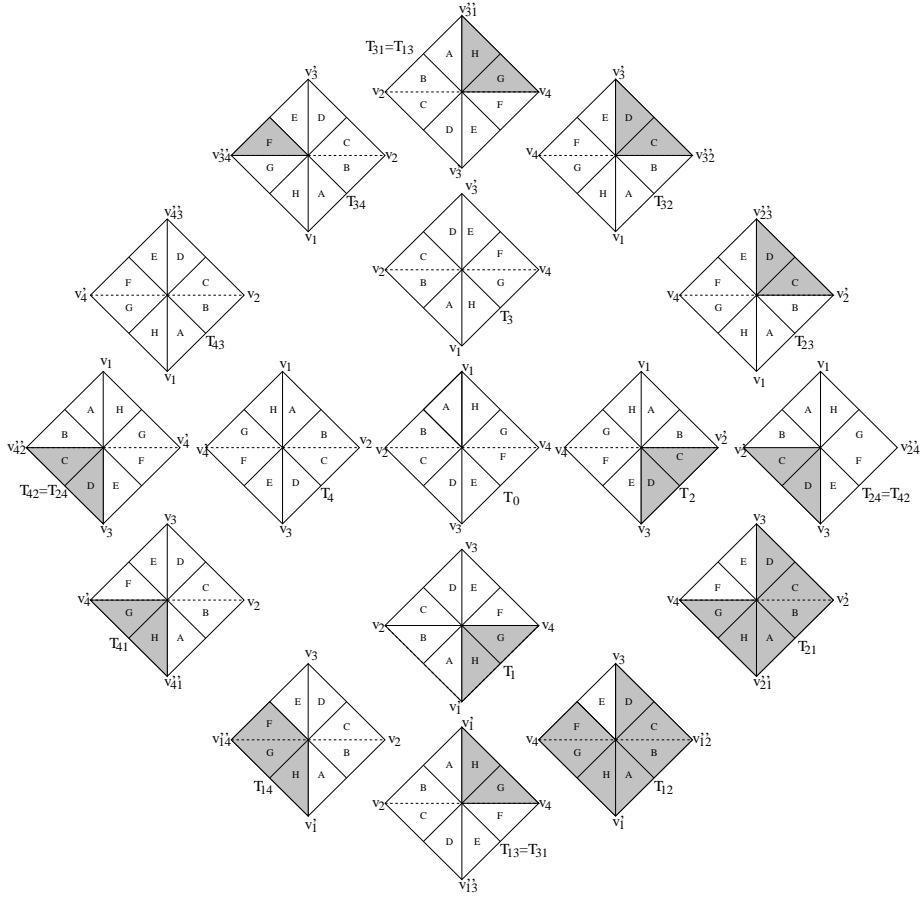}
\end{center}
\caption{The influence region for groups that contain the transversal rotation of order two.
\label{fig:influence-order2}}
\end{figure}

Figure~\ref{fig:VorExt-order2} shows the extended Voronoi region. The rotations used
to cut it from the union of $T$ and its four neighbors $T_i$ are:
\begin{itemize}
\item The transversal rotation $\rho_0$ itself, which excludes the following
fundamental subdomains:
\begin{align*}
\text{from $T$:\ }&\quad T^D, T^E, T^F, \\
\text{from $T_1$:}&\quad T_1^A, T_1^D,T_1^E,T_1^F,T_1^G,T_1^H,\\
\text{from $T_2$:}&\quad T_2^C, T_2^D, T_2^E, T_2^F.
\end{align*}
\item The rotation of order two on the axis $v_1v_3$, which excludes
the following fundamental subdomains from $T_2$:
\[
T_2^A, T_2^B,T_2^C, T_2^D.
\]
\item The rotation of order two on the axis $v_2v_4$, which excludes
the following fundamental subdomains from $T_1$:
\[
T_1^A, T_1^B,T_1^G, T_1^H.
\]
\item The rotation of  order two on the axis perpendicular to $v_1v_3$ and $v_2v'_4$;
that is, the conjugate of $\rho_0$ with axis crossing $T_4$ and $T_2$. It excludes 
the following fundamental subdomains from $T_4$:
\[
T_4^E, T_4^F.
\]
\item The rotation of  order two on the axis perpendicular to $v_1v'_3$ and $v_2v_4$;
that is, the conjugate of $\rho_0$ with axis crossing $T_1$ and $T_3$. It excludes 
the following fundamental subdomains from $T_3$:
\[
T_3^D, T_3^E, T_3^F.
\]
We leave it to the reader to verify the above assertions. There are nineteen fundamental subdomains not excluded by any of the above, namely those not shadowed in 
Figure~\ref{fig:VorExt-order2}.
\end{itemize}
We also leave out the details of the computation of the influence region.
The interested reader, should verify that sending $\VorExt(T^A)$ to be based
at each of the 31 shadowed fundamental subdomains of Figure~\ref{fig:influence-order2}
produces a region $ \VorExt(S)$ that does not intersect $\VorExt(T^A)$ itself.
Hence, the influence region of $T^A$ consists of the 89 not shadowed fundamental subdomains.

These 89 subdomains are divided as follows, according to the action of the group $R_1$.
Representatives for the sixteen orbits of fundamental subdomains via this action are,
for example, the eight subdomains contained in $T$ and the eight contained in any 
specific neighbor $T_i$. The numbers of them are given in the following tables:
\\
\begin{center}
\begin{tabular}{|ccc|}
\hline
Orbit   &      Subdomains    &  Bound on     \\
             & in $\Infl(T^A)$ & nbr. of facets  \\
\hline
$T^A$ & 9 & 8  \\
$T^B$ & 9 & 6   \\
$T^C$ & 6 & 4 \\
$T^D$ & 6 & 4  \\
$T^E$ & 11 & 7 \\
$T^F$ & 8 & 7  \\
$T^G$ & 6 & 4  \\
$T^H$ & 6 & 4 \\
\hline
Total & 61 & 44 \\
\hline
\end{tabular}
\quad
\begin{tabular}{|cc|}
\hline
Orbit   &      Subdomains      \\
             & in $\Infl(T^A)$  \\
\hline
$T_i^A$ &4  \\
$T_i^B$ & 4 \\
$T_i^C$ & 3\\
$T_i^D$ &3 \\
$T_i^E$ &4\\
$T_i^F$ & 4 \\
$T_i^G$ &3 \\
$T_i^H$ & 3\\
\hline
  Total & 28\\
\hline
\end{tabular}
\end{center}
\medskip
In the case of the subdomains in $T$, the table contains two numbers. The first number
is the number of fundamental subdomains of that orbit contained in $\Infl(T^A)$. That is,
a mere counting of the ``white triangles'' of each type in Figure~\ref{fig:influence-order2}.
The second number is the amount by which the subdomains of a given type contribute 
to the upper bound on the number of facets. It is obtained from the first number by
subtracting 1 to $T^A$ because of the base subdomain $p$, plus one to the
corresponding class by the presence of the following pairs of subdomains in $\Infl(T^A)$.
Of each pair, only one subdomain may produce a facet, by Lemma~\ref{lemma:rotation}:
{
\small
\[
\begin{array}{cccc}
 \smallskip
&
(T^B_{23},T^B_{32}) &
(T^B_{34},T^B_{43}) &
(T^B_{14},T^B_{41}) \\ \smallskip
&
&
(T^C_{34},T^C_{43}) &
(T^C_{14},T^C_{41}) \\ \smallskip
&
&
(T^D_{34},T^D_{43}) &
(T^D_{14},T^D_{41}) \\ \smallskip
(T^E_{12},T^E_{21}) &
(T^E_{23},T^E_{32}) &
(T^E_{34},T^E_{43}) &
(T^E_{14},T^E_{41}) \\ \smallskip
&
&
&
(T^F_{23},T^F_{32}) \\ \smallskip
&
(T^G_{23},T^G_{32}) &
(T^G_{34},T^G_{43}) 
& \\ 
&
(T^H_{23},T^H_{32}) &
(T^H_{34},T^H_{43}) 
&\\
\end{array}
\]
}

In order to obtain the specific bound for each of the five groups we are interested in,
we just need to add, for each of them, the bounds corresponding to the fundamental subdomains that contain a point of $Gp$, for each specific $G$.
The result is in the following table. 
\\
\begin{center}
\begin{tabular}{|c|c|c|}
\hline
Group $G$ & subdomains containing & Bound on the\\
                     &   points of $Gp$ &  number of facets\\
\hline
$P23$ & 
$T^A$, $T^E$ &
8 + 7 = 15 \\
$I23$  & 
$T^A$, $T^E$, $T_i^C$, $T_i^G$ &
8 + 7 + 3 + 3 = 21 \\
$P\frac{2}{n}\overline{3}$ & 
$T^A$, $T^E$, $T_i^B$, $T_i^F$ &
8 + 7 + 4 + 4 = 23 \\
$P4_{2}32$ & 
$T^A$, $T^E$, $T^B$, $T^F$ &
8 + 7 + 6 + 7 = 28 \\
$P\overline{4}3n$ & 
$T^A$, $T^E$, $T^C$, $T^G$ &
8 + 7 + 4 + 4 = 23 \\
\hline
\end{tabular}
\end{center}

\section{The group $P4_232$}
  \label{sec:p4232}

In this section we analyze in a bit more detail the group on which the methods so far give the worst
upper bound, namely $P4_232$. Summing up, what we have said in the previous sections is that
for every choice of a base point $p$ in the fundamental subdomain $T^A_0$:
\begin{itemize}
\item  Only the 37 orbit points lying in the following fundamental subdomains can possibly be neighbors of $p$ in the Voronoi diagram of the orbit $Gp$:
\[
\begin{array}{ccccccccccccc}
\medskip
T^B_0, &  T^E_0, &  T^F_0, \quad &
 T^A_{13}, &  T^B_{13}, &  T^E_{13}, &  T^F_{13}, &&
 T^A_{24}, &  T^B_{24}, &  T^E_{24}, &  T^F_{24}, \\
\medskip
&&&
 & &  T^E_{12}, & &&
 & &  T^E_{21}, &T^F_{21},  \\
\medskip
&&&
  T^A_{14}, &  T^B_{14}, &  T^E_{14}, &  T^F_{14}, &&
  T^A_{41}, &  T^B_{41}, &  T^E_{41}, &  T^F_{41}, \\
\medskip
&&&
  T^A_{23}, &  T^B_{23}, &  T^E_{23}, &  T^F_{23}, &&
  T^A_{32}, &  T^B_{32}, &  T^E_{32}, &  T^F_{32}, \\
\medskip
&&&
  T^A_{34}, &  T^B_{34}, &  T^E_{34}, &  & &
  T^A_{43}, &  T^B_{43}, &  T^E_{43}, &  T^F_{43}, \\
 \end{array}
\]
\item The count gets down to 28 since the two orbit points in  $T^x_{ij}$ and 
$T^x_{ji}$, with $x\in\{B,E,F\}$ and $|i-j|\ne 2$, cannot be neighbors for the same
choice of $p$ (Lemma~\ref{lemma:rotation}).
\end{itemize}

In what follows we investigate in further detail which fundamental subdomains produce 
a neighbor for all choices of $p$, or for none, and also which further dependences
are there between fundamental subdomains that prevent them to produce neighbors simultaneously.
Our main goal is to lower the upper bound a bit for this particular group, but also this analysis
gives an idea of how far our bound is from being tight.

We start with the following relatively easy observations:
\begin{lemma}
\label{lemma:p4232}
For any base point $p$ in general position in the region $T^A_0$:
\begin{enumerate}
\item The orbit points in
$T^A_0$, $T^A_{34}$, $T^A_{43}$, $T^B_{0}$, $T^B_{34}$, and $T^B_{43}$ form a Delaunay cell,
with the combinatorial type of an octahedron in which the opposite point to $T^A_0$ is 
either $T^B_{34}$ or $T^B_{43}$. 

\item Similarly, the orbit points in
$T^A_0$, $T^A_{23}$, $T^A_{32}$, $T^F_{0}$, $T^F_{23}$, and $T^F_{32}$ form a Delaunay cell,
with the combinatorial type of an octahedron in which the opposite point to $T^A_0$ is 
either $T^F_{23}$ or $T^F_{32}$. 

\item The four orbit points $T^A_{0}$, $T^B_{0}$, $T^E_{0}$, and $T^F_{0}$
 also form a Delaunay tetrahedron. 

\item Let $\pi$ be the plane that bisects $T_0$ along the mid-points of the four edges
$v_iv_{i+1}$. Let $\pi^+$ and $\pi^-$ be the half-spaces containing $v_1$ (and $v_3$) and $v_2$
(and $v_4$), respectively. If $p$ is in $\pi^+$ then $T^A_0$, $T^E_0$, $T^A_{13}$ and $T^E_{13}$
form a Delaunay tetraedron. Otherwise, $T^A_0$, $T^E_0$, $T^A_{24}$ and $T^E_{24}$ form one. That is,
for every base point either both $T^A_{13}$ and $T^E_{13}$ or both $T^A_{24}$ and $T^E_{24}$ are neighbors.
This is not exclusive, in principle the four of them could be neighbors simutaneously.
\end{enumerate}
\end{lemma}

\begin{proof}
In part one, the six points in question are equidistant to the mid-point of the segment $v_1v_2$, and
closer to this mid-point than any other orbit point. Actually, the points in $T^A_0$, $T^A_{34}$, $T^A_{43}$ form an
equilateral triangle on the side of $v_1$ and those in $T^B_{0}$, $T^B_{34}$, and $T^B_{43}$ form another, on the side of $v_2$. Depending on the direction of the plane $v_1v_2p$, the points of the second triangle closer to $p$
will be $T^B_{0}$ and $T^B_{34}$ or $T^B_0$ and $T^B_{43}$.

The same arguments, on the edge $v_1v_4$, prove part 2.

Part 3 holds since these four points are
equidistant and closest to the centroid of $T$.
For part (4), the four points in question are equidistant and closest to the midpoint of $v_1v_3$
or to the mid-point of $v_2v_4$, depending on whether 
 $p$ lies in the half-space $\pi^+$ or $\pi^-$.
\end{proof}

\begin{corollary}
\label{coro:p4232}
For every base point in general position, the Dirichlet stereohedron formed by the 
crystallographic group $P4_232$ has at least 11 facets. More precisely:
\begin{itemize}
\item $T^A_{34}$, $T^A_{43}$, $T^B_{0}$, $T^A_{23}$, $T^A_{32}$, $T^F_{0}$, and $T^E_0$
are always neighbors of $T^A_0$.
\item Exactly one of $T^B_{34}$ or $T^B_{43}$ and one of $T^F_{23}$ or $T^F_{32}$ is a neighbor of $T^A_0$.
\item At least one of $T^A_{13}$ or $T^E_{13}$ and at least one of $T^A_{24}$ or $T^E_{24}$ is a neighbor of $T^A_0$.
\end{itemize}
\end{corollary}

We have done experiments with 150 base points in general position in $T^A_0$ 
and the result is that the
number of facets has always been between 14 and 17. This shows that the lower bound
of the corollary is fairly good, while the upper bound of 28 can perhaps be lowered.
In what follows we actually lower it to 25. 

\begin{figure}[htb]
\centerline{
\includegraphics[width=13cm]{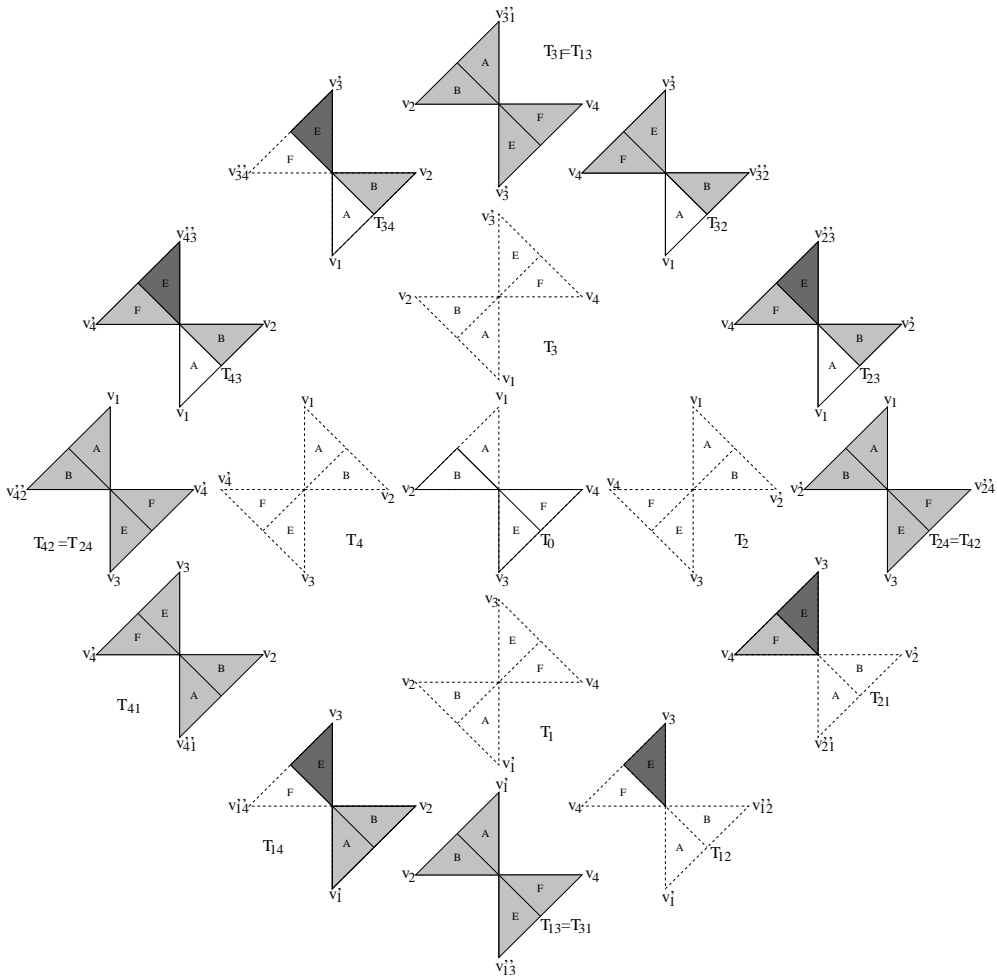}
}
\caption{Experiments with the group $P4_232$. Arbitrary base point.
\label{fig:P4_232-Completo}}
\end{figure}

\begin{figure}[htb]
\centerline{
\includegraphics[width=13cm]{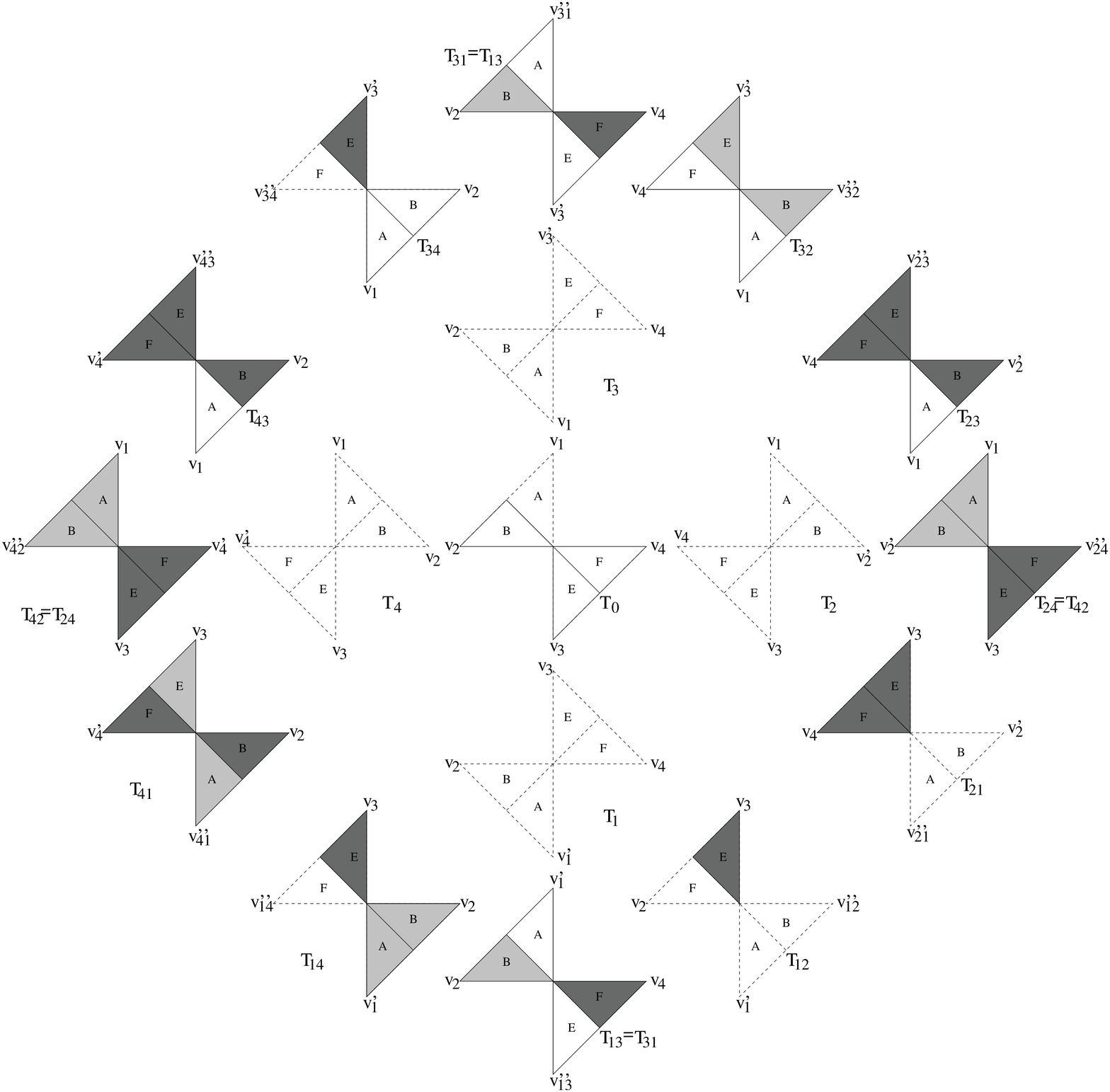}
}
\caption{Experiments with the group $P4_232$. Base point in lower half-space.
\label{fig:P4_232-Inferior}}
\end{figure}

\begin{figure}[htb]
\centerline{
\includegraphics[width=13cm]{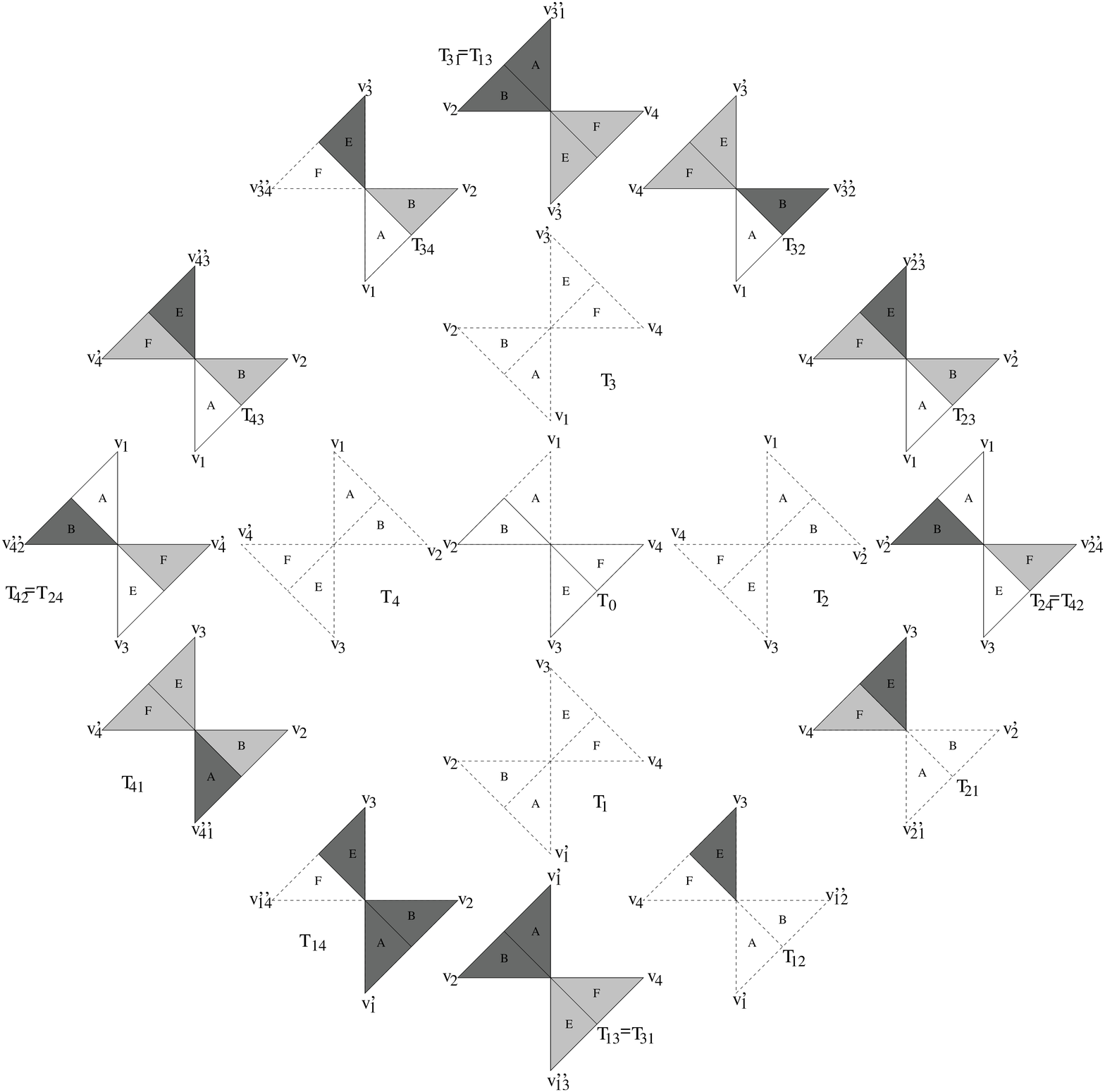}
}
\caption{Experiments with the group $P4_232$. Base point in upper half-space.
\label{fig:P4_232-Superior}}
\end{figure}

Our experimental results are shown grafically in
Figures~\ref{fig:P4_232-Completo},~\ref{fig:P4_232-Inferior} and~\ref{fig:P4_232-Superior}. In each picture, the 37 fundamental subdomains that can possibly produce neighbors is drawn. Subdomains are colored white, grey or black according to whether the observed behaviour is that
they \emph{always}, \emph{sometimes}, or \emph{never} produce a neighbor.
We do this for an arbitrary base point (Figure~\ref{fig:P4_232-Completo}) and also
separately for base points chosen in the lower and upper halfspaces $\pi^-$ and $\pi^+$
of part 4 of Lemma~\ref{lemma:p4232}.

The first observation is that the statement in Corollary~\ref{coro:p4232} is tight in the sense that all the
subdomains that are white in the pictures are mentioned in that statement. In fact, the only whites that does not explicitly follows from Lemma~\ref{lemma:p4232} is that, apparently, when
$p$ is in the lower half-space automatically $T^B_{43}$ and $T^F_{23}$ become the
opposite vertex to $T^A_0$ in the two octahedra referred to in parts 1 and 2 of the Lemma.

But more important for our purposes is to understand the black subdomains in the figures. What
we experimentally observe is:

\begin{enumerate}
\item Both in the upper and lower case, only five of the eight possible subdomains in the 
tetrahedra $T_{13}$ and $T_{24}$ produce neighbors. More precisely, if $p$ is in the lower
half-space then $T^E_{24}$, $T^F_{24}$, and $T^F_{13}$ are never neighbors, while
if $p$ is in the upper half-space then $T^A_{13}$, $T^B_{13}$, and $T^B_{24}$ are never neighbors.

\item The following six subdomains never produce a neighbor: $T^E_{34}$, 
 $T^E_{43}$,  $T^E_{12}$,  $T^E_{21}$,  $T^E_{14}$,  $T^E_{23}$.

\item There are another four and three, respectively, subdomains, that never produce a neighbor for
a base point in the lower and upper half-spaces.

\end{enumerate}

In our next and last result we prove part 1 of the above list, which decreases our
upper bound for the group $P4_2 32$ from 28 to 25. Proving part 2 would
decrease it to 23 (subtract one for each pair $T^E_{34}--T^E_{43}$,  $T^E_{12}--T^E_{21}$)
and proving part 3 would further decrease it to 21.

To lower the bound we consider the subgroup of motions of $G$ that send the 
vertical axis of $T_0$ (the line passing through the midpoints of $v_1v_3$ and $v_2v_4$) 
to itself.  This subgroup can be generated, for example,
by the three rotations of order two that stabilize $T_0$ plus those on the edges $v_1v_3$ and
$v_2v_4$ of $T_0$ (two of these five rotations are redundant, but not any two).

Let $G'$ denote this subgroup, and let $G'p$ denote the orbit of an arbitrarily chosen point in the
fundamental subdomain $T^A_0$. 

In order to simplify our arguments, we slightly translate the system of coordinates that we have been 
using so far to one in which the vertical axis of $T_0$ is actually a coordinate axis and with the origin
in the centroid of $T_0$. More precisely, we take as coordinate system the one in which
\[
v_1=(0,-1/2,1/4),v_2=(1/2,0,-1/4), v_3=(0,1/2,1/4), v_4=(-1/2,0,-1/4).
\]
The generating rotations of $G'$ then are on the following axes:
\begin{itemize}
\item The vertical coordinate axis $\{x=y=0\}$.
\item The diagonal horizontal axes at height zero $\{z=0, x+y=0\}$ and $\{z=0, x-y=0\}$
\item The horizontal axes in the coordinate directions at height $1/4$ and/or height $-1/4$.
\end{itemize}
More generally, $G'$ contains rotations of order four in the ``diagonal horizontal'' lines at height 
every integer or half-integer, and the ``coordinate horizontal'' axes at every height $n/4$ for
an odd integer $n$. It also contains the two screw rotations of order four with axis $\{x=y=0\}$ and
translation vector $(0,0,\pm1/2)$

We can also relax a bit our conditions on the base point $p$, so that for the purposes of the
following theorem we can consider $T_0$ to be the whole strip delimited by the planes
$z=\pm 1/4$, and divide it into eight sectors such as the following one that we call
$T^A_0$. 
\[
T^A_0:=\{(x,y,z)\in \R^3 : -1/4 < z < 1/4, x>0, x+y < 0\}.
\]
We also use the rest of labels $T^x_0$ or $T^x_{ij}$ in this extended sense. 

One nice feature of this subgroup $G'$, highlighted by our choice of coordinates, 
 is that the reflection that interchanges
the two half-spaces $\pi^+$ and $\pi^-$ (now the reflection on the plane $z=0$)
 is in its normalizer. In particular,  the choices of base point in the
  ``upper'' and ``lower' part of $T^A_0$ produce equivalent  Dirichlet tesselations. In what follows we assume that the base point is chosen in the upper part of $T^A_0$, but what we say is true for the lower part too except there is the
following involution that has to be applied to the names of subdomains when passing from 
an upper to a lower base point, or vice-versa: those in $T_0$ remain the same, while those
in $T_{13}$ and $T_{24}$ are exchanged as follows:
\[
T^A_{13}\leftrightarrow T^E_{24},\quad
T^B_{13}\leftrightarrow T^F_{24},\quad
T^E_{13}\leftrightarrow T^A_{24},\quad
T^F_{13}\leftrightarrow T^B_{24}.
\]
Of course, only the subdomains of types $A$, $B$, $E$ or $F$ will contain orbit points when 
the base point is in $T^A_0$.
Figure~\ref{fig:espiral} schematically shows a vertical projection of these subdomains together with (part of) an orbit for this group. The figure is designed as an aid in the proof of the following result:

\begin{figure}[htb]
\begin{center}
\includegraphics[width=7cm]{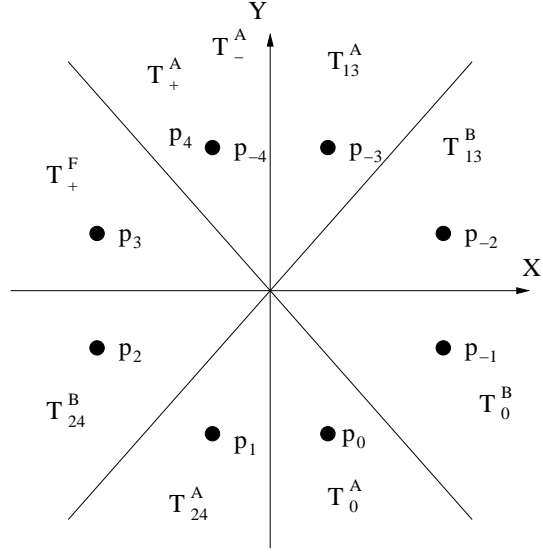}
\end{center}
\caption{(Part of) an orbit of the subgroup $G'$ of the group $G=P4_232$.
\label{fig:espiral}}
\end{figure}

\begin{theorem}
All Voronoi diagrams for a generic  orbit of the group $G'$ are equivalent, and their 
Voronoi cells have nine facets. More precisely, the neighbors of a generic base point in the
upper half of $T^A_0$ are the orbit points in the following eleven fundamental subdomains:
\[
T^E_{13}, T^F_{13}, T^ B_0, T^E_0, T^F_0,T^A_{24}, T^E_{24} , T^F_{24}, 
T^A_+, T^A_-, \text{ and } T^F_+,
\]
where 
$T^A_+= T^A_0 + (0,0,1)$, $T^A_-= T^A_0 - (0,0,1)$ and
$T^F_+= T^F_0 + (0,0,1)$ denote the translates of $T^A_0$ and $T^F_0$ one unit up or down.
\end{theorem}

Observe that we have three ``extra'' neighbors, such as  $T^A_+$, in this statement which do not appear 
anywhere else in this paper. The reason is that this is a neighbor in the Dirichlet tesselation
of the group $G'$ (not a crystallographic one) that is never a neighbor in a crystallographic full cubic group. Observe also that the notation $T^A_+$ is a bit misleading, since this fundamental subdomain should be considered ``ot type $E$'' rather than $A$. Indeed, the screw rotation or order two with translation vector $(0,0,1)$ is an element in the odd subgroup, and gives 
$T^A_+$ as the transformed of $T^E_0$, while the 
translation $(0,0,1)$ itself is not in the odd subgroup.

\begin{proof}
Let $p=p_0=(\alpha,-\beta,h)$ be our base point, with $0<h<1/4$, $0<\alpha<\beta$.
We introduce the following names for the other relevant points in the orbit. The following list
also  shows the 
fundamental subdomain they lie in. As in the statement, the subdomain $T^x_\pm$ denotes
the subdomain $T^x_0$ translated one unit up or down:
\[
\begin{array}{llcll}
\medskip
p_{-4}=(-\alpha,\beta, h- 1)&\in T^E_{-}, && q_{-4}=(\alpha,-\beta, h- 1)&\in T^A_{-}, \\
\medskip
p_{-3}=(\alpha,\beta,-h-1/2)&\in T^A_{13}, && q_{-3}=(-\alpha,-\beta,-h-1/2)&\in T^E_{13}, \\
\medskip
p_{-2}=(\beta,\alpha, h- 1/2)&\in T^B_{13}, && q_{-2}=(-\beta,-\alpha, h- 1/2)&\in T^F_{13}, \\
\medskip
p_{-1}=(\beta,-\alpha,-h) &\in T^B_{0}, && q_{-1}=(-\beta,\alpha,-h)&\in T^F_{0}, \\
\medskip
p_{0}=(\alpha,-\beta,h)  &\in T^A_{0}, && q_{0}=(-\alpha,\beta,h)&\in T^E_{0}, \\
\medskip
p_{1}=(-\alpha,-\beta,-h+1/2)&\in T^A_{24}, && q_{1}=(\alpha,\beta,-h+1/2)&\in T^E_{24}, \\
\medskip
p_{2}=(-\beta,-\alpha, h+ 1/2)&\in T^B_{24}, && q_{2}=(\beta,\alpha, h+ 1/2)&\in T^F_{24}, \\
\medskip
p_{3}=(-\beta,\alpha,-h+1)&\in T^F_{+}, && q_{3}=(\beta,-\alpha,-h+1)&\in T^F_{+}, \\
\medskip
p_{4}=(-\alpha,\beta, h+ 1)&\in T^E_{+}, && q_{4}=(\alpha,-\beta, h+ 1)&\in T^A_{+}, \\
\end{array}
\]

That the neighbors of $p_0$ must be contained in this list follows from the fact that the
translations $(0,0,\pm1)$ lie in $G'$.

The structure of the orbit is a double helix: The points $p_i$ are cyclically ordered
on the edges of an infinite prism over a (non-regular) octagon, and each point $q_i$ 
is ``opposite'' and at the same height as the corresponding $p_i$. The pairs of
points have their heights ordered by $i$, since   $0<h<1/4$.
Figure~\ref{fig:espiral} shows (the projection of) the spiral of points $p_i$, with $p_4$ superimposed over $p_{-4}$, indicating for each point the subdomain in which it lies.

From this double helix structure we can derive some properties of  the
Dirichlet tesselation or, rather, its dual Delaunay triangulation:

\begin{itemize}
\item For every $i$, the tetrahedron $\{p_{i-1},q_{i-1}, p_{i}, q_{i}\}$ is a Delaunay tetrahedron.
Indeed, its four points are equidistant and closest to the point $(0,0,i/4)$.

\item For every $i$, the tetrahedron $\{q_{i-2},q_{i-1}, p_{i}, p_{i+1}\}$ is a Delaunay tetrahedron.
Here we work out the details for $i=0$, the other cases being symmetric to it. 

Consider the empty sphere with center at the point $o=(-c,-c,0)$, for a positive $c$. 
The orbit points closest to it for $c$ sufficiently large are clearly $p_1$ and $q_{-2}$. 
On the other hand, when $c$ approaches zero,
the closest points to $o$ are those with minimum (in absolut value) height and lying in the half-space
$x+y<0$; that is, $p_0$ and $q_{-1}$. Moreover, it is easy to show that at any intermediate $c$ the closest points to $o$  must be one of these two pairs. Hence, there must be a $c$ where these
four points are equidistant and closest to $o$, hence they form a Delaunay tetrahedron.

\item Similarly, for every $i$, the tetrahedron $\{p_{i-2},p_{i-1}, q_{i}, q_{i+1}\}$ is a Delaunay tetrahedron.

\end{itemize}

These three classes of tetrahedra do not cover completely the convex hull of the orbit
(which is an infinite prism over an octagon) but almost: indeed, for each $i$ these three tetrahedra
form a ``layer''
separate the prism into two parts, so that every other Delaunay tetrahedron must lie
in between two layers. See Figure~\ref{fig:layers}, where (the vertical projection of)
the layers with $i=0$ and $i=1$ are drawn.

\begin{figure}[htb]
\begin{center}
\includegraphics[width=4cm]{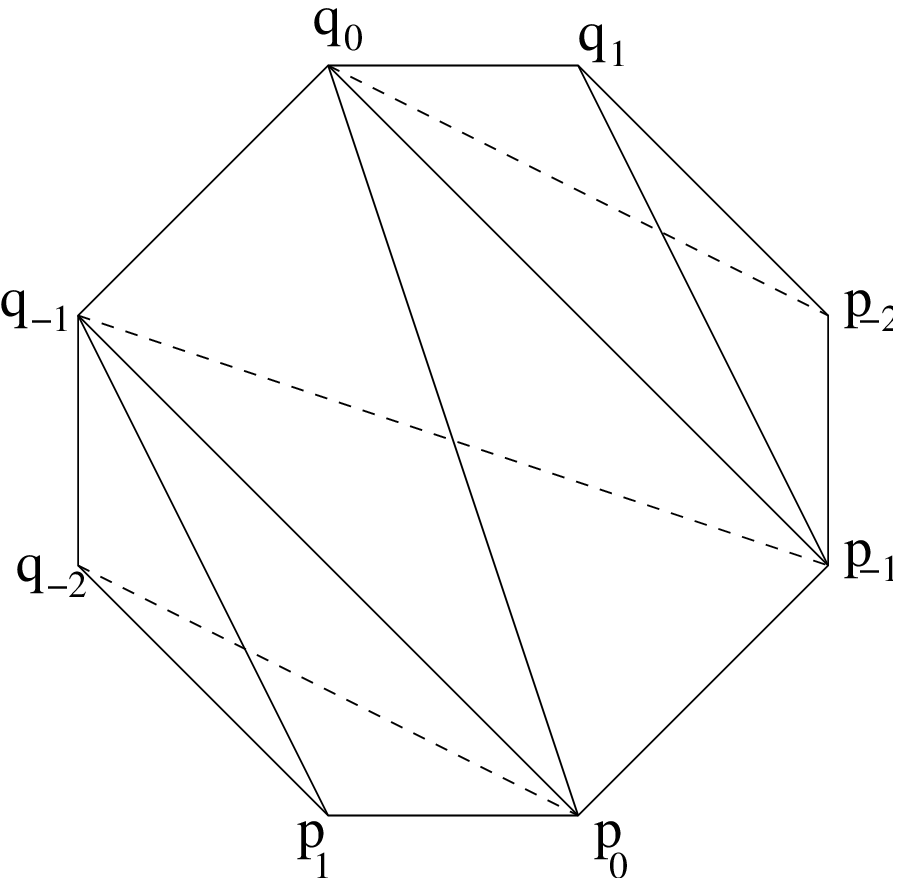}
\qquad
\includegraphics[width=4cm]{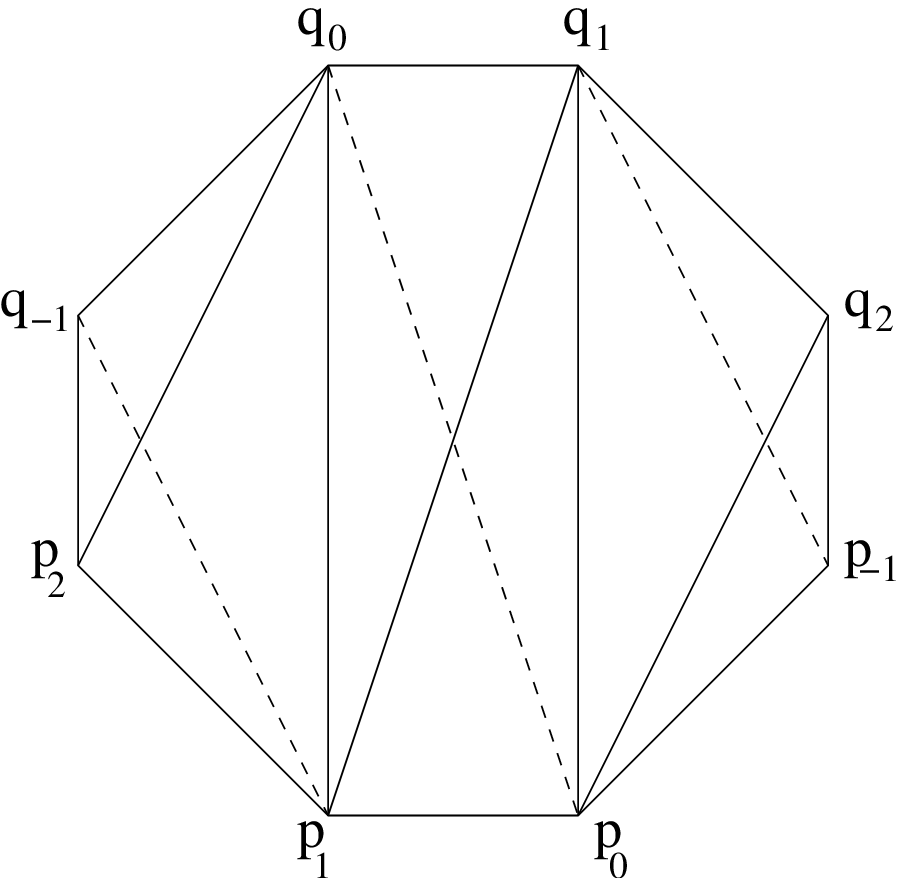}
\end{center}
\caption{Two consecutive layers (three tetrahedra each) in $Vor_{G'p_0}$. In between the two layers lie the four tetrahedra $q_{-2}q_{-1}p_1p_2$, $q_{-1}q_{0}p_0p_1$,  $p_{-1}p_{0}q_0q_1$, and 
 $p_{-2}p_{-1}q_1q_2$.
\label{fig:layers}}
\end{figure}
Moreover, the $i-1$th and the $i$th layers have three edges in common,
so that the space between them  consists of four non-overlapping and triangle-disjoint
tetrahedra. These four tetrahedra must then also be Delaunay (and add no Delaunay edges
except for $p_iq_{i\pm 4}$ ), which finishes the description
of the Delaunay tesselation. 

Hence, the list of neighbors of $p_0$ in the Delaunay tesselation is
\[q_{-4}, q_{-3}, q_{-2}, q_{-1}, q_{0}, q_{1}, q_{2}, q_{3}, q_{4}, p_{-1}, p_{1},\]
as claimed.
\end{proof}

We remark that a ``single helix'' analogue of this last statement was proved as
Lemma 4.2 in \cite[Lemma~2.8]{Bochis-Santos-2006} and, independently, in \cite[Lemma~2.1]{Erickson}.

%
%

\appendix

\section*{Appendix and erratum. Orbits with non-trivial stabilizer}

\setcounter{section}{1}
\setcounter{theorem}{0}

In this and the previous papers in this series we have always assumed that the base point $p$ for the crystallographic orbit $Gp$ under consideration has trivial stabilizer in $G$. 
That is, there is no non-trivial element $g\in G$ with $gp=p$. We believed this was no loss of generality   \emph{``since otherwise  $Gp$ is also an orbit of a proper subgroup of $G$''} (cf.~\cite[p.~92]{Bochis-Santos-2006}).
This sentence is actually false. In this appendix we show that:

\begin{itemize}
\item The sentence only fails in tetragonal and cubic groups, and only for orbits with stabilizer of order two.

\item In the tetragonal case, the sentence is almost true: every tetragonal crystallographic orbit $Gp$ is also an orbit with trivial stabilizer for a different group $G'$, except $G'$ may not be a subgroup of $G$.

\item In the cubic case, there are exactly four (one-parameter families of) cubic orbits that are not orbits with trivial stabilizer for any crystallographic group. Two arise in full groups and produce stereohedra with at most 12 facets, and two in quarter groups. For the latter, results of Koch~\cite{Koch73} give a bound of 23 facets.

\end{itemize}

\subsection{Non-trivial stabilizers have at most order two}

Let $D$ be a stereohedron for a crystallographic group $G$ (not necessarily a Dirichlet one). 
Suppose that the stabilizer of $D$ in $G$  is a non-trivial subgroup $H_0\le G$ but that, still,  $D$ is not a stereohedron for any proper subgroup of $G$. That is, there is no proper subgroup $F$ such that $\{fD:f\in F\} = \{ gD : g\in G\}$. 

Observe that $H_0$ can naturally be considered a subgroup of the point group $G_0$ of $G$:
 every element of $H_0$ fixes the centroid of $D$, which we take as the origin.
 
 \begin{lemma}
 \label{lemma:appendix-subgroup}
 There is no proper subgroup $F_0\le G_0$ such that $F_0 H_0 = G_0$.
 \end{lemma} 

Here $F_0 H_0 = \{ hf : h\in H_0, f\in F_0\}$ is not necessarily a subgroup of $G_0$. It is a union of right cosets of $F_0$ and also a union of left cosets of $H_0$. It is guaranteed to be a subgroup only if one (or both) of $F_0$ and $H_0$ are normal.

 \begin{proof}
 Let $H$ and $F$ denote the preimages of $H_0$ and $F_0$ under the canonical quotient map $G\to G_0$. Suppose that $F_0  H_0 = G_0$ and let $gD$ be a copy of $D$ obtained by some element $g\in G$. We want to show that $gD=fD$ for some $f\in F$.
 
 Let $g_0\in G_0$ be the point transformation corresponding to $g$. Since $F_0  H_0 = G_0$, there are elements $h_0\in H_0$ and $f_0 \in F_0$ such that $f_0 h_0 = g_0$. Since $H_0\le G$, we can consider $h_0$ as an element of $G$, not only a point transformation. Let $f\in F$ be any element that projects to $f_0$. Then, $fh_0 D = fD$ is a translated copy of $gD$. That is, $gD= tfD$ for some translation $t$. But $f\in F$ implies $tf\in F$, which finishes the proof.
\end{proof}

 \begin{corollary}
  \label{coro:appendix-subgroup}
$H_0$ is contained in all the subgroups of $G_0$ of index two.
 \end{corollary}
 
 \begin{proof}
If $F_0\le G_0$ has index two and $H_0$ is not contained in $F_0$ then $F_0H_0=G_0$.
 \end{proof}
 
 \begin{theorem}
   \label{thm:appendix-subgroup}
   \label{thm:stabilizer-subgroup}
 If $H_0$ is not trivial then it has order two. More precisely, either

 \begin{enumerate}
 \item $G$ is tetragonal, and $H_0$ is generated by a vertical rotation of order two; or
 
 
 \item $G$ is cubic, and $H_0$ is  generated by a rotation of order two in one of the coordinate directions.
 \end{enumerate}
  \end{theorem}

 \begin{proof}
%
(1) Let $P_{2i}$, with $i\in\{1,2,3\}$ denote the point group  $\frac{2i}{m}\frac{2}{m}\frac{2}{m}$. Put differently, $P_6$, $P_4$ and $P_2$ are the symmetry groups of a hexagonal prism, square prism, and parallelepiped, respectively. The last one is contained in the other two.

Every non-cubic point group is contained in one of these three groups. 
More precisely: triclinic, monoclinic and orthorhombic point groups are contained in $P_2$; tetragonal point groups are those contained in $P_4$ but not on $P_2$; and trigonal and hexagonal are those contained in $P_6$ but not on $P_2$. So, let $i$ be such that $G_0\le P_{2i}$.

For each index two subgroup $N$  of $P_{2i}$, $N\cap G_0$ is either $G_0$ or an index two subgroup of $G_0$. Hence, by the previous corollary, $H_0$ is contained in the intersection of the following three 
subgroups of $P_{2i}$, of index two: The group of orientation-preserving symmetries, the group of transformations sending the positive vertical direction to itself, and the group $P_i$ of  ``symmetries of the $i$-gonal prism'' ($P_1$ is the order four group generated by the mirror reflections with respect to a vertical and the horizontal plane).

In the three cases this intersection is the group of order $i$ generated by a rotation around the vertical axis. (Of course, the ``rotation of order one'' is the identity). If $i=1$ then the intersection is trivial, and the case $i=2$ is the one allowed in the statement. In the case $i=3$, if $H_0$ is not trivial (hence equal to the $C_3$ generated by a vertical rotation of order three) we can apply Lemma~\ref{lemma:appendix-subgroup} to the index three subgroup $F_0:=P_2\cap G_0$ of $G_0$. The three elements of $C_3$ belong to the three different cosets  of $P_2$ in $P_6$, hence to the three  of $F_0$ in $G_0$. 

(2) If $G$ is cubic, then $G_0$ is contained in the symmetry group of the regular cube. By the corollary, 
$H_0$ must be contained in the subgroups of orientation preserving transformations and of symmetries of the inscribed regular tetrahedron. 
The intersection of them is the group $23$ of order twelve consisting of the identity, the three rotations of order two in coordinate directions, and the eight rotations of order three in the four diagonal directions.
But $H_0$ cannot contain a rotation $\rho$ of order 3. If it did, let $F_0$ be the subgroup of $G_0$ consisting of transformations that send the $z$ axis to itself. Since $1$, $\rho$ and $\rho^2$ lie in the three cosets of $F_0$ in $G_0$, we have that $F_0H_0=G_0$.

So, $H_0$ must be contained in the group $222$ of order four generated by the rotations of order two in the three coordinate directions. It cannot equal that group either, or otherwise we could take $F_0$ the subgroup of $G_0$ that sends the diagonal axis $\{x=y=z\}$ to itself. Since the four diagonal axes are permuted transitively by the group $222$, $F_0H_0$ would equal $G_0$. 
\end{proof}

%

\subsection{Non-trivial stabilizers exist only in cubic groups}

Here we assume that $p$ is a base point with stabilizer generated by a rotation of order two in the vertical direction of a tetragonal group $G$. We prove that $Gp$ is also an orbit with trivial stabilizer
for another group $G'$. 

Since $G$ is tetragonal and $Gp$ lies in order two rotation axes, the horizontal projection of $Gp$ is a square lattice. For the sake of concreteness, let us introduce coordinates in which this lattice is the integer two-dimensional lattice $\Z^2$.

Let $G_h$ be the subgroup of $G$ generated by the vertical rotations of order two on 
axes that contain orbit points.  Let $T_h\le G_h$ be the subgroup of translations in  $G_h$. 
$T_h$ is the group of horizontal translations by vectors with even coordinates. Put differently, 
 $T_h$  splits  the rotation axes in $G_h$ into four orbits, with representatives in the vertical lines
 $l_{(0,0)}$, $l_{(1,0)}$, $l_{(1,1)}$ and $l_{(0,1)}$. ($l_{(x,y)}$ denotes the vertical line through $(x,y,0)$).

Also, let $T_v$ be the minimal vertical translation under which the orbit $Gp$ is invariant, and let $c$ be its length. On each vertical axis, $Gp$ splits into either one or two orbits of $T_v$ and it thus is an orbit of $T_v$ together with perhaps a reflection. We choose the vertical coordinate so that the orbit along the line $l_{(0,0)}$ is symmetric under the change of sign in this coordinate.

Only the three possibilities are compatible with $G$ being a tetragonal group:

\begin{itemize}
\item On all the $l_{(x,y)}$'s, the orbit points lie at the same height. Then, we let $G'$ be generated by the translations with integer horizontal coordinates and vertical coordinates multiple of $c$, together with a mirror reflection in the horizontal plane at height zero, if needed.

\item On all the $l_{(x,y)}$'s the orbit points lie at different heights. This can happen only if $G$ contains a screw rotation $\rho$ of type $4_1$ (or $4_3$) on the vertical axis through the point $(1/2,1/2,0)$.
The translation part of $\rho$ is of length $c/4$. We let $G'$ be generated by $T_h$, $T_v$, $\rho$ and a rotation of order two on the horizontal axis passing through $(0,0,0)$ and $(1,1,0)$. By our choice of the third coordinate this rotation keeps the orbit $Gp$ invariant. Hence $Gp$ is still an orbit of $G'$ and it has trivial stabilizer in it.

\item On $l_{(0,0)}$ and $l_{(1,1)}$ the orbit points are at the same height, but different to that of $l_{(1,0)}$ and $l_{(0,1)}$ (and the latter are equal to one another). If this happens, then either along each vertical axis there is only one translational orbit, or there are two but then the orbit on the axes  $l_{(0,0)}$ and $l_{(1,1)}$ is displaced a vertical length of exactly $c/2$ with respect of the orbit on $l_{(1,0)}$ and $l_{(0,1)}$ (if none of this happens, then there is no transformation of 3-space that sends $Gp$ to $Gp$ and switches the two translational orbits along the axis $l_{(0,0)}$). In the first case we let $G'$ be generated by $T_h$, $T_v$, the translation with vector $(0,1,1)$, and a rotation of order two in a horizontal axis through the points $(0,1/2,\alpha)$ and $(1,1/2,\alpha)$. Here, $\alpha$ is half the vertical displacement between the orbits in $l_{(0,0)}$ and $l_{(1,0)}$. In the second case, we let $G'$ be generated by $T_h$, $T_v$, and the translations with vectors
$(1,1,0)$ and $(0,1,c/2)$.

\end{itemize}

\subsection{Dirichlet stereohedra with non-trivial stabilizer, for cubic groups}

Fischer~\cite{Fischer73} completely classified the cubic orbits with less than three degrees of freedom. This includes the case of all cubic orbits with non-trivial stabilizer. His findings in this respect are
as follows (see also~\cite{Fischer80}):

\begin{proposition}[Fischer] There exist exactly five cubic groups with orbits that are not orbits with trivial stabilizer for non-cubic groups: the quarter groups $I\overline{4}3d$ and  $I \frac{4_1}{g} \overline{3} \frac{2}{d}$, and the full groups $F\frac{2}{d}\overline{3}$, $F 4_1 32$, and $F\frac{4}{m}\overline{3}\frac{2}{n}$.
\end{proposition}

Here we recover this result and show that, in the full case, Dirichlet stereohedra produced by these orbits have at most 12 facets.
The group $F\frac{4}{m}\overline{3}\frac{2}{n}$ is superfluous for us, since the ``bad'' orbits in it are also orbits of its 
proper subgroup $F 4_1 32$.
%


\begin{theorem}
\label{thm:stabilizer-full}
Let $p$ be a base point with non-trivial stabilizer in a full cubic group $G$, and suppose $Gp$ is not an orbit with smaller stabilizer for any other crystallographic group.
Then, $G$ is one of the groups $F4_132$ or $F\frac{2}{d}\overline{3}$ and the Dirichlet stereohedra of this orbit 
have at most twelve facets.
\end{theorem}

\begin{proof}
By Theorem~\ref{thm:appendix-subgroup}, $p$ lies in an order-two coordinate-parallel rotation axis.
There are two types of coordinate-parallel rotations of order two in full cubic groups. Those on the long sides of the fundamental subdomains and those in the common perpendicullar line to the two long sides of one fundamental subdomain.

 If $p$ lies in a rotation axis of the second type, then $Gp$ is an orbit for some proper subgroup of $G$: If there is one orbit point in each fundamental subdomain, then $Gp$ is an orbit with trivial stabilizer for either $F23$, $F432$ or $F\frac{2}{d}\overline{3}$. If there are two (symmetrically placed) orbit points in each fundamental subdomain, then $Gp$ is an orbit for either 
 $F4_1 32$, or $F\frac{4_1}{d}\overline {3}\frac{2}{n}$. See Table~2.

If $p$ lies in a long side of the fundamental subdomain, suppose that it lies on the edge $v_1v_3$ of $T$ and closer to $v_1$ than to $v_3$. We say that $p$ lies ``near $v_1$''. The whole orbit $Gp$ decomposes into at most four orbits of $F23$, one ``near $v_i$'' for each $i=1,2,3,4$.
There are four cases, depending on the action of $G$ on the four symbols $1$, $2$, $3$ and $4$ of the lattice points. That is, on how $G$ mixes the four sublattices of type $F$ that make the $I$ lattice whose Delone tetrahedra we call fundamental subdomains: 
\begin{itemize}

\item If $G$ does not mix the sublattices, then the orbit consists of six points near each center of the lattice $F$, lying in the two directions of the three coordinate axes from that center. This is also an orbit for the semi-direct product of the translational lattice $F$ and the order six point group $\overline{3}$ generated by the inversion and a diagonal rotation of order three.

\item If $G$ mixes the four sublattices, then the same is true except now for the lattice $I$.

\item If $G$ mixes $1$ only with $3$ and $2$ only with $4$, then the same holds for the lattice $P$.

\item So, the only problematic case is when $G$ mixes the sublattice $1$ with one of $2$ and $4$ and the sublattice $3$ with the other. Looking at Table 2, we see that the only full groups that do this are  $F\frac{2}{d}\overline{3}$, $F4_1 32$ and $F\frac{4}{m}\overline{3}\frac{2}{n}$. The latter is discarded since then the stabilizer of $p$ contains not  only the rotation on the axis $v_1v_3$, but also a mirror reflection.

In this case we consider $Gp$ as the union of two orbits of the index two subgroup $F23$. The Dirichlet stereohedra for each $F23$-orbit are the octahedra obtained glueing the four fundamental subdomain that share a long edge, the edge being $v_1v_3$ for one orbit and $v_2v_4$ for the other. Then, $p$ has at most eight neighbors in its own $F23$-orbit (corresponding to the eight facets of the octahedron), and at most four in the other $F23$-orbit (because each octahedron from one orbit intersects four from the other).
\end{itemize}

\end{proof}

For quarter groups,
a quick look at the list of them (Figure~2 in that paper) easily implies:

\begin{lemma}
If $p$ has non-trivial stabilizer in a quarter group $G$ and $Gp$ is not an orbit for a proper subgroup of $G$, then $G$ is either $I\overline{4}3d$ or $I\frac{4_1}{g}\overline{3}\frac{2}{d}$.
\end{lemma}

\begin{proof}
Of the other six quarter groups, three do not have rotations of order two and the other three contain the first three as subgroups of index two.
\end{proof}

\small
\bibliographystyle{abbrv}

\end{document}